\documentclass[12pt]{article}

 \usepackage{amsfonts}

 \usepackage{amsmath}
 \usepackage{amsthm}

 \usepackage{amssymb}

 \allowdisplaybreaks[4] %gong shi fen ye

 \newcommand{\me}{\mathrm{e}}

 \DeclareMathAlphabet{\mathsfsl}{OT1}{cmss}{m}{sl}

 \newtheorem{theorem}{Theorem}[section]

 \title{Carleman estimates and some inverse problems for the coupled quantitative thermoacoustic equations by boundary data. Part I: Carleman estimates}

 \author{Yunxia Shang\thanks{Mathematics \& Science college, Shanghai Normal University, 100 Guilin Road, Shanghai, 200234, P.R. China({\tt syx@shnu.edu.cn}).}
 \thanks{School of Mathematical Sciences,
 University of Science and Technology of China, 96 Jinzhai Road Baohe,
  Hefei, Anhui, 230026,  P.R. China.}\ \
 and\
 Shumin Li\thanks{CAS Wu Wen-Tsun Key Laboratory of Mathematics, University of Science and Technology of China, 96 Jinzhai Road Baohe, Hefei, Anhui, 230026, P.R. China ({\tt
 shuminli@ustc.edu.cn}).}
 }

 \date{}

\begin{document}

 \maketitle
%\begin{abstract}

\noindent
 {\bf Abstract.} In this paper, we consider Carleman estimates and inverse problems for the coupled quantitative thermoacoustic equations. In Part I, we establish Carleman estimates for the coupled quantitative thermoacoustic equations
 by assuming that the coefficients satisfy suitable conditions and taking the usual weight function $\varphi(x,t)={\rm e}^{\lambda\psi(x,t)}$, $\psi(x,t)=\left|x-x_{0}\right|^{2}-\beta\left(t-t_0\right)^{2}+\beta t_0^{2}$ for $x$ in a bounded domain in $\mathbb{R}^{n}$ with $C^{3}$-boundary and $t\in(0, T)$, where $t_0=T/2$.
 We will  discuss applications of the Carleman estimates to some inverse problems for the coupled quantitative thermoacoustic equations in
 the succeeding Part II paper \cite{part II}.
 \\
 \\
 {\bf 2010 Mathematics Subject Classfication.}  35M10, 35B45, 35R30.
%\end{abstract}

 %\begin{keywords}
 %\end{keywords}

 %\begin{AMS}
 %Mathematics Subject Classification 35L55, 35R30, 74E10
 %\end{AMS}

 \pagestyle{myheadings}
 \thispagestyle{plain}

 \section{Introduction and main results}
  \setcounter{equation}{0}

 We investigate the so called quantitative thermoacoustic tomography process (e.g.
 \cite{Akhouayri:2016,Cox:2009,Patch:2007,Stefanov:2009} and their references).
 According to \cite{Cox:2009}, assuming that the variations
in temperature and pressure are weak and neglecting the nonlinear effects, we obtain the system
\begin{equation}\label{cp}
  \left\{\begin{aligned}
  &\partial_{t}^{2}p-\rho v^{2}_{s}{\rm div}\left(\frac{1}{\rho}\nabla p\right)-\Gamma\partial_{t}\{{\rm div}\left(\kappa\nabla\theta\right)\}=\Gamma\partial_{t}\Pi_a, \\
 &\partial_{t}\theta-\frac{1}{\rho C_{p}}{\rm div}\left(\kappa\nabla\theta\right)-\frac{\theta_{0}\varsigma}{\rho C_{p}}\partial_{t}p=\frac{\Pi_a}{\rho C_{p}},\qquad \qquad\mbox{in}\ Q,
 \end{aligned}\right.
  \end{equation}
  for the temperature rise $\theta$ and the pressure perturbation $p$ from the equilibrium steady state depending on
  $(x, t)=(x_1,\cdots, x_n, t)\in Q\triangleq\Omega\times(0, T)$. Here $n\in {\mathbb{N}^{*}}$ and $\Omega$ is a bounded domain in $\mathbb{R}^{n}$ with the boundary $\partial\Omega \in C^{3}$.
 Throughout this paper, we set $\partial_j=\frac{\partial}{\partial x_j}$,
 $\partial_t=\frac{\partial}{\partial t}$,  $\partial_j^2=\frac{\partial^2}{\partial x_j^2}$, $\partial_j\partial_k=\frac{\partial^2}{\partial x_j\partial x_k}$, $\triangle=\sum_{j=1}^{n}\frac{\partial^{2}}{\partial x_{j}^{2}}$, $\partial_t^2=\frac{\partial^2}{\partial t^2}$, $1\leq j, k\leq n$.
 We assume that the mass density at steady state $\rho$, the thermal conductivity $\kappa$, the acoustic wave velocity $v_{s}$ and the isobar specific heat capacity $C_{p}$ are given strictly positive functions of  $x$ and independent of $t$, the Gr\"{u}neisen parameter $\Gamma$,  the background temperature
 $\theta_{0}$ and the volume thermal expansivity $\varsigma$ are given
 non-negative functions of  $x$ and independent of $t$, and the absorbed energy $\Pi_a$ is an unknown function of $x$ and $t$.

  We set
 \begin{equation}\label{Theta}
 \Theta(x,t)=\partial_t\theta(x,t).
 \end{equation}
 Differentiating the second equation in (\ref{cp}) with respect to $t$, we obtain
 \begin{equation}\label{equation}\left\{\begin{aligned}
 &\partial^{2}_tp-\rho_{1}(x){\rm div}\left(q(x)\nabla p\right)-\rho_{2}(x){\rm div}\left(\kappa(x)\nabla\Theta\right)=\Gamma(x)\partial_{t}\Pi_a,
  \\
 &\partial_t\Theta-\rho_{3}(x){\rm div}\left(\kappa(x)\nabla\Theta\right)-\rho_{4}(x)
  \partial^{2}_tp=\rho_3(x)\partial_{t}\Pi_a,\qquad  \mbox{in}\ Q,
 \end{aligned}\right.
 \end{equation}
  where $\rho_{1}(x)=\rho(x) v^{2}_{s}(x)$, $q(x)=\frac{1}{\rho(x)}$, $\rho_{2}(x)=\Gamma(x)$, $\rho_{3}(x)=\frac{1}{\rho(x)C_p(x)}$, and
  $\rho_{4}(x)=\frac{\theta_{0}(x)\varsigma(x)}{\rho(x)C_{p}(x)}$. Let
  $a_1(x)=\rho_{1}(x)q(x)$, $a_2(x)=\rho_{2}(x)\kappa(x)$, $a_3(x)=\rho_{3}(x)\kappa(x)$, and
  $a_4(x)=\rho_{4}(x)$.

  Therefore, in this paper, we consider Carleman estimates
  for the following strongly coupled hyperbolic-parabolic system
 \begin{equation}\label{fg}
  \left\{\begin{aligned}
  &\partial_{t}^{2}p(x,t)-a_{1}(x)\Delta p(x,t)-a_{2}(x)\Delta\Theta(x,t)=f(x,t), \\
 &\partial_{t}\Theta(x,t)-a_{3}(x)\Delta\Theta(x,t)-a_{4}(x)\partial_{t}^{2}p(x,t)=g(x,t), \quad \ \mbox{in} \quad Q,
 \end{aligned}\right.
  \end{equation}
 where $f(x,t), \ g(x,t)\in L^{2}(Q)$, $a_{j}(x)\in C^{2}(\overline{\Omega})$ $(j=1,2,3,4)$ are real-valued functions.

  In order to state our main results, we introduce more notations.
 Let $\left(x\cdot x'\right)$ denote the scalar product in $\mathbb{R}^{n}$. Let $\nu=\nu(x)=\left(\nu_1(x), \cdots,
  \nu_n(x)\right)$ denote the outward unit normal vector to $\partial\Omega$ at $x$.
  We assume that $\omega\subset\Omega$ is a  subdomain of $\Omega$ satisfy
  \begin{equation}\label{omega}
  \overline{\partial\Omega\setminus\partial\omega}\subset\big{\{}x\in\partial\Omega\big{|} \left((x-x_0)\cdot\nu(x)\right)<0\big{\}}
  \end{equation}
 with some $x_0=\left(x_0^1, x_0^2, \cdots, x_0^n\right)\in {\mathbb{R}}^n\setminus\overline{\Omega}$. Let $T>0$ be given.
 We shall establish a Carleman estimate with
 the measurement in a subboundary layer
  $$\Theta(x,t),\qquad (x,t)\in Q_\omega\triangleq\omega\times(0,T).$$

 We introduce two sets which are concerned with the coefficients $a_j(x)$, $j=1$, $2$, $3$, $4$:
 \begin{align}\label{con1}&\nonumber \mathcal{U}=\mathcal{U}_{\sigma_0, \sigma_1, M_0, M_1, M_2}=\Bigg{\{}\left.\big(a_1(x), a_2(x), a_3(x), a_4(x)\big)\in \left(C^{2}(\overline{\Omega})\right)^4\right|
 \\\nonumber&a_{1}(x)\geq \sigma_{1}, \ a_{3}(x)\geq \sigma_{1},\ a_{2}(x)\geq0, \ a_{4}(x)\geq0,\ \forall x\in\overline{\Omega},
\\\nonumber&\|a_{j}\|_{C(\overline{\Omega})}\leq M_{0}, \ \|a_{j}\|_{C^{1}(\overline{\Omega})}\leq M_{1},
\|a_{j}\|_{C^{2}(\overline{\Omega})}\leq M_{2}, \ j=1,2,3,4,
\\&3a_{1}-2\left((x-x_{0})\cdot\nabla a_{1}\right)+\frac{a_{1}}{(1+\frac{a_{2}a_{4}}{a_{3}})}\left(2(x-x_{0})\cdot\nabla\left(\frac{a_{2}a_{4}}{a_{3}}\right)\right)
\geq\sigma_{0}\Bigg{\}}
\end{align}
where the constants $M_0>1$, $M_1>0$, $M_2>0$, $\sigma_0>0$, $M_0>\sigma_1>0$ are given.
 Denote
 \begin{equation}
 m=\inf_{x\in\overline{\Omega}}\left|x-x_{0}\right|^{2}, \ \ M=\sup_{x\in\overline{\Omega}}\left|x-x_{0}\right|^{2},\ \ \mbox{and}  \ \
 {\mathcal{D}}=\sqrt{M-m}.
 \end{equation}
 We assume that $\big(a_1(x), a_2(x), a_3(x), a_4(x)\big)\in \mathcal{U}=\mathcal{U}_{\sigma_0, \sigma_1, M_0, M_1, M_2}$.
 Denote
 \begin{align}\label{alpha}
 \nonumber \alpha_{1}=&228nM^{3}_{0}M_{1}M^{\frac{1}{2}}+20M^{4}_{0}-8\sigma_1^2,
 \quad \alpha_{2}=\left(132nM_{1}M^{\frac{1}{2}}+16M_{0}\right)M^{3}_{0},
 \\\nonumber \alpha_{3}=&\min\left\{\frac{\sigma_{0}\sigma_{1}}{\left(\frac{2M^{4}_{0}\alpha_{1}}{\sigma^{2}_{1}}+\alpha_{2}\right)},
 \frac{\sigma_{1}}{8M},
 \frac{\sigma^{4}_{1}}{16M^{8}_{0}}\right\},
 \\\nonumber
 \alpha_{4}=&6M_0^2\alpha_3+\frac{2M_0^2}{\sigma_1},\ \alpha_{5}=2M_0^2\alpha_3
 +2\left(1+\frac{M^{2}_{0}}{\sigma_{1}}+\frac{3M^{3}_{0}}{\sigma^{2}_{1}}\right)M_{1}+\frac{2M_0^2}{\sigma_{1}},
 \\\nonumber
 \alpha_{6}=&\frac{4M_0^3\alpha_3}{\sigma_{1}}+2\left(1+\frac{M^{2}_{0}}{\sigma_{1}}
 +\frac{3M^{3}_{0}}{\sigma^{2}_{1}}\right)nM_{1},
 \\\nonumber \alpha_{7}=&\frac{16M^{8}_{0}}{\sigma^{7}_{1}}\left(1+\frac{M^{2}_{0}}{\sigma_{1}}\right)^{2},
 \
 \alpha_{8}=2\left(1+\frac{M^{2}_{0}}{\sigma_{1}}\right)^{2},
 \\ \alpha_{9}=&\frac{-{\mathcal{D}} \alpha_{6}+\sqrt{{\mathcal{D}}^2 \alpha_{6}^2
 +\sigma_{1}\left({\mathcal{D}}^2 \alpha_{7}+\alpha_{8}\right)}}{4\left({\mathcal{D}}^2 \alpha_{7}+\alpha_{8}\right)}.
 \end{align}
 We choose $\beta>0$ such that
 \begin{align}\label{beta}
 0<\beta<\min\left\{\alpha_{9}^{2}, \ \frac{\alpha^{2}_{1}\alpha^{2}_{3}}{16\alpha_{4}^2{\mathcal{D}}^{2}},\
 \frac{\sigma^{2}_{0}\sigma^{2}_{1}}{16\alpha_{5}^2{\mathcal{D}}^{2}},\ \frac{m^{2}\sigma^{3}_{1}}{2M_{0}\left(\sigma_{1}+M^{2}_{0}\right){\mathcal{D}}^{2}}\right\}.
 \end{align}

 Denote $t_0=\frac{T}{2}$. We will prove a Carleman estimate for (\ref{fg}) with the exponential weight function ${\rm e}^{2s\varphi}$ where
 \begin{equation}\label{0varphi}
 \varphi(x,t)={\rm e}^{\lambda\psi(x,t)}, \ \ \psi(x,t)=\left|x-x_{0}\right|^{2}-\beta\left(t-t_0\right)^{2}+\beta t_0^{2}, \ \ \forall(x,t)\in Q,
 \end{equation}
 and $\lambda>0$ is a suitably large constant.

 We set
 $\nabla=(\partial_1, \cdots, \partial_n)$, $\nabla_{x, t}=(\partial_1, \cdots, \partial_n, \partial_t),$
   $|\nabla{w}|^2=\sum_{k=1}^n{|\partial_kw|^2},$
   $|\nabla_{x, t}{w}|^2=|\nabla{w}|^2+|\partial_tw|^2,$
 and so on. $L^{2}(Q)$, $H^{2}(\Omega)$, etc. denote usual Sobolev spaces.
 We  further set
\begin{align*}
&\mathcal{H}^{2,1}(Q)=\left\{u\in L^{2}(Q); \partial_{j}u, \partial^{2}_{j}u, \partial_{j}\partial_{k}u, \partial_{t}u\in L^{2}(Q), \ j, k=1, \cdots,n \right\},
\\&\mathcal{H}^{2,2}(Q)=\left\{u\in L^{2}(Q); \partial_{j}u, \partial^{2}_{j}u, \partial_{j}\partial_{k}u, \partial_{t}u, \partial_{t}^{2}u\in L^{2}(Q), \ j, k=1, \cdots,n \right\},
\end{align*}
and $\mathcal{W}=\mathcal{H}^{2,1}(Q)\times \mathcal{H}^{2,2}(Q)$.

 \begin{theorem}\label{carleman}
 Let $(\Theta,p)\in\mathcal{W}$ satisfy (\ref{fg}) and
 \begin{equation}\label{Db2}
 \Theta(x,t)=0,\quad p(x,t)=0,\qquad (x,t)\in\Sigma\triangleq\partial\Omega\times(0,T)
 \end{equation}
 \begin{equation}\label{Tb2}
 \Theta(x,0)=\Theta(x,T)=0, \  \partial_t^j{p}(x,0)=\partial_t^j{p}(x,T)=0,\  x\in\overline{\Omega}, \ j=0,1.
 \end{equation}
 We assume that $\big(a_1, a_2, a_3, a_4\big)\in \mathcal{U}$, and that (\ref{beta}) holds.
 Then there exists a constant $\eta(\beta)>0$ such that for any
 $T\in\left(0,\frac{2({\mathcal{D}}+\eta)}{\sqrt{\beta}}\right)$, there exists a constant $\lambda_{0}>0$ such that for all $\lambda>\lambda_{0}$, there exist constants $s_{0}(\lambda)>0$ and $K_1=K_1(s_{0}$,
     $\lambda_{0}$, $\beta$, $\Omega$, $T$, $m$, $M$, $M_0$, $M_1$, $M_{2}$, $\sigma_0$, $\sigma_1)>0$  such that
 \begin{align}\label{k}
 &\nonumber\int_{Q}\left(s^{3}\lambda^{4}\varphi^{3}\Theta^{2}+\frac{1}{s\varphi}|\partial_{t}\Theta|^{2}
 +s\lambda\varphi|\nabla\Theta|^{2}+s^{3}\lambda^{4}\varphi^{3}p^{2}+s\lambda\varphi|\nabla_{x,t}p|^{2}\right){\rm e}^{2s\varphi}{\rm d}x{\rm d}t
 \\\nonumber &\leq K_1\int_{Q}\left(f^2+g^2\right){\rm e}^{2s\varphi}{\rm d}x{\rm d}t
 \\&\hspace{0.2cm}+K_1\int_{Q_{\omega}}\bigg\{s^{3}\lambda^{4}\varphi^{3}\left({\Theta}^{2}+{p}^{2}\right)+\frac{1}{s\varphi}|\partial_{t}{\Theta}|^{2}
 +s^3\lambda^3\varphi^3\left|\partial_t{p}\right|^{2}\bigg\}{\rm e}^{2s\varphi}{\rm d}x{\rm d}t,
  \end{align}
 for all $s\geq s_{0}$.
 \end{theorem}

 Let $\Gamma_{0}\subset\partial\Omega$ satisfy
 \begin{equation}\label{epsilon}
  \overline{\partial\Omega\setminus\Gamma_{0}}\subset\big{\{}x\in\partial\Omega\big{|} \left((x-x_0)\cdot\nu(x)\right)<0\big{\}}
  \end{equation}
   Let $\frac{\partial u}{\partial\nu}$ denote the normal derivative: $\frac{\partial u}{\partial\nu}(x)=\sum_{j=1}^n\nu_j(x)\partial_ju(x)$, $x\in\partial\Omega$.
 The following Carleman estimate plays an important role in the proof of Theorem \ref{carleman}:
 \begin{theorem}\label{carleman2}
 Let $(\Theta,p)\in\mathcal{W}$ satisfy (\ref{fg}), (\ref{Db2}) and (\ref{Tb2}).
 We assume that $\big(a_1, a_2, a_3, a_4\big)\in \mathcal{U}$, and that (\ref{beta}) holds.
Then there exists a constant $\eta(\beta)>0$ such that for any
$T\in\left(0,\frac{2({\mathcal{D}}+\eta)}{\sqrt{\beta}}\right)$, there exists a constant $\lambda_{0}>0$ such that for all $\lambda>\lambda_{0}$, there exist constants $s_{0}(\lambda)>0$ and $K_2=K_2(s_{0}$,
     $\lambda_{0}$, $\beta$, $\Omega$, $T$, $m$, $M$, $M_0$, $M_1$, $M_{2}$, $\sigma_0$, $\sigma_1$, $\epsilon)>0$  such that
 \begin{align}\label{k2}
 &\nonumber\int_{Q}\left(s^{3}\lambda^{4}\varphi^{3}\Theta^{2}+\frac{1}{s\varphi}|\partial_{t}\Theta|^{2}
 +s\lambda\varphi|\nabla\Theta|^{2}+s^{3}\lambda^{4}\varphi^{3}p^{2}+s\lambda\varphi|\nabla_{x,t}p|^{2}\right){\rm e}^{2s\varphi}{\rm d}x{\rm d}t
 \\&\leq K_2\int_{Q}\left(f^2+g^2\right){\rm e}^{2s\varphi}{\rm d}x{\rm d}t
 +K_2s\lambda\int_{\Gamma_0\times(0,T)}\varphi\left(\left|\frac{\partial \Theta}{\partial\nu}\right|^{2}+\left|\frac{\partial p}{\partial\nu} \right|^{2}\right){\rm d}\sigma{\rm d}t,
 \end{align}
for all $s\geq s_{0}$.
\end{theorem}

 The Carleman estimate was first introduced by Carleman \cite{carl:1939} in 1939 for the purpose of proving uniqueness in the general Cauchy problem of first-order systems  in two dimensions with simple characteristics. After that the Carleman estimate and applications were greatly developed. We refer to Bellassoued and Yamamoto \cite{bell:2018}, Chae, Imanuvilov and Kim \cite{chae:1996}, Egorov \cite{egor:1986}, Eller and Isakov \cite{elle:2000}, Fursikov and Imanuvilov \cite{furs:2000}, H$\ddot{\rm o}$rmander \cite{horm:1963,horm:1985}, Imanuvilov \cite{Iman2002}, Imanuvilov and Yamamoto \cite{iman:2004}, Isakov \cite{isak:1986}--\cite{isak:2006},  Klibanov \cite{klib:1992}, Klibanov and Timonov \cite{klib:2004}, Lavrent'ev, Romanov and Shishatskii \cite{lavr:1986}, Li \cite{li:2015},  Romanov \cite{roma:2006}, Tataru \cite{tata:1996}, Tr$\grave{\rm e}$ves \cite{trev:1970},
 Yamamoto \cite{Yamamoto:2009}, Yuan and Yamamoto \cite{yuan:2009},
 and related references in them. In this references, Carleman estimates are proved by general methods or direct methods.

 In this paper, we concentrate on the proof of Carleman estimates for (\ref{fg}), that is, Theorem \ref{carleman} and \ref{carleman2}, by a direct method. We will  discuss applications of Theorem \ref{carleman} to several inverse problems in the domain of thermoacoustic tomography (TAT) in
 the succeeding Part II paper \cite{part II}. To our best knowledge, inverse problems in the domain of TAT are very few addressed until now because the models of the literature do not consider the heating process.

 The paper consists of three sections. In Section \ref{carproof2}, we will give a direct derivation of Theorem \ref{carleman2} by apply the argument in e.g. \cite{klib:2004,lavr:1986,roma:2006,Yamamoto:2009}. In Section \ref{carproof1}, we will derive Theorem \ref{carleman} by applying Theorem \ref{carleman2} and the argument in \cite{Iman2002}.

 \section{Proof of Theorem \ref{carleman2}}\label{carproof2}
 \setcounter{equation}{0}

 By a usual density argument, we can assume that $(\Theta,p)\in C^{\infty}(Q)\times C^{\infty}(Q)$. Let $\lambda>1$ and $s>1$. Eliminating $\partial^2_t p$ from the two equations in (\ref{fg}), we have
\begin{equation}\label{a5}
\partial_{t}\Theta-a_{5}\Delta\Theta-a_{6}\Delta p=a_{4}f+g,
\end{equation}
where
\begin{equation}\label{a_56}
a_{5}=a_{3}+a_{2}a_{4} \qquad \mbox{and} \qquad a_{6}=a_{1}a_{4}.
\end{equation}
 Denote
\begin{equation}\label{l0}
{L}_{0}\left(\Theta,p\right)=\partial_{t}\Theta-a_{5}\Delta\Theta-a_{6}\Delta p.
\end{equation}
We set $\left({\vartheta},{w}\right)={\rm e}^{s\varphi}\left(\Theta,p\right)$, $\mathcal{L}\left({\vartheta},{w}\right)={\rm e}^{s\varphi}{L}_{0}\left({\rm e}^{-s\varphi}\left({\vartheta},{w}\right)\right)$.
By (\ref{0varphi}), we have
\begin{align}
&\nonumber\partial_{t}\psi=-2\beta\left(t-t_0\right), \ \partial^{2}_{t}\psi=-2\beta,\ \partial_{j}\psi=2(x_{j}-x^{j}_{0}),\
\partial_{j}\partial_{k}\psi=2\delta_{jk},
\\\nonumber&\partial_{t}\varphi=\lambda\varphi(\partial_{t}\psi), \ \partial^{2}_{t}\varphi=\lambda\varphi(\partial^{2}_{t}\psi)+\lambda^{2}\varphi(\partial_{t}\psi)^{2},\
\partial_{j}\varphi=\lambda\varphi(\partial_{j}\psi),
\\&\partial_{j}\partial_{k}\varphi=\lambda\varphi(\partial_{j}\partial_{k}\psi)
+\lambda^{2}\varphi(\partial_{j}\psi)(\partial_{k}\psi), \ \mbox{for all} \ j,k=1,\cdots,n.
\end{align}
By direct calculations, we have
\begin{align*}
&{\rm e}^{s\varphi}\left(a_{4}f+g\right)=\mathcal{L}\left({\vartheta},{w}\right)={\rm e}^{s\varphi}\left\{\partial_{t}\left({\rm e}^{-s\varphi}{\vartheta}\right)-a_{5}\Delta\left({\rm e}^{-s\varphi}{\vartheta}\right)-a_{6}\Delta\left({\rm e}^{-s\varphi}{w}\right)\right\}
\\&=-a_{5}\Delta{\vartheta}-a_{6}\Delta{w}-s^{2}\lambda^{2}\varphi^{2}|\nabla\psi|^{2}a_{5}{\vartheta}
-s^{2}\lambda^{2}\varphi^{2}|\nabla\psi|^{2}a_{6}{w}+\partial_{t}{\vartheta}
\\&\hspace{0.5cm}+2s\lambda\varphi\left(\nabla\psi\cdot\left(a_{5}\nabla{\vartheta}+a_{6}\nabla{w}\right)\right)
+A({\vartheta},{w}), \quad \mbox{in} \ Q,
\end{align*}
where
\begin{align}\label{a}
A({\vartheta},{w})=s\lambda^{2}\varphi|\nabla\psi|^{2}\left(a_{5}{\vartheta}+a_{6}{w}\right)
+2ns\lambda\varphi\left(a_{5}{\vartheta}+a_{6}{w}\right)-s\lambda\varphi\left(\partial_{t}\psi\right){\vartheta}.
\end{align}
Then taking into consideration the orders of $(s,\lambda,\varphi)$, we divide $\mathcal{L}\left({\vartheta},{w}\right)$ as follows:
\begin{align}\label{l11}
\mathcal{L}\left({\vartheta},{w}\right)=L_{1}\left({\vartheta},{w}\right)+L_{2}\left({\vartheta},{w}\right),
\end{align}
where
\begin{align}\label{l12}
&\nonumber L_{1}\left({\vartheta},{w}\right)=-a_{5}\Delta{\vartheta}-a_{6}\Delta{w}
-s^{2}\lambda^{2}\varphi^{2}|\nabla\psi|^{2}a_{5}{\vartheta}
-s^{2}\lambda^{2}\varphi^{2}|\nabla\psi|^{2}a_{6}{w},
\\&L_{2}\left({\vartheta},{w}\right)=\partial_{t}{\vartheta}+2s\lambda\varphi
\left(\nabla\psi\cdot\left(a_{5}\nabla{\vartheta}+a_{6}\nabla{w}\right)\right)
+A({\vartheta},{w}).
\end{align}
 Therefore we have
 \begin{align}\label{L12}
 &\left\|\left(a_{4}f+g\right){\rm e}^{s\varphi}\right\|^2_{L^2(Q)}
 =\left\|\mathcal{L}\left({\vartheta},{w}\right)\right\|^{2}_{L^2(Q)}=\left\|L_{1}\left({\vartheta},{w}\right)+L_{2}\left({\vartheta},{w}\right)\right\|^{2}_{L^2(Q)}
 \nonumber\\&=\left\|L_{1}\left({\vartheta},{w}\right)\right\|^{2}_{L^2(Q)}+2\int_QL_{1}\left({\vartheta},{w}\right)
L_{2}\left({\vartheta},{w}\right){\rm d}x{\rm d}t+\left\|L_{2}\left({\vartheta},{w}\right)\right\|^{2}_{L^2(Q)}.
\end{align}
 Direct calculations yeild
\begin{align}\label{jk}
&\nonumber L_{1}\left({\vartheta},{w}\right)L_{2}\left({\vartheta},{w}\right)
=-a_{5}\left(\Delta{\vartheta}\right)\left(\partial_{t}{\vartheta}\right)
-2a_{5}s\lambda\varphi\left(\Delta{\vartheta}\right)\left(\nabla\psi\cdot\left(a_{5}\nabla{\vartheta}
+a_{6}\nabla{w}\right)\right)
\\\nonumber&\hspace{0.3cm}-a_{5}\left(\Delta{\vartheta}\right)A({\vartheta},{w})
-a_{6}\left(\Delta{w}\right)\left(\partial_{t}{\vartheta}\right)
-2s\lambda\varphi a_{6}\left(\Delta{w}\right)\left(\nabla\psi\cdot\left(a_{5}\nabla{\vartheta}
+a_{6}\nabla{w}\right)\right)
\\\nonumber&\hspace{0.3cm}-a_{6}\left(\Delta{w}\right)A({\vartheta},{w})
-s^{2}\lambda^{2}\varphi^{2}|\nabla\psi|^{2}a_{5}{\vartheta}\left(\partial_{t}{\vartheta}\right)
\\\nonumber&\hspace{0.3cm}-2s^{3}\lambda^{3}\varphi^{3}|\nabla\psi|^{2}a_{5}\left(\nabla\psi\cdot\left(a_{5}\nabla{\vartheta}
+a_{6}\nabla{w}\right)\right){\vartheta}
-s^{2}\lambda^{2}\varphi^{2}|\nabla\psi|^{2}a_{5}{\vartheta}A({\vartheta},{w})
\\\nonumber&\hspace{0.3cm}-s^{2}\lambda^{2}\varphi^{2}|\nabla\psi|^{2}a_{6}{w}\left(\partial_{t}{\vartheta}\right)
-2s^{3}\lambda^{3}\varphi^{3}|\nabla\psi|^{2}a_{6}\left(\nabla\psi\cdot\left(a_{5}\nabla{\vartheta}
+a_{6}\nabla{w}\right)\right){w}
\\&\hspace{0.3cm}-s^{2}\lambda^{2}\varphi^{2}|\nabla\psi|^{2}a_{6}{w}A({\vartheta},{w})
\triangleq\sum_{k=1}^{12}J_{k}.
\end{align}
 We integrate $J_{k}$, $k=1,2,...,12$ over $Q$. We calculate them by applying the integration by parts and apply the conditions (\ref{Db2}), (\ref{Tb2}), and $\big(a_1, a_2, a_3, a_4\big)\in \mathcal{U}$.
\begin{align}\label{j1}
 \nonumber &\int_Q J_{1}{\rm d}x{\rm d}t
 =-\int_Qa_{5}\left(\Delta{\vartheta}\right)\left(\partial_{t}{\vartheta}\right){\rm d}x{\rm d}t
 \\\nonumber&=
 \int_Q\left(\nabla a_{5}\cdot \nabla {\vartheta}\right)\left(\partial_{t}{\vartheta}\right)
 {\rm d}x{\rm d}t
 +\int_Qa_{5}\partial_{t}\left(\frac{1}{2}\left|\nabla{\vartheta}\right|^{2}\right)
 {\rm d}x{\rm d}t
  \\&=
 \int_Q\left(\nabla a_{5}\cdot \nabla {\vartheta}\right)\left(\partial_{t}{\vartheta}\right)
 {\rm d}x{\rm d}t.
 \end{align}
 \begin{align*}\label{j2}
 &\int_Q J_{2}{\rm d}x{\rm d}t
 =-2 s\lambda\int_Q\varphi\left(\Delta{\vartheta}\right)\left(\nabla\psi\cdot\left(a_{5}^2
 \nabla{\vartheta}+a_{5}a_{6}\nabla{w}\right)\right){\rm d}x{\rm d}t
 \\  &=-2 s\lambda\int_Q\varphi a_{5}^2\left(\Delta{\vartheta}\right)\left(\nabla\psi\cdot
 \nabla{\vartheta}\right){\rm d}x{\rm d}t
 -2 s\lambda\int_Q\varphi a_{5}a_{6}\left(\Delta{\vartheta}\right)\left(\nabla\psi\cdot\nabla{w}\right){\rm d}x{\rm d}t
 \\  &=-2 s\lambda\int_\Sigma \varphi a^{2}_{5}\left(\nabla\psi\cdot\nabla{\vartheta}\right)
 \left(\nabla{\vartheta}\cdot \nu\right){\rm d}\sigma{\rm d}t+2 s\lambda^{2}\int_Q \varphi a^{2}_{5}\left(\nabla\psi\cdot\nabla{\vartheta}\right)^2{\rm d}x{\rm d}t\\
 &\hspace{0.3cm}+4 s\lambda\int_Q\varphi a_{5}\left(\nabla a_{5}\cdot\nabla{\vartheta}\right)
 \left(\nabla\psi\cdot\nabla{\vartheta}\right){\rm d}x{\rm d}t
 +4 s\lambda\int_Q\varphi a^{2}_{5}\left|\nabla{\vartheta}\right|^{2}{\rm d}x{\rm d}t
 \\ &\hspace{0.3cm}+s\lambda\int_Q \varphi a^{2}_{5}\left(\nabla\psi\cdot
  \nabla\left(\left|\nabla{\vartheta}\right|^2\right)\right){\rm d}x{\rm d}t
 -2s\lambda\int_{\Sigma}\varphi a_{5}a_{6}\left(\nabla\psi\cdot\nabla{w}\right)
 \left(\nabla{\vartheta}\cdot \nu\right){\rm d}\sigma{\rm d}t
 \\&\hspace{0.3cm}+2 s\lambda^{2}\int_Q\varphi a_{5}a_{6}
 \left(\nabla\psi\cdot\nabla{\vartheta}\right)\left(\nabla\psi\cdot\nabla{w}\right){\rm d}x{\rm d}t
 \\ &\hspace{0.3cm}+2 s\lambda\int_Q \varphi\left(\nabla\left(a_{5}a_{6}\right)\cdot\nabla{\vartheta}\right)\left(\nabla\psi\cdot\nabla{w}\right)
 {\rm d}x{\rm d}t
 +4 s\lambda\int_Q \varphi a_{5}a_{6}\left(\nabla{\vartheta}\cdot\nabla{w}\right){\rm d}x{\rm d}t
 \\ &\hspace{0.3cm}+2 s\lambda\int_Q \varphi a_{5}a_{6}\sum_{j,k=1}^{n}\left(\partial_{j}\psi\right)
 \left(\partial_{k}{\vartheta}\right)\left(\partial_{j}\partial_{k}{w}\right){\rm d}x{\rm d}t.
 \end{align*}
 Calculating the fifth term by the integration by parts, we have
 \begin{align*}
   (\mbox{the fifth term})=s\lambda\int_\Sigma \varphi a^{2}_{5}\left(\nabla\psi\cdot\nu\right)
  \left|\nabla{\vartheta}\right|^2{\rm d}\sigma{\rm d}t-s\lambda^2\int_Q \varphi\left|\nabla\psi\right|^2 a^{2}_{5}
  \left|\nabla{\vartheta}\right|^2{\rm d}x{\rm d}t&
  \\ \nonumber-2s\lambda\int_Q \varphi a_{5}\left(\nabla\psi\cdot\nabla a_{5}\right)
  \left|\nabla{\vartheta}\right|^2{\rm d}x{\rm d}t-s\lambda\int_Q \varphi\left(\triangle\psi\right)a^{2}_{5}
  \left|\nabla{\vartheta}\right|^2{\rm d}x{\rm d}t.&
 \end{align*}
 We note that $\triangle\psi=2n$ on $\overline{\Omega}$. By (\ref{Db2}), we have ${\vartheta}=0$ and ${w}=0$ on $\Sigma$, so that we have
 $\nabla{\vartheta}=\frac{\partial {\vartheta}}{\partial\nu}\nu$ and $\nabla{w}=\frac{\partial {w}}{\partial\nu}\nu$ on $\Sigma$.
 Therefore we have
 \begin{align}
   \nonumber &\int_Q J_{2}{\rm d}x{\rm d}t=- s\lambda\int_\Sigma \varphi\left(\nabla\psi\cdot\nu\right)a^{2}_{5}
 \left|\frac{\partial {\vartheta}}{\partial\nu}\right|^2{\rm d}\sigma{\rm d}t
 -2s\lambda\int_{\Sigma}\varphi\left(\nabla\psi\cdot\nu\right)a_{5}a_{6}
 \frac{\partial {\vartheta}}{\partial\nu}\frac{\partial {w}}{\partial\nu}{\rm d}\sigma{\rm d}t
 \\ \nonumber &\hspace{0.3cm}+2 s\lambda^{2}\int_Q \varphi a^{2}_{5}\left(\nabla\psi\cdot\nabla{\vartheta}\right)^2{\rm d}x{\rm d}t
 +4 s\lambda\int_Q\varphi a_{5}\left(\nabla a_{5}\cdot\nabla{\vartheta}\right)
 \left(\nabla\psi\cdot\nabla{\vartheta}\right){\rm d}x{\rm d}t\\ \nonumber &\hspace{0.3cm}
 +4 s\lambda\int_Q\varphi a^{2}_{5}\left|\nabla{\vartheta}\right|^{2}{\rm d}x{\rm d}t
 -s\lambda^2\int_Q \varphi\left|\nabla\psi\right|^2 a^{2}_{5}
  \left|\nabla{\vartheta}\right|^2{\rm d}x{\rm d}t
  \\ \nonumber &\hspace{0.3cm}
 -2s\lambda\int_Q \varphi a_{5}\left(\nabla\psi\cdot\nabla a_{5}\right)
  \left|\nabla{\vartheta}\right|^2{\rm d}x{\rm d}t-2ns\lambda\int_Q \varphi a^{2}_{5}
  \left|\nabla{\vartheta}\right|^2{\rm d}x{\rm d}t
  \\ \nonumber&\hspace{0.3cm}+2 s\lambda^{2}\int_Q\varphi a_{5}a_{6}
 \left(\nabla\psi\cdot\nabla{\vartheta}\right)\left(\nabla\psi\cdot\nabla{w}\right){\rm d}x{\rm d}t
 \\ \nonumber&\hspace{0.3cm}+2 s\lambda\int_Q \varphi\left(\nabla\left(a_{5}a_{6}\right)\cdot\nabla{\vartheta}\right)\left(\nabla\psi\cdot\nabla{w}\right)
 {\rm d}x{\rm d}t
 +4 s\lambda\int_Q \varphi a_{5}a_{6}\left(\nabla{\vartheta}\cdot\nabla{w}\right){\rm d}x{\rm d}t
 \\ &\hspace{0.3cm}+2 s\lambda\int_Q \varphi a_{5}a_{6}\sum_{j,k=1}^{n}\left(\partial_{j}\psi\right)
 \left(\partial_{k}{\vartheta}\right)\left(\partial_{j}\partial_{k}{w}\right){\rm d}x{\rm d}t.
 \end{align}
\begin{align}\label{j3}
\nonumber \int_Q &J_{3}{\rm d}x{\rm d}t=-\int_Q a_{5}\left(\Delta{\vartheta}\right)A({\vartheta},{w}){\rm d}x{\rm d}t
\\\nonumber&=-s\lambda^{2}\int_Q\varphi\left|\nabla\psi\right|^{2}a^{2}_{5}\left(\Delta{\vartheta}\right){\vartheta}{\rm d}x{\rm d}t
 -s\lambda^{2}\int_Q\varphi\left|\nabla\psi\right|^{2}a_{5}a_{6}\left(\Delta{\vartheta}\right){w}{\rm d}x{\rm d}t
 \\\nonumber&\hspace{0.3cm}-2n s\lambda\int_Q\varphi a^{2}_{5}\left(\Delta{\vartheta}\right){\vartheta}{\rm d}x{\rm d}t
 -2n s\lambda\int_Q\varphi a_{5}a_{6}\left(\Delta{\vartheta}\right){w}{\rm d}x{\rm d}t
 \\\nonumber&\hspace{0.3cm}+s\lambda\int_Q\varphi\left(\partial_{t}\psi\right)a_{5}\left(\Delta{\vartheta}\right){\vartheta}{\rm d}x{\rm d}t
 \\\nonumber&=s\lambda^{3}\int_Q\varphi\left|\nabla\psi\right|^{2}a^{2}_{5}
 \left(\nabla\psi\cdot\nabla{\vartheta}\right){\vartheta}{\rm d}x{\rm d}t
 +s\lambda^{2}\int_Q\varphi\left(\nabla\left(a^{2}_{5}\left|\nabla\psi\right|^{2}\right)\cdot
 \nabla{\vartheta}\right){\vartheta}{\rm d}x{\rm d}t
 \\\nonumber&\hspace{0.3cm}+s\lambda^{2}\int_Q\varphi\left|\nabla\psi\right|^{2}a^{2}_{5}
 \left|\nabla{\vartheta}\right|^{2}{\rm d}x{\rm d}t
 +s\lambda^{3}\int_Q \varphi \left|\nabla\psi\right|^{2}a_{5}a_{6}
 \left(\nabla\psi\cdot\nabla{\vartheta}\right){w}{\rm d}x{\rm d}t
 \\\nonumber&\hspace{0.3cm}+s\lambda^{2}\int_Q\varphi\left(\nabla\left(\left|\nabla\psi\right|^{2}a_{5}a_{6}\right)\cdot
 \nabla{\vartheta}\right){w}{\rm d}x{\rm d}t
 \\\nonumber&\hspace{0.3cm}+s\lambda^{2}\int_Q\varphi \left|\nabla\psi\right|^{2}a_{5}a_{6}\left(\nabla{\vartheta}\cdot\nabla{w}\right){\rm d}x{\rm d}t
 +2n s\lambda^{2}\int_Q \varphi a^{2}_{5}\left(\nabla\psi\cdot\nabla{\vartheta}\right){\vartheta}{\rm d}x{\rm d}t
 \\\nonumber&\hspace{0.3cm}+4n s\lambda\int_Q\varphi
  a_{5}\left(\nabla a_{5}\cdot\nabla{\vartheta}\right){\vartheta}{\rm d}x{\rm d}t
 +2n s\lambda\int_Q \varphi a^{2}_{5}\left|\nabla{\vartheta}\right|^{2}{\rm d}x{\rm d}t
 \\\nonumber&\hspace{0.3cm}+2n s\lambda^{2}\int_Q \varphi a_{5}a_{6}
 \left(\nabla\psi\cdot\nabla{\vartheta}\right){w}{\rm d}x{\rm d}t
 +2ns\lambda\int_Q\varphi\left(\nabla\left(a_{5}a_{6}\right)\cdot
 \nabla{\vartheta}\right){w}{\rm d}x{\rm d}t
 \\\nonumber&\hspace{0.3cm}+2n s\lambda\int_Q \varphi a_{5}a_{6}\left(\nabla{\vartheta}\cdot\nabla{w}\right){\rm d}x{\rm d}t
 -s\lambda^{2}\int_Q \varphi\left(\partial_{t}\psi\right)a_{5}
 \left(\nabla\psi\cdot\nabla{\vartheta}\right){\vartheta}{\rm d}x{\rm d}t
 \\\nonumber&\hspace{0.3cm}-s\lambda\int_Q\varphi\left(\partial_{t}\psi\right)\left(\nabla a_{5}\cdot\nabla{\vartheta}\right){\vartheta}{\rm d}x{\rm d}t
 -s\lambda\int_Q\varphi\left(\partial_{t}\psi\right)a_{5}\left|\nabla{\vartheta}\right|^{2}{\rm d}x{\rm d}t
\\\nonumber&\geq s\lambda^{2}\int_Q\varphi\left|\nabla\psi\right|^{2}a^{2}_{5}
 \left|\nabla{\vartheta}\right|^{2}{\rm d}x{\rm d}t
 +s\lambda^{2}\int_Q\varphi \left|\nabla\psi\right|^{2}a_{5}a_{6}\left(\nabla{\vartheta}\cdot\nabla{w}\right){\rm d}x{\rm d}t
 \\\nonumber&\hspace{0.3cm}
 +2n s\lambda\int_Q \varphi a^{2}_{5}\left|\nabla{\vartheta}\right|^{2}{\rm d}x{\rm d}t
 +2n s\lambda\int_Q \varphi a_{5}a_{6}\left(\nabla{\vartheta}\cdot\nabla{w}\right){\rm d}x{\rm d}t
 \\&\hspace{0.3cm}
 -s\lambda\int_Q\varphi\left(\partial_{t}\psi\right)a_{5}\left|\nabla{\vartheta}\right|^{2}{\rm d}x{\rm d}t
 -C_1s\lambda^3\int_Q\varphi\left|\nabla {\vartheta}\right|(|{\vartheta}|+|{w}|){\rm d}x{\rm d}t.
 \end{align}
 Here and henceforth, $C_k$ ($k=1, 2, \cdots$) denote generic positive constants which may depend on $n$, $\Omega$, $T$, $\beta$, $m$, $M$, $M_{j}$ $(j=0, 1, 2)$, $\sigma_{0}$, $\sigma_{1}$, $s_{0}$, $\lambda_{0}$, but are independent of $s$ and $\lambda$.
 Integrating by parts several times, we have
 \begin{align}\label{j40}
 \nonumber \int_Q & J_{4}{\rm d}x{\rm d}t=-\int_Qa_{6}\left(\Delta{w}\right)\left(\partial_{t}{\vartheta}\right){\rm d}x{\rm d}t
 \\\nonumber&=
 \int_Q\left(\partial_{t}{\vartheta}\right)\left(\nabla a_{6}\cdot\nabla{w}\right){\rm d}x{\rm d}t
 +\int_Q a_{6}\left(\nabla{w}\cdot\nabla(\partial_{t}{\vartheta})\right){\rm d}x{\rm d}t
 \\\nonumber&=
-\int_Q\left(\nabla a_{6}\cdot\nabla(\partial_{t}{w})\right){\vartheta}{\rm d}x{\rm d}t
-\int_Q a_{6}\left(\nabla(\partial_{t}{w})\cdot\nabla{\vartheta}\right){\rm d}x{\rm d}t
\\\nonumber&=\int_Q\left(\Delta a_{6}\right)\left(\partial_{t}{w}\right){\vartheta}{\rm d}x{\rm d}t
+2\int_Q \left(\nabla a_{6}\cdot\nabla{\vartheta}\right)\left(\partial_{t}{w}\right){\rm d}x{\rm d}t
\\&\hspace{0.3cm}
+\int_Q a_{6}\left(\partial_{t}{w}\right)\left(\Delta{\vartheta}\right){\rm d}x{\rm d}t.
\end{align}
We calculate the third term by applying the second equation in (\ref{fg}).
Since
\begin{align*}
&\nonumber\Delta{\vartheta}=\Delta\left({\rm e}^{s\varphi}\Theta\right)
=\left(s^{2}\lambda^{2}\varphi^{2}\left|\nabla\psi\right|^{2}+s\lambda^{2}\varphi\left|\nabla\psi\right|^{2}
+2ns\lambda\varphi\right)\Theta{\rm e}^{s\varphi}
\\\nonumber&\hspace{0.3cm}+\left(\Delta\Theta\right){\rm e}^{s\varphi}
+2s\lambda\varphi\left(\nabla\psi\cdot\nabla{\Theta}\right){\rm e}^{s\varphi},
\\&\partial_{t}{w}=\partial_{t}\left({\rm e}^{s\varphi}p\right)
=s\lambda\varphi\left(\partial_{t}\psi\right)p{\rm e}^{s\varphi}+\left(\partial_{t}p\right){\rm e}^{s\varphi},
\end{align*}
we have
\begin{align}\label{j41}
&\nonumber a_{6}\left(\partial_{t}{w}\right)\left(\Delta{\vartheta}\right)
\\\nonumber&=s\lambda\varphi\left(\partial_{t}\psi\right)
\left(s^{2}\lambda^{2}\varphi^{2}\left|\nabla\psi\right|^{2}+s\lambda^{2}\varphi\left|\nabla\psi\right|^{2}
+2ns\lambda\varphi\right)a_{6}p\Theta{\rm e}^{2s\varphi}
\\\nonumber&\hspace{0.5cm}+s\lambda\varphi\left(\partial_{t}\psi\right)a_{6}\left(\Delta\Theta\right)p{\rm e}^{2s\varphi}+2s^{2}\lambda^{2}\varphi^{2}\left(\partial_{t}\psi\right)a_{6}
\left(\nabla\psi\cdot\nabla\Theta\right)p{\rm e}^{2s\varphi}
\\\nonumber&\hspace{0.5cm}+\left(s^{2}\lambda^{2}\varphi^{2}\left|\nabla\psi\right|^{2}+s\lambda^{2}\varphi\left|\nabla\psi\right|^{2}
+2ns\lambda\varphi\right)a_{6}\Theta\left(\partial_{t}p\right){\rm e}^{2s\varphi}
\\&\hspace{0.5cm}+a_{6}\left(\partial_{t}p\right)\left(\Delta{\Theta}\right){\rm e}^{2s\varphi}
+2s\lambda\varphi a_{6}\left(\nabla\psi\cdot\nabla\Theta\right)\left(\partial_{t}p\right){\rm e}^{2s\varphi}.
\end{align}
Integrating by parts, we have
 \begin{align}\label{j42}
 &\nonumber \int_Qs\lambda\varphi\left(\partial_{t}\psi\right)a_{6}\left(\Delta\Theta\right)p
 {\rm e}^{2s\varphi}{\rm d}x{\rm d}t
 \\ \nonumber&=-s\lambda\int_Q\varphi\left(\partial_{t}\psi\right)\left(\nabla a_{6}\cdot \nabla\Theta\right)
 p{\rm e}^{2s\varphi}{\rm d}x{\rm d}t
 \\ \nonumber&\hspace{0.2cm}-s\lambda^{2}\int_Q\varphi \left(\partial_{t}\psi\right)a_{6}\left(\nabla\psi\cdot\nabla\Theta\right)p
 {\rm e}^{2s\varphi}{\rm d}x
 {\rm d}t
 -s\lambda\int_Q \varphi\left(\partial_{t}\psi\right)a_{6}\left(\nabla\Theta\cdot
 \nabla p\right){\rm e}^{2s\varphi}{\rm d}x{\rm d}t
 \\&\hspace{0.2cm}-2 s^{2}\lambda^{2}\int_Q \varphi^{2}\left(\partial_{t}\psi\right)a_{6}\left(\nabla\psi\cdot\nabla\Theta\right)p{\rm e}^{2s\varphi}{\rm d}x{\rm d}t.
 \end{align}
 By the second equation in (\ref{fg}), we have
 \begin{align*}
 a_{6}\left(\partial_{t}p\right)&\left(\Delta{\Theta}\right){\rm e}^{2s\varphi}
 =\frac{a_{6}}{a_{3}}\left(\partial_{t}p\right)\left(\partial_{t}\Theta-a_{4}\partial^{2}_{t}p-g\right)
{\rm e}^{2s\varphi}\\
&=\frac{a_{6}}{a_{3}}\left(\partial_{t}p\right)\left(\partial_{t}\Theta\right)
{\rm e}^{2s\varphi}-\frac{a_{4}a_{6}}{2a_{3}}\partial_t\left(\left|\partial_{t}p\right|^2\right)
{\rm e}^{2s\varphi}-\frac{a_{6}}{a_{3}}\left(\partial_{t}p\right)g
{\rm e}^{2s\varphi}.
\end{align*}
Therefore we have
\begin{align}\label{j43}
\nonumber \int_Q a_{6}&\left(\partial_{t}p\right)\left(\Delta{\Theta}\right){\rm e}^{2s\varphi}{\rm d}x{\rm d}t
=\int_Q \frac{a_{6}}{a_{3}}\left(\partial_{t}p\right)\left(\partial_{t}\Theta\right){\rm e}^{2s\varphi}{\rm d}x{\rm d}t
\\&+s\lambda\int_Q\varphi\left(\partial_{t}\psi\right)\frac{ a_{4}a_{6}}{a_{3}}\left|\partial_{t}p\right|^{2}
{\rm e}^{2s\varphi}{\rm d}x{\rm d}t
-\int_Q\frac{a_{6}}{a_{3}}\left(\partial_{t}p\right)g{\rm e}^{2s\varphi}{\rm d}x{\rm d}t.
\end{align}
 We note that $s>1$, $\lambda>1$ and $\varphi>1$ on $\overline{Q}$. Combinating of (\ref{j40})-(\ref{j43}) and noting that the integration of the third term in (\ref{j41}) over $Q$ plus the last term in (\ref{j42}) equals to $0$,
 we obtain
 \begin{align}\label{j4}
 \nonumber &\int_Q J_{4}{\rm d}x{\rm d}t\geq-s\lambda\int_Q \varphi\left(\partial_{t}\psi\right)a_{6}\left(\nabla\Theta\cdot
 \nabla p\right){\rm e}^{2s\varphi}{\rm d}x{\rm d}t
 \\ \nonumber&\hspace{0.2cm}+\int_Q \frac{a_{6}}{a_{3}}\left(\partial_{t}p\right)\left(\partial_{t}\Theta\right){\rm e}^{2s\varphi}{\rm d}x{\rm d}t
 +s\lambda\int_Q\varphi\left(\partial_{t}\psi\right)\frac{ a_{4}a_{6}}{a_{3}}\left|\partial_{t}p\right|^{2}
 {\rm e}^{2s\varphi}{\rm d}x{\rm d}t
  \\ \nonumber&\hspace{0.2cm}+2 s\lambda\int_Q \varphi a_{6}\left(\nabla\psi\cdot\nabla\Theta\right)\left(\partial_{t}p\right){\rm e}^{2s\varphi}{\rm d}x{\rm d}t
 \\ \nonumber&\hspace{0.2cm} -{C_2}\int_Q\left|\partial_{t}{w}\right||{\vartheta}|{\rm d}x{\rm d}t-{C_2}\int_Q\left|\partial_{t}{w}\right|\left|\nabla{\vartheta}\right|{\rm d}x{\rm d}t-{C_2}s^3\lambda^3\int_Q \varphi^3|p||\Theta|{\rm e}^{2s\varphi}{\rm d}x{\rm d}t
 \\ \nonumber&\hspace{0.2cm} -{C_2}s\lambda^{2}\int_Q\varphi |p|\left|\nabla\Theta\right|
 {\rm e}^{2s\varphi}{\rm d}x
 {\rm d}t
  -{C_2}s^{2}\lambda^{2}\int_Q\varphi^{2}\left|\partial_{t}p\right||\Theta|{\rm e}^{2s\varphi}{\rm d}x{\rm d}t
  \\&\hspace{0.2cm}
 -{C_2}\int_Q\left|\partial_{t}p\right||g|{\rm e}^{2s\varphi}{\rm d}x{\rm d}t.
 \end{align}
 Similarly to $J_2$, we have
 \begin{align*}
 &\int_Q J_{5}{\rm d}x{\rm d}t=-2 s\lambda\int_Q \varphi a_{6}\left(\Delta {w}\right)\left(\nabla\psi\cdot\left(a_{5}\nabla{\vartheta}+
 a_{6}\nabla{w}\right)\right){\rm d}x{\rm d}t
 \\&=-2 s\lambda\int_Q \varphi a_{5}a_{6}\left(\nabla\psi \cdot
 \nabla{\vartheta}\right)\left(\Delta {w}\right){\rm d}x{\rm d}t-2 s\lambda\int_Q \varphi a^{2}_{6}\left(\nabla\psi\cdot\nabla{w}\right)\left(\Delta {w}\right){\rm d}x{\rm d}t
 \\&=-2s\lambda\int_\Sigma \varphi\left(\nabla\psi\cdot \nu\right)a_{5}a_{6}
 \frac{\partial{\vartheta}}{\partial \nu}\frac{\partial{w}}{\partial \nu}{\rm d}\sigma{\rm d}t+4 s\lambda\int_Q \varphi a_{5}a_{6}\left(\nabla{\vartheta}\cdot\nabla{w}\right){\rm d}x{\rm d}t
 \\&\hspace{0.5cm}+2s\lambda^{2}\int_Q\varphi a_{5}a_{6}\left(\nabla\psi\cdot \nabla {\vartheta}\right)
 \left(\nabla\psi\cdot\nabla{w}\right){\rm d}x{\rm d}t
 \\&\hspace{0.5cm}+2s\lambda\int_Q\varphi\left(\nabla\left(a_{5}a_{6}\right)\cdot\nabla w\right)
 \left(\nabla\psi\cdot\nabla{\vartheta}\right){\rm d}x{\rm d}t
 \\&\hspace{0.5cm}+2 s\lambda\int_Q \varphi a_{5}a_{6}\sum_{j,k=1}^{n}\left(\partial_{j}\psi\right)
 \left(\partial_{j}\partial_{k}{\vartheta}\right)\left(\partial_{k}{w}\right){\rm d}x{\rm d}t
 \\&\hspace{0.5cm}-2 s\lambda\int_\Sigma \varphi\left(\nabla\psi\cdot\nu\right)a^{2}_{6}
 \left|\frac{\partial{w}}{\partial \nu}\right|^2{\rm d}\sigma{\rm d}t
 +2 s\lambda^{2}\int_Q \varphi a^{2}_{6}\left(\nabla\psi\cdot\nabla{w}\right)^2{\rm d}x{\rm d}t
 \\&\hspace{0.5cm}+4 s\lambda\int_Q \varphi a_{6}\left(\nabla a_{6}\cdot\nabla{w}\right)
 \left(\nabla\psi\cdot\nabla{w}\right){\rm d}x{\rm d}t
 +4 s\lambda\int_Q \varphi a^{2}_{6}\left|\nabla{w}\right|^{2}{\rm d}x{\rm d}t
 \\&\hspace{0.5cm}
 +s\lambda\int_Q\varphi a^{2}_{6}\left(\nabla\psi\cdot\nabla
 \left(\left|\nabla{w}\right|^2\right)\right){\rm d}x{\rm d}t.
 \end{align*}
 Calculating the fifth term by the integration by parts, we have
 \begin{align*}
 &(\mbox{the fifth term})=2 s\lambda\int_\Sigma\varphi\left(\nabla\psi\cdot \nu\right)a_{5}a_{6}
  \frac{\partial{\vartheta}}{\partial\nu}\frac{\partial{w}}{\partial \nu}{\rm d}\sigma{\rm d}t
 \\&\hspace{0.5cm}
 -2 s\lambda^{2}\int_Q\varphi\left|\nabla\psi\right|^{2}a_{5}a_{6}\left(\nabla{\vartheta}\cdot\nabla{w}\right){\rm d}x{\rm d}t
 \\&\hspace{0.5cm}-2s\lambda\int_Q\varphi\left(\nabla\left(a_{5}a_{6}\right)\cdot\nabla\psi\right)
 \left(\nabla{\vartheta}\cdot\nabla{w}\right){\rm d}x{\rm d}t
 -4n s\lambda\int_Q\varphi a_{5}a_{6}\left(\nabla{\vartheta}\cdot\nabla{w}\right){\rm d}x{\rm d}t
 \\&\hspace{0.5cm}-2 s\lambda\int_Q\varphi a_{5}a_{6}\sum_{j,k=1}^{n}\left(\partial_{j}\psi\right)
 \left(\partial_{k}{\vartheta}\right)\left(\partial_{j}\partial_{k}{w}\right){\rm d}x{\rm d}t.
 \end{align*}
 Calculating the tenth term by the integration by parts, we have
 \begin{align*}
 &(\mbox{the tenth term})=s\lambda\int_\Sigma\varphi\left(\nabla\psi\cdot\nu\right)a^{2}_{6}
 \left|\frac{\partial{w}}{\partial\nu}\right|^{2}{\rm d}\sigma{\rm d}t
 -s\lambda^{2}\int_Q \varphi\left|\nabla\psi\right|^{2}a^{2}_{6}\left|\nabla{w}\right|^{2}{\rm d}x{\rm d}t
 \\&\hspace{0.5cm}-2 s\lambda\int_Q \varphi a_{6}\left(\nabla a_{6}\cdot \nabla\psi\right)
 \left|\nabla{w}\right|^{2}{\rm d}x{\rm d}t
 -2n s\lambda\int_Q \varphi a^{2}_{6}\left|\nabla{w}\right|^{2}{\rm d}x{\rm d}t.
 \end{align*}
 Therefore we obtain
 \begin{align}\label{j5}
 \nonumber &\int_Q J_{5}{\rm d}x{\rm d}t
 =4(1-n) s\lambda\int_Q \varphi a_{5}a_{6}\left(\nabla{\vartheta}\cdot\nabla{w}\right){\rm d}x{\rm d}t
 \\ \nonumber&\hspace{0.5cm}+2s\lambda^{2}\int_Q\varphi a_{5}a_{6}\left(\nabla\psi\cdot \nabla {\vartheta}\right)
 \left(\nabla\psi\cdot\nabla{w}\right){\rm d}x{\rm d}t
 \\ \nonumber&\hspace{0.5cm}+2s\lambda\int_Q\varphi\left(\nabla\left(a_{5}a_{6}\right)\cdot\nabla w\right)
 \left(\nabla\psi\cdot\nabla{\vartheta}\right){\rm d}x{\rm d}t
 \\ \nonumber&\hspace{0.5cm}
 -2 s\lambda^{2}\int_Q\varphi\left|\nabla\psi\right|^{2}a_{5}a_{6}\left(\nabla{\vartheta}\cdot\nabla{w}\right){\rm d}x{\rm d}t
 \\ \nonumber&\hspace{0.5cm}-2s\lambda\int_Q\varphi\left(\nabla\left(a_{5}a_{6}\right)\cdot\nabla\psi\right)
 \left(\nabla{\vartheta}\cdot\nabla{w}\right){\rm d}x{\rm d}t
 \\ \nonumber&\hspace{0.5cm}-2 s\lambda\int_Q\varphi a_{5}a_{6}\sum_{j,k=1}^{n}\left(\partial_{j}\psi\right)
 \left(\partial_{k}{\vartheta}\right)\left(\partial_{j}\partial_{k}{w}\right){\rm d}x{\rm d}t \\ \nonumber&\hspace{0.5cm}-s\lambda\int_\Sigma \varphi\left(\nabla\psi\cdot\nu\right)a^{2}_{6}
 \left|\frac{\partial{w}}{\partial \nu}\right|^2{\rm d}\sigma{\rm d}t
 +2 s\lambda^{2}\int_Q \varphi a^{2}_{6}\left(\nabla\psi\cdot\nabla{w}\right)^2{\rm d}x{\rm d}t
 \\ \nonumber&\hspace{0.5cm}+4 s\lambda\int_Q \varphi a_{6}\left(\nabla a_{6}\cdot\nabla{w}\right)
 \left(\nabla\psi\cdot\nabla{w}\right){\rm d}x{\rm d}t
 +2(2-n) s\lambda\int_Q \varphi a^{2}_{6}\left|\nabla{w}\right|^{2}{\rm d}x{\rm d}t
 \\ &\hspace{0.5cm}
 -s\lambda^{2}\int_Q \varphi\left|\nabla\psi\right|^{2}a^{2}_{6}\left|\nabla{w}\right|^{2}{\rm d}x{\rm d}t
 -2 s\lambda\int_Q \varphi a_{6}\left(\nabla\psi\cdot\nabla a_{6}\right)
 \left|\nabla{w}\right|^{2}{\rm d}x{\rm d}t.
 \end{align}
 \begin{align}\label{j6}
  \nonumber&\int_QJ_{6}{\rm d}x{\rm d}t=-\int_Qa_{6}\left(\Delta{w}\right)A\left({\vartheta},{w}\right){\rm d}x{\rm d}t
 \\\nonumber&=-s\lambda^{2}\int_Q\varphi\left|\nabla\psi\right|^{2}a_{5}a_{6}\left(\Delta{w}\right){\vartheta}{\rm d}x{\rm d}t
 -s\lambda^{2}\int_Q \varphi\left|\nabla\psi\right|^{2}a^{2}_{6}\left(\Delta{w}\right){w}{\rm d}x{\rm d}t
 \\\nonumber&\hspace{0.2cm}-2n s\lambda\int_Q\varphi a_{5}a_{6}\left(\Delta{w}\right){\vartheta}{\rm d}x{\rm d}t
 -2n s\lambda\int_Q\varphi a^{2}_{6}\left(\Delta{w}\right){w}{\rm d}x{\rm d}t
 \\\nonumber&\hspace{0.2cm}+s\lambda\int_Q \varphi\left(\partial_{t}\psi\right)a_{6}\left(\Delta{w}\right){\vartheta}{\rm d}x{\rm d}t
 \\\nonumber&=s\lambda^{3}\int_Q\varphi\left|\nabla\psi\right|^{2}a_{5}a_{6}\left(\nabla\psi\cdot
 \nabla{w}\right){\vartheta}{\rm d}x{\rm d}t
 \\\nonumber&\hspace{0.2cm}+s\lambda^{2}\int_Q\varphi\left(\nabla\left(\left|\nabla\psi\right|^{2}a_{5}a_{6}
 \right)\cdot\nabla{w}\right){\vartheta}{\rm d}x{\rm d}t+s\lambda^{2}\int_Q\varphi\left|\nabla\psi\right|^{2}a_{5}a_{6}
 \left(\nabla{w}\cdot\nabla{\vartheta}\right){\rm d}x{\rm d}t
 \\\nonumber&\hspace{0.2cm}
 +s\lambda^{3}\int_Q\varphi\left|\nabla\psi\right|^{2}a^{2}_{6}\left(\nabla\psi\cdot
 \nabla{w}\right){w}{\rm d}x{\rm d}t
 +s\lambda^{2}\int_Q\varphi\left(\nabla\left(\left|\nabla\psi\right|^{2}a^{2}_{6}
 \right)\cdot\nabla{w}\right){w}{\rm d}x{\rm d}t
  \\\nonumber&\hspace{0.2cm}+s\lambda^{2}\int_Q\varphi\left|\nabla\psi\right|^{2}a^{2}_{6}
 \left|\nabla{w}\right|^{2}{\rm d}x{\rm d}t
 +2n s\lambda^{2}\int_Q\varphi a_{5}a_{6}\left(\nabla\psi\cdot\nabla{w}\right){\vartheta}{\rm d}x{\rm d}t
 \\\nonumber&\hspace{0.2cm}+2ns\lambda\int_Q\varphi\left(\nabla\left(a_{5}a_{6}\right)
 \cdot\nabla{w}\right){\vartheta}{\rm d}x{\rm d}t
 +2n s\lambda\int_Q\varphi a_{5}a_{6}\left(\nabla{w}\cdot\nabla{\vartheta}\right){\rm d}x{\rm d}t
 \\\nonumber&\hspace{0.2cm}
 +2n s\lambda^{2}\int_Q\varphi a^{2}_{6}\left(\nabla\psi\cdot\nabla{w}\right){w}{\rm d}x{\rm d}t
 +4n s\lambda\int_Q\varphi a_{6}\left(\nabla a_{6}\cdot\nabla{w}\right){w}{\rm d}x{\rm d}t
  \\\nonumber&\hspace{0.2cm}+2n s\lambda\int_Q\varphi a^{2}_{6}\left|\nabla{w}\right|^{2}{\rm d}x{\rm d}t
 -s\lambda^{2}\int_Q\varphi\left(\partial_{t}\psi\right)a_{6}\left(\nabla\psi\cdot
 \nabla{w}\right){\vartheta}{\rm d}x{\rm d}t
 \\\nonumber&\hspace{0.2cm}-s\lambda\int_Q\varphi\left(\partial_{t}\psi\right)\left(\nabla a_{6}\cdot
 \nabla{w}\right){\vartheta}{\rm d}x{\rm d}t-s\lambda\int_Q\varphi\left(\partial_{t}\psi\right)a_{6}
 \left(\nabla{w}\cdot\nabla{\vartheta}\right){\rm d}x{\rm d}t
 \\\nonumber &\geq
  s\lambda^{2}\int_Q\varphi\left|\nabla\psi\right|^{2}a^{2}_{6}
 \left|\nabla{w}\right|^{2}{\rm d}x{\rm d}t
 +2n s\lambda\int_Q\varphi a^{2}_{6}\left|\nabla{w}\right|^{2}{\rm d}x{\rm d}t
 \\\nonumber &\hspace{0.2cm}+s\lambda^{2}\int_Q\varphi\left|\nabla\psi\right|^{2}a_{5}a_{6}
 \left(\nabla{w}\cdot\nabla{\vartheta}\right){\rm d}x{\rm d}t
 +2n s\lambda\int_Q\varphi a_{5}a_{6}\left(\nabla{w}\cdot\nabla{\vartheta}\right){\rm d}x{\rm d}t
  \\&\hspace{0.2cm}-s\lambda\int_Q\varphi\left(\partial_{t}\psi\right)a_{6}
 \left(\nabla{w}\cdot\nabla{\vartheta}\right){\rm d}x{\rm d}t
 -{C_3}s\lambda^{3}\int_Q\varphi\left|\nabla{w}\right|(|{\vartheta}|+|{w}|){\rm d}x{\rm d}t.
 \end{align}
 \begin{align}\label{j7}
 \nonumber &\int_QJ_{7}{\rm d}x{\rm d}t=-\frac{1}{2}s^{2}\lambda^{2}\int_Q\varphi^{2}\left|\nabla\psi\right|^{2}a_{5}
 \left\{\partial_{t}\left({\vartheta}^2\right)\right\}{\rm d}x{\rm d}t
 \\&=s^{2}\lambda^{3}\int_Q \varphi^{2}\left|\nabla\psi\right|^{2}\left(\partial_{t}\psi\right)a_{5}{\vartheta}^{2}{\rm d}x{\rm d}t
 \geq -{C_4}s^{2}\lambda^{3}\int_Q\varphi^{2}{\vartheta}^{2}{\rm d}x{\rm d}t.
 \end{align}
 \begin{align}\label{j8}
 \nonumber &\int_Q J_{8}{\rm d}x{\rm d}t=-2 s^{3}\lambda^{3}\int_Q\varphi^{3}\left|\nabla\psi\right|^{2}a_{5}\left(\nabla\psi\cdot\left(a_{5}
 \nabla{\vartheta}+a_{6}\nabla{w}\right)\right){\vartheta}{\rm d}x{\rm d}t
 \\\nonumber&=3 s^{3}\lambda^{4}\int_Q\varphi^{3}\left|\nabla\psi\right|^{4}a^{2}_{5}{\vartheta}^{2}{\rm d}x{\rm d}t
 +s^{3}\lambda^{3}\int_Q\varphi^{3}\sum_{j=1}^{n}\partial_{j}\left\{\left(\partial_{j}\psi\right)
 \left|\nabla\psi\right|^{2}a^{2}_{5}\right\}{\vartheta}^{2}{\rm d}x{\rm d}t
 \\ \nonumber&\hspace{0.5cm}-2 s^{3}\lambda^{3}\int_Q\varphi^{3}\left|\nabla\psi\right|^{2}a_{5}a_{6}\left(\nabla\psi\cdot
 \nabla{w}\right){\vartheta}{\rm d}x{\rm d}t
 \\ \nonumber&\geq 3 s^{3}\lambda^{4}\int_Q\varphi^{3}\left|\nabla\psi\right|^{4}a^{2}_{5}{\vartheta}^{2}{\rm d}x{\rm d}t
 -2 s^{3}\lambda^{3}\int_Q\varphi^{3}\left|\nabla\psi\right|^{2}a_{5}a_{6}\left(\nabla\psi\cdot
 \nabla{w}\right){\vartheta}{\rm d}x{\rm d}t
 \\&\hspace{0.5cm}-{C_5}s^{3}\lambda^{3}\int_Q\varphi^{3}{\vartheta}^{2}{\rm d}x{\rm d}t.
\end{align}
 \begin{align}\label{j9}
 \nonumber \int_Q&J_{9}{\rm d}x{\rm d}t
 =-s^{2}\lambda^{2}\int_{Q}\varphi^{2}\left|\nabla\psi\right|^{2}a_{5}{\vartheta}A\left({\vartheta},{w}\right){\rm d}x{\rm d}t
 \\\nonumber&=-s^{3}\lambda^{4}\int_Q\varphi^{3}\left|\nabla\psi\right|^{4}a^{2}_{5}{\vartheta}^{2}{\rm d}x{\rm d}t
 -s^{3}\lambda^{4}\int_Q\varphi^{3}\left|\nabla\psi\right|^{4}a_{5}a_{6}{\vartheta}{w}{\rm d}x{\rm d}t
 \\&\hspace{0.5cm}-{C_6} s^{3}\lambda^{3}\int_Q\varphi^{3}{\vartheta}^{2}{\rm d}x{\rm d}t
 -{C_6} s^{3}\lambda^{3}\int_Q\varphi^{3}|{\vartheta}||{w}|{\rm d}x{\rm d}t.
 \end{align}
 \begin{align}\label{j10}
 \nonumber &\int_QJ_{10}{\rm d}x{\rm d}t
 =-s^{2}\lambda^{2}\int_Q\varphi^{2}\left|\nabla\psi\right|^{2}a_{6}{w}\left(\partial_{t}{\vartheta}\right){\rm d}x{\rm d}t
 \\\nonumber&=s^{2}\lambda^{2}\int_Q\varphi^{2}\left|\nabla\psi\right|^{2}a_{6}{\vartheta}\left(\partial_{t}{w}\right){\rm d}x{\rm d}t
 +2 s^{2}\lambda^{3}\int_Q\varphi^{2}\left(\partial_{t}\psi\right)\left|\nabla\psi\right|^{2}a_{6}{\vartheta}{w}{\rm d}x{\rm d}t
 \\&\geq -{C_7}s^{2}\lambda^{2}\int_Q\varphi^{2}|{\vartheta}|\left|\partial_{t}{w}\right|{\rm d}x{\rm d}t
 -{C_7} s^{2}\lambda^{3}\int_Q\varphi^{2}|{\vartheta}||{w}|{\rm d}x{\rm d}t.
 \end{align}
 \begin{align}\label{j11}
 \nonumber &\int_QJ_{11}{\rm d}x{\rm d}t
 =-2 s^{3}\lambda^{3}\int_Q\varphi^{3}\left|\nabla\psi\right|^{2}a_{6}\left(\nabla\psi\cdot\left(a_{5}
 \nabla{\vartheta}+a_{6}\nabla{w}\right)\right){w}{\rm d}x{\rm d}t
 \\\nonumber&=
 2s^{3}\lambda^{3}\int_Q\varphi^{3}\sum_{j=1}^{n}\partial_{j}\left\{\left|\nabla\psi\right|^{2}
 \left(\partial_{j}\psi\right)a_{5}a_{6}\right\}{\vartheta}{w}{\rm d}x{\rm d}t
 \\\nonumber&\hspace{0.2cm}
 +6 s^{3}\lambda^{4}\int_Q\varphi^{3}\left|\nabla\psi\right|^{4}a_{5}a_{6}{\vartheta}{w}{\rm d}x{\rm d}t
 +2 s^{3}\lambda^{3}\int_Q\varphi^{3}\left|\nabla\psi\right|^{2}a_{5}a_{6}
 \left(\nabla\psi\cdot\nabla{w}\right){\vartheta}{\rm d}x{\rm d}t
 \\\nonumber&\hspace{0.2cm}
 +3 s^{3}\lambda^{4}\int_Q\varphi^{3}\left|\nabla\psi\right|^{4}a^{2}_{6}{w}^{2}{\rm d}x{\rm d}t
 +s^{3}\lambda^{3}\int_Q\varphi^{3}\sum_{j=1}^{n}\partial_{j}\left\{\left|\nabla\psi\right|^{2}
 \left(\partial_{j}\psi\right)a^{2}_{6}\right\}{w}^{2}{\rm d}x{\rm d}t
\\\nonumber&\geq 6 s^{3}\lambda^{4}\int_Q\varphi^{3}\left|\nabla\psi\right|^{4}a_{5}a_{6}{\vartheta}{w}{\rm d}x{\rm d}t
 + 2 s^{3}\lambda^{3}\int_Q\varphi^{3}\left|\nabla\psi\right|^{2}a_{5}a_{6}
 \left(\nabla\psi\cdot\nabla{w}\right){\vartheta}{\rm d}x{\rm d}t
 \\\nonumber&\hspace{0.2cm}
 +3 s^{3}\lambda^{4}\int_Q\varphi^{3}\left|\nabla\psi\right|^{4}a^{2}_{6}{w}^{2}{\rm d}x{\rm d}t.
 \\&\hspace{0.2cm}-{C_8}s^{3}\lambda^{3}\int_Q\varphi^{3}|{\vartheta}||{w}|{\rm d}x{\rm d}t
 -{C_8}s^{3}\lambda^{3}\int_Q\varphi^{3}{w}^{2}{\rm d}x{\rm d}t.
 \end{align}
 \begin{align}\label{j12}
 \nonumber \int_Q &J_{12}{\rm d}x{\rm d}t
 =-s^{2}\lambda^{2}\int_Q\varphi^{2}\left|\nabla\psi\right|^{2}a_{6}{w}A\left({\vartheta},{w}\right){\rm d}x{\rm d}t
 \\\nonumber&\geq-s^{3}\lambda^{4}\int_Q\varphi^{3}\left|\nabla\psi\right|^{4}a_{5}a_{6}{\vartheta}{w}{\rm d}x{\rm d}t
 -s^{3}\lambda^{4}\int_Q\varphi^{3}\left|\nabla\psi\right|^{4}a^{2}_{6}{w}^{2}{\rm d}x{\rm d}t
 \\&\hspace{0.5cm}-{C_9}s^{3}\lambda^{3}\int_Q\varphi^{3}|{\vartheta}||{w}|{\rm d}x{\rm d}t
 -{C_9}s^{3}\lambda^{3}\int_Q\varphi^{3}{w}^{2}{\rm d}x{\rm d}t.
 \end{align}
 We combine (\ref{jk})-(\ref{j3}) and (\ref{j4})-(\ref{j12}). The sum of the third and ninth terms in the right hand side of (\ref{j2}) and the second and eighth terms in the right hand side of (\ref{j5}) equals to
 $$
 2s\lambda^{2}\int_Q\varphi
 \left|a_{5}\left(\nabla\psi\cdot\nabla{\vartheta}\right)+a_{6}\left(\nabla\psi\cdot\nabla{w}\right)\right|^{2}{\rm d}x{\rm d}t.
 $$
 The sum of the sixth term in the right hand side of (\ref{j2}) and the first term in the right hand side of (\ref{j3}) equals to $0$. The sum of the eighth term in the right hand side of (\ref{j2}) and the third term in the right hand side of (\ref{j3}) equals to $0$. The sum of the twelfth term in the right hand side of (\ref{j2}) and the sixth term in the right hand side of (\ref{j5}) equals to $0$. The sum of the second term in the right hand side of (\ref{j3}), the fourth term in the right hand side of (\ref{j5}), and the third term in the right hand side of (\ref{j6}) equals to $0$. The sum of the eleventh term in the right hand side of (\ref{j5}) and the first term in the right hand side of (\ref{j6}) equals to $0$. The sum of the second term in the right hand side of (\ref{j8}) and the second term in the right hand side of (\ref{j11}) equals to $0$. Furthermore, we combine the eleventh term in right hand side of (\ref{j2}), the fourth term in the right term in the right hand of (\ref{j3}), the first term in the right hand side of (\ref{j5}), and the fourth term in the right hand side of (\ref{j6}). We combine the tenth term in the right hand side of (\ref{j5}) and the second term in the right hand side of (\ref{j6}). We combine the first term in the right hand side of (\ref{j8}) and the first term in the right hand side of (\ref{j9}). We combine the second term in the right hand side of (\ref{j9}), the first term in the right hand side of (\ref{j11}), and the first term in the right hand side of (\ref{j12}).  We combine the third term in the right hand side of (\ref{j11}) and the second term in the right hand side of (\ref{j12}).
  Therefore, we can obtain
 \begin{align}\label{wg}
 \nonumber &\int_QL_{1}\left({\vartheta},{w}\right)L_{2}\left({\vartheta},{w}\right){\rm d}x{\rm d}t
 \geq-s\lambda\int_\Sigma \varphi\left(\nabla\psi\cdot\nu\right)\left(a_{5}
 \frac{\partial {\vartheta}}{\partial\nu}+a_{6}
 \frac{\partial{w}}{\partial \nu}\right)^2{\rm d}\sigma{\rm d}t
 \\\nonumber&\hspace{0.2cm}+2s\lambda^{2}\int_Q\varphi
 \left|a_{5}\left(\nabla\psi\cdot\nabla{\vartheta}\right)+a_{6}\left(\nabla\psi\cdot\nabla{w}\right)\right|^{2}{\rm d}x{\rm d}t+\int_Q\left(\nabla a_{5}\cdot \nabla {\vartheta}\right)\left(\partial_{t}{\vartheta}\right)
 {\rm d}x{\rm d}t
 \\\nonumber&\hspace{0.2cm}+4 s\lambda\int_Q\varphi a_{5}\left(\nabla a_{5}\cdot\nabla{\vartheta}\right)
 \left(\nabla\psi\cdot\nabla{\vartheta}\right){\rm d}x{\rm d}t
 +4 s\lambda\int_Q\varphi a^{2}_{5}\left|\nabla{\vartheta}\right|^{2}{\rm d}x{\rm d}t
 \\ \nonumber &\hspace{0.2cm}-2s\lambda\int_Q \varphi a_{5}\left(\nabla\psi\cdot\nabla a_{5}\right)
  \left|\nabla{\vartheta}\right|^2{\rm d}x{\rm d}t
 +2 s\lambda\int_Q \varphi\left(\nabla\left(a_{5}a_{6}\right)\cdot\nabla{\vartheta}\right)\left(\nabla\psi\cdot\nabla{w}\right)
 {\rm d}x{\rm d}t
 \\ \nonumber &\hspace{0.2cm}+8 s\lambda\int_Q \varphi a_{5}a_{6}\left(\nabla{\vartheta}\cdot\nabla{w}\right){\rm d}x{\rm d}t
 -s\lambda\int_Q\varphi\left(\partial_{t}\psi\right)
 a_{5}\left|\nabla{\vartheta}\right|^{2}{\rm d}x{\rm d}t
 \\\nonumber&\hspace{0.2cm}-s\lambda\int_Q \varphi\left(\partial_{t}\psi\right)a_{6}\left(\nabla\Theta\cdot
 \nabla p\right){\rm e}^{2s\varphi}{\rm d}x{\rm d}t
 +\int_Q \frac{a_{6}}{a_{3}}\left(\partial_{t}p\right)\left(\partial_{t}\Theta\right){\rm e}^{2s\varphi}{\rm d}x{\rm d}t
 \\ \nonumber&\hspace{0.2cm}+s\lambda\int_Q\varphi\left(\partial_{t}\psi\right)\frac{ a_{4}a_{6}}{a_{3}}\left|\partial_{t}p\right|^{2}
 {\rm e}^{2s\varphi}{\rm d}x{\rm d}t
 +2 s\lambda\int_Q \varphi a_{6}\left(\nabla\psi\cdot\nabla\Theta\right)
 \left(\partial_{t}p\right){\rm e}^{2s\varphi}{\rm d}x{\rm d}t
 \\ \nonumber&\hspace{0.2cm}+2s\lambda\int_Q\varphi\left(\nabla\left(a_{5}a_{6}\right)\cdot\nabla w\right)
 \left(\nabla\psi\cdot\nabla{\vartheta}\right){\rm d}x{\rm d}t
 \\ \nonumber&\hspace{0.2cm}
 -2s\lambda\int_Q\varphi\left(\nabla\left(a_{5}a_{6}\right)\cdot\nabla\psi\right)
 \left(\nabla{\vartheta}\cdot\nabla{w}\right){\rm d}x{\rm d}t
 \\ \nonumber&\hspace{0.2cm}+4 s\lambda\int_Q \varphi a_{6}\left(\nabla a_{6}\cdot\nabla{w}\right)
 \left(\nabla\psi\cdot\nabla{w}\right){\rm d}x{\rm d}t+4 s\lambda\int_Q \varphi a^{2}_{6}\left|\nabla{w}\right|^{2}{\rm d}x{\rm d}t
 \\ \nonumber&\hspace{0.2cm}
 -2 s\lambda\int_Q \varphi a_{6}\left(\nabla \psi\cdot \nabla a_{6}\right)
 \left|\nabla{w}\right|^{2}{\rm d}x{\rm d}t
 -s\lambda\int_Q\varphi\left(\partial_{t}\psi\right)a_{6}
 \left(\nabla{w}\cdot\nabla{\vartheta}\right){\rm d}x{\rm d}t
 \\ \nonumber&\hspace{0.2cm}+2 s^{3}\lambda^{4}\int_Q\varphi^{3}\left|\nabla\psi\right|^{4}a^{2}_{5}{\vartheta}^{2}{\rm d}x{\rm d}t
 +4s^{3}\lambda^{4}\int_Q\varphi^{3}\left|\nabla\psi\right|^{4}a_{5}a_{6}{\vartheta}{w}{\rm d}x{\rm d}t
 \\ \nonumber&\hspace{0.2cm}+2s^{3}\lambda^{4}\int_Q\varphi^{3}\left|\nabla\psi\right|^{4}a^{2}_{6}{w}^{2}{\rm d}x{\rm d}t
 -{C_{10}}\int_Q\left|\partial_{t}{w}\right|\left(|{\vartheta}|+\left|\nabla{\vartheta}\right|\right){\rm d}x{\rm d}t
 \\ \nonumber&\hspace{0.2cm} -{C_{10}}s^3\lambda^3\int_Q \varphi^3|p||\Theta|{\rm e}^{2s\varphi}{\rm d}x{\rm d}t
-{C_{10}}s^{2}\lambda^{2}\int_Q\varphi^{2}\left|\partial_{t}p\right||\Theta|{\rm e}^{2s\varphi}{\rm d}x{\rm d}t
 \\  \nonumber&\hspace{0.2cm}-{C_{10}}s\lambda^{2}\int_Q\varphi |p|\left|\nabla\Theta\right|
 {\rm e}^{2s\varphi}{\rm d}x
 {\rm d}t
 -{C_{10}}\int_Q\left|\partial_{t}p\right||g|{\rm e}^{2s\varphi}{\rm d}x{\rm d}t
 \\&\nonumber\hspace{0.2cm}-{C_{10}}s\lambda^3\int_Q\varphi\left(\left|\nabla {\vartheta}\right|+\left|\nabla{w}\right|\right)\left(\left|{\vartheta}\right|+\left|{w}\right|\right){\rm d}x{\rm d}t
 \\&\hspace{0.2cm}
 -{C_{10}}s^{3}\lambda^{3}\int_Q\varphi^{3}\left({\vartheta}^2+{w}^{2}\right){\rm d}x{\rm d}t
 -{C_{10}}s^{2}\lambda^{2}\int_Q\varphi^{2}|{\vartheta}|\left|\partial_{t}{w}\right|{\rm d}x{\rm d}t.
 \end{align}

For $s>1$, we have, $\left(\Theta, p\right)={\rm e}^{-s\varphi}\left({\vartheta}, {w}\right)$,
\begin{align}\label{trans}
\nonumber\partial_{t}\Theta&=\partial_{t}\left({\rm e}^{-s\varphi}{\vartheta}\right)
=-s\lambda\varphi\left(\partial_{t}\psi\right){\rm e}^{-s\varphi}{\vartheta}+{\rm e}^{-s\varphi}\partial_{t}{\vartheta}
\\\nonumber&=2s\lambda\varphi\beta\left(t-t_0\right){\rm e}^{-s\varphi}{\vartheta}+{\rm e}^{-s\varphi}\partial_{t}{\vartheta},
\\\nonumber \partial_{j}\Theta&=\partial_{j}\left({\rm e}^{-s\varphi}{\vartheta}\right)
=-s\lambda\varphi\left(\partial_{j}\psi\right){\rm e}^{-s\varphi}{\vartheta}+{\rm e}^{-s\varphi}\partial_{j}{\vartheta}
\\\nonumber&=-2s\lambda\varphi\left(x_{j}-x^{j}_{0}\right){\rm e}^{-s\varphi}{\vartheta}+{\rm e}^{-s\varphi}\partial_{j}{\vartheta},
\\\nonumber \partial_{t}p&=\partial_{t}\left({\rm e}^{-s\varphi}{w}\right)
=-s\lambda\varphi\left(\partial_{t}\psi\right){\rm e}^{-s\varphi}{w}+{\rm e}^{-s\varphi}\partial_{t}{w}
\\\nonumber&=2s\lambda\varphi\beta\left(t-t_0\right){\rm e}^{-s\varphi}{w}+{\rm e}^{-s\varphi}\partial_{t}{w},
\\\nonumber \partial_{j}p&=\partial_{j}\left({\rm e}^{-s\varphi}{w}\right)
=-s\lambda\varphi\left(\partial_{j}\psi\right){\rm e}^{-s\varphi}{w}+{\rm e}^{-s\varphi}\partial_{j}{w}
\\&=-2s\lambda\varphi\left(x_{j}-x^{j}_{0}\right){\rm e}^{-s\varphi}{w}+{\rm e}^{-s\varphi}\partial_{j}{w},\ \ \mbox{in}\ Q.
\end{align}
By the Cauchy-Schwarz inequality, $\big(a_1, a_2, a_3, a_4\big)\in \mathcal{U}$, and (\ref{trans}), we estimate every term of (\ref{wg}).
 We note that
 $\|a_{5}\|_{C(\overline{\Omega})}\leq2M^{2}_{0}$,
 $\|a_{6}\|_{C(\overline{\Omega})}\leq M^{2}_{0}$,
 $\|\partial_{j}a_{5}\|_{C(\overline{\Omega})}\leq3M_{0}M_{1}$, $\|\partial_{j}a_{6}\|_{C(\overline{\Omega})}\leq2M_{0}M_{1}$,
 and
 $\|\partial_{j}(a_{5}a_{6})\|_{C(\overline{\Omega})}
 =\|a_{5}\partial_{j}a_{6}+a_{6}\partial_{j}a_{5}\|_{C(\overline{\Omega})}\leq7M^{3}_{0}M_{1}$, $j=1$, $2$, $\cdots$, $n$.
 By the Cauchy-Schwarz inequality, we can get, for any $\varepsilon_{1}>0$,
\begin{align}\label{1}
\left|\left(\nabla a_{5}\cdot \nabla {\vartheta}\right)\left(\partial_{t}{\vartheta}\right)\right|
\leq\frac{9M^{2}_{0}M^{2}_{1}}{2}\varepsilon^{-1}_{1}s\varphi\left|\nabla{\vartheta}\right|^{2}
+\frac{n\varepsilon_{1}}{2s\varphi}\left|\partial_{t}{\vartheta}\right|^{2},\ \ \mbox{in}\ Q.
\end{align}
\begin{align}
\nonumber\left|4 s\lambda\varphi a_{5}\left(\nabla a_{5}\cdot\nabla{\vartheta}\right)
 \left(\nabla\psi\cdot\nabla{\vartheta}\right)
 \right|
 &=\left|8s\lambda\varphi a_{5}\left(\nabla a_{5}\cdot\nabla{\vartheta}\right)
 \left(\left(x-x_{0}\right)\cdot\nabla{\vartheta}\right)\right|
\\&\leq48nM^{3}_{0}M_{1}M^{\frac{1}{2}}s\lambda\varphi\left|\nabla{\vartheta}\right|^{2},\ \ \mbox{in}\ Q.
\end{align}
\begin{align}\label{esti11}
 4s\lambda\varphi a^{2}_{5}\left|\nabla{\vartheta}\right|^{2}\geq 4\sigma_1^{2}s\lambda\varphi\left|\nabla{\vartheta}\right|^{2},\ \ \mbox{in}\ Q.
 \end{align}
 \begin{align}
 &\nonumber\left|-2s\lambda\varphi a_{5}\left(\nabla\psi\cdot \nabla a_{5}\right)
 \left|\nabla{\vartheta}\right|^{2}\right|
 =\left|4s\lambda\varphi a_{5}\left((x-x_{0})\cdot\nabla a_{5}\right)\right|
 \left|\nabla{\vartheta}\right|^{2}
 \\&\hspace{0.5cm}\leq24nM^{3}_{0}M_{1}M^{\frac{1}{2}}s\lambda\varphi\left|\nabla{\vartheta}\right|^{2},\ \ \mbox{in}\ Q.
 \end{align}
Similarly,
\begin{align}
&\left|-2s\lambda\varphi a_{6}\left(\nabla\psi\cdot\nabla a_{6}\right)
\left|\nabla{w}\right|^{2}\right|
\leq8nM^{3}_{0}M_{1}M^{\frac{1}{2}}
s\lambda\varphi\left|\nabla{w}\right|^{2},\ \ \mbox{in}\ Q.
\end{align}
\begin{align}
\left|4 s\lambda \varphi a_{6}\left(\nabla a_{6}\cdot\nabla{w}\right)
 \left(\nabla\psi\cdot\nabla{w}\right)
\right|
\leq16nM^{3}_{0}M_{1}M^{\frac{1}{2}}s\lambda\varphi\left|\nabla{w}\right|^{2},\ \ \mbox{in}\ Q.
\end{align}
 Similarly, we have
 \begin{align}
 \nonumber&\left|2 s\lambda \varphi\left(\nabla\left(a_{5}a_{6}\right)\cdot\nabla{\vartheta}\right)\left(\nabla\psi\cdot\nabla{w}\right)\right|
 +\left|2s\lambda\varphi\left(\nabla\left(a_{5}a_{6}\right)\cdot\nabla w\right)
 \left(\nabla\psi\cdot\nabla{\vartheta}\right)\right|
 \\ \nonumber&\hspace{0.5cm}\hspace{0.2cm}+\left|-2s\lambda\varphi\left(\nabla\left(a_{5}a_{6}\right)\cdot\nabla\psi\right)
 \left(\nabla{\vartheta}\cdot\nabla{w}\right)\right|
 \\ &\hspace{0.5cm}\leq 42 nM^{3}_{0}M_{1}M^{\frac{1}{2}}s\lambda\varphi
 \left(\left|\nabla{\vartheta}\right|^{2}+\left|\nabla{w}\right|^{2}\right),\ \ \mbox{in}\ Q.
 \end{align}
 \begin{align}
 \left|8 s\lambda \varphi a_{5}a_{6}\left(\nabla{\vartheta}\cdot\nabla{w}\right)\right|
 \leq8M^{4}_{0}s\lambda\varphi
 \left(\left|\nabla{\vartheta}\right|^{2}+\left|\nabla{w}\right|^{2}\right),\ \ \mbox{in}\ Q.
 \end{align}
 \begin{align}
 &\left|-s\lambda\varphi\left(\partial_{t}\psi\right)a_{5}\left|\nabla{\vartheta}\right|^{2}\right|
 =\left|2s\lambda\varphi\beta a_{5}\left(t-t_0\right)\right|\left|\nabla{\vartheta}\right|^{2}
 \leq2M^{2}_{0}s\lambda\varphi \beta T\left|\nabla{\vartheta}\right|^{2},
 \end{align}
 in $Q$.
 By (\ref{trans}), we have
 \begin{align}\label{trans1}
 &\nonumber \Big|2s\lambda\varphi a_{6}\left(\nabla\psi\cdot\nabla{\Theta}\right)
 \left(\partial_{t}{p}\right){\rm e}^{2s\varphi}\Big|
 \\\nonumber&=\Big|2s\lambda\varphi a_{6}\left(\nabla\psi\cdot\nabla{{\vartheta}}\right)
 \left(\partial_{t}{{w}}\right)-2s^{2}\lambda^{2}\varphi^{2}\left(\partial_{t}\psi\right)a_{6}
 \left(\nabla\psi\cdot\nabla{{\vartheta}}\right){w}
 \\\nonumber&\hspace{0.5cm}-2s^{2}\lambda^{2}\varphi^{2}\left|\nabla\psi\right|^{2}a_{6}\left(\partial_{t}{{w}}\right){\vartheta}
 +2s^{3}\lambda^{3}\varphi^{3}\left|\nabla\psi\right|^{2}\left(\partial_{t}\psi\right)a_{6}{w}{\vartheta}\Big|
 \\\nonumber&\leq\left|4s\lambda\varphi\left(\left(\partial_{t}{{w}}\right)\left(x-x_0\right)\cdot
 a_{6}\nabla{{\vartheta}}\right)\right|
 \\\nonumber&\hspace{0.5cm}+
 \left|8\beta\left(t-t_0\right)a_{6}\left(\left(x-x_{0}\right)\cdot\left(s^{\frac{1}{2}}\lambda^{\frac{1}{4}}
 \varphi^{\frac{1}{2}}\nabla{{\vartheta}}\right)\right)
 s^{\frac{3}{2}}\lambda^{\frac{7}{4}}\varphi^{\frac{3}{2}}{w}\right|
 \\ \nonumber &\hspace{0.5cm}+\left|\nabla\psi\right|^{2}\left|2a_{6}
 \left(s^{\frac{1}{2}}\lambda^{\frac{1}{4}}\varphi^{\frac{1}{2}}\partial_{t}{{w}}\right)
 s^{\frac{3}{2}}\lambda^{\frac{7}{4}}\varphi^{\frac{3}{2}}{\vartheta}\right|
 +s^{3}\lambda^{3}\varphi^{3}\left|\nabla\psi\right|^{2}\left|\partial_{t}\psi\right|\left|2a_{6}
 {w}{\vartheta}\right|
 \\\nonumber&\leq2s\lambda\varphi\left(M\left|\partial_{t}{{w}}\right|^{2}
 +M^{4}_{0}\left|\nabla{{\vartheta}}\right|^{2}\right)
 + C_{11}s\lambda^{\frac{1}{2}}\varphi\left(\left|\nabla{\vartheta}\right|^{2}
 +\left|\partial_{t}{{w}}\right|^{2}\right)
 \\&\hspace{0.5cm}+C_{11}s^{3}\lambda^{\frac{7}{2}}\varphi^{3}\left({\vartheta}^{2}+{w}^{2}\right),\ \ \mbox{in}\ Q.
 \end{align}
 Similarly to (\ref{trans1}), we have
 \begin{align}\label{ssss}
 &\nonumber\left|-s\lambda\varphi\left(\partial_{t}\psi\right)a_{6}\left(\nabla\Theta\cdot\nabla p\right){\rm e}^{2s\varphi}\right|
 \\\nonumber&=\Big|-s\lambda\varphi\left(\partial_{t}\psi\right)a_{6}\left(\nabla{\vartheta}\cdot\nabla{w}\right)
 +s^{2}\lambda^{2}\varphi^{2}\left(\partial_{t}\psi\right)a_{6}
 \left(\nabla\psi\cdot\nabla{\vartheta}\right){w}
 \\\nonumber&\hspace{0.5cm}+s^{2}\lambda^{2}\varphi^{2}\left(\partial_{t}\psi\right)a_{6}
 \left(\nabla\psi\cdot\nabla{w}\right){\vartheta}
 -s^{3}\lambda^{3}\varphi^{3}\left|\nabla\psi\right|^{2}\left(\partial_{t}\psi\right)a_{6}{\vartheta}{w}\Big|
 \\\nonumber&\leq\frac{1}{2}M^{2}_{0}s\lambda\varphi\beta T\left(\left|\nabla{\vartheta}\right|^{2}
 +\left|\nabla{{w}}\right|^{2}\right)
 +C_{12}s\lambda^{\frac{1}{2}}\varphi\left(\left|\nabla{\vartheta}\right|^{2}+\left|\nabla{w}\right|^{2}\right)
 \\&\hspace{0.5cm}+C_{12}s^{3}\lambda^{\frac{7}{2}}\varphi^{3}\left({\vartheta}^{2}+{w}^{2}\right),\ \ \mbox{in}\ Q.
 \end{align}
 Similarly, we can obtain
 \begin{equation}\label{ss}
 s^{2}\lambda^{2}\varphi^{2}|{\vartheta}|\left|\partial_{t}{w}\right|\leq C_{13}s\lambda^{\frac{1}{2}}\varphi\left|\partial_t{w}\right|^{2}
 +C_{13}s^{3}\lambda^{\frac{7}{2}}\varphi^{3}{\vartheta}^{2},
 \end{equation}
 \begin{align}\label{sss}
 \nonumber s^{2}\lambda^{2}\varphi^{2}\left|\partial_{t}p\right||\Theta|{\rm e}^{2s\varphi}&\leq s^{2}\lambda^{2}\varphi^{2}\left(\left|\partial_{t}{w}\right||{\vartheta}|
 +\left|-s\lambda\varphi\left(\partial_t\psi\right){w}\right||{\vartheta}|\right)
 \\&\leq C_{14}s\lambda^{\frac{1}{2}}\varphi\left|\partial_t{w}\right|^{2}
 +C_{14}s^{3}\lambda^{\frac{7}{2}}\varphi^{3}\left({\vartheta}^{2}+{w}^{2}\right),\ \ \mbox{in}\ Q.
 \end{align}
 \begin{align}
 &\nonumber\left|\frac{a_{6}}{a_{3}}\left(\partial_{t}{\Theta}\right)\left(\partial_{t}{p}\right){\rm e}^{2s\varphi}\right|
 =\Big|\frac{a_{6}}{a_{3}}\left(\partial_{t}{{\vartheta}}\right)\left(\partial_{t}{{w}}\right)
 -s\lambda\varphi\left(\partial_{t}\psi\right)\frac{a_{6}}{a_{3}}\left(\partial_{t}{{\vartheta}}\right){w}
 \\\nonumber&\hspace{0.5cm}-s\lambda\varphi\left(\partial_{t}\psi\right)\frac{a_{6}}{a_{3}}\left(\partial_{t}{{w}}\right){\vartheta}
 +s^{2}\lambda^{2}\varphi^{2}\left(\partial_{t}\psi\right)^{2}\frac{a_{6}}{a_{3}}{\vartheta}{w}\Big|
 \\\nonumber&\leq\left|\left(\frac{1}{\sqrt{6s\varphi M}}
 \partial_{t}{{\vartheta}}\right)\left(\frac{M^{2}_{0}}{\sigma_{1}}\sqrt{6s\varphi M}\partial_{t}{{w}}\right)\right|
 \\\nonumber&\hspace{0.5cm}+\frac{M^{2}_{0}\beta T}{\sigma_{1}}\left|
 \left(\sqrt{\frac{1}{s\lambda\varphi}}\partial_{t}{{\vartheta}}\right)
 \left(s\lambda\varphi\sqrt{s\lambda\varphi}{w}\right)\right|
 +\frac{M^{2}_{0}\beta T}{\sigma_{1}}\left|\left(\partial_{t}{{w}}\right)
 \left(s\lambda\varphi{\vartheta}\right)\right|
 \\\nonumber&\hspace{0.5cm}+\frac{M^{2}_{0}\beta^2 T^2}
 {\sigma_{1}}\left|s^{2}\lambda^{2}\varphi^{2}{\vartheta}{w}\right|
 \\\nonumber&\leq\frac{1}{12s\varphi M}\left|\partial_{t}{{\vartheta}}\right|^{2}
 +\frac{3M^{4}_{0}M}{\sigma^{2}_{1}}s\varphi\left|\partial_{t}{{w}}\right|^{2}
 \\&\hspace{0.5cm}+C_{15}\left(\frac{1}{s\lambda\varphi}\left|\partial_{t}{{\vartheta}}\right|^{2}
 +s^{3}\lambda^{3}\varphi^{3}{w}^{2}+\left|\partial_{t}{{w}}\right|^{2}
 +s^{2}\lambda^{2}\varphi^{2}{\vartheta}^{2}\right),\ \ \mbox{in}\ Q.
 \end{align}
 \begin{align}
 &\nonumber\left|s\lambda\varphi\left(\partial_{t}\psi\right)\frac{a_{4}a_{6}}{a_{3}}
 \left|\partial_{t}p\right|^{2}{\rm e}^{2s\varphi}\right|
 =\left|2s\lambda\varphi\beta\left(t-t_0\right)\frac{a_{4}a_{6}}{a_{3}}\right|
 \left|\partial_{t}{w}+2s\lambda\varphi\beta\left(t-t_0\right){w}\right|^{2}
 \\&\leq\frac{2 M^{3}_{0}\beta T}{\sigma_{1}}s\lambda\varphi\left|\partial_{t}{w}\right|^{2}
 +C_{16}s^{3}\lambda^{3}\beta^3\varphi^{3}{w}^{2},\ \ \mbox{in}\ Q.
 \end{align}
 \begin{align}\label{esti12}
 4s\lambda\varphi a^{2}_{6}\left|\nabla{w}\right|^{2}\geq 0,\ \ \mbox{in}\ Q.
 \end{align}
 \begin{align}
 \left|-s\lambda\varphi\left(\partial_{t}\psi\right)a_{6}\left(\nabla{w}\cdot
 \nabla{\vartheta}\right)\right|
 \leq\frac{1}{2}M^{2}_{0}s\lambda \beta T\varphi\left(\left|\nabla{\vartheta}\right|^{2}+\left|\nabla{w}\right|^{2}\right),\ \ \mbox{in}\ Q.
 \end{align}
 \begin{align}
 2s^{3}\lambda^{4}\varphi^{3}\left|\nabla\psi\right|^{4}a^{2}_{5}{\vartheta}^{2}\geq2\sigma^{2}_{1}
 s^{3}\lambda^{4}\varphi^{3}\left|\nabla\psi\right|^{4}{\vartheta}^{2},\ \ \mbox{in}\ Q.
 \end{align}
 By the Cauchy-Schwarz inequality, we have, for any $\varepsilon_{2}>0$,
 \begin{align}
 \left|4s^{3}\lambda^{4}\varphi^{3}\left|\nabla\psi\right|^{4}a_{5}a_{6}{\vartheta}{w}\right|
  \leq4M^{4}_{0}s^{3}\lambda^{4}\varphi^{3}\left|\nabla\psi\right|^{4}
  \left(\varepsilon_{2}{\vartheta}^{2}+\varepsilon^{-1}_{2}{w}^{2}\right),\ \ \mbox{in}\ Q.
 \end{align}
 \begin{align}
 2s^{3}\lambda^{4}\varphi^{3}\left|\nabla\psi\right|^{4}a^{2}_{6}{w}^{2}\geq0,\ \ \mbox{in}\ Q.
 \end{align}
   Moreover, by (\ref{trans}), we have
 \begin{align}
 \nonumber s\lambda^{2}\varphi |p|\left|\nabla\Theta\right|{\rm e}^{2s\varphi}&\leq s\lambda^{2}\varphi\left(\left|\nabla{\vartheta}\right|+\left|-s\lambda\varphi\left(\nabla\psi\right){\vartheta}\right|\right)|{w}|
 \\&\leq C_{17}\left|\nabla{\vartheta}\right|^{2}
 +C_{17}s^2\lambda^4\varphi^2\left({\vartheta}^{2}+{w}^{2}\right),\ \ \mbox{in}\ Q.
 \end{align}
 \begin{align}
 \nonumber  |\partial_t p|\left|g\right|{\rm e}^{2s\varphi}&\leq \left(\left|\partial_t{w}\right|+\left|-s\lambda\varphi\left(\partial_t\psi\right){w}\right|\right)\left(|{g}|{\rm e}^{s\varphi}\right)
 \\&\leq C_{18}\left|\partial_t{w}\right|^{2}
 +C_{18}s^2\lambda^2\varphi^2{w}^{2}+C_{18}g^2{\rm e}^{2s\varphi},\ \ \mbox{in}\ Q.
 \end{align}
 \begin{align}
 \nonumber s\lambda^{3}\varphi &\left(\left|\nabla{\vartheta}\right|+\left|\nabla{w}\right|\right)\left(|{\vartheta}|+|{w}|\right)
 \\&\leq  C_{19}\left(\left|\nabla{\vartheta}\right|^{2}+\left|\nabla{w}\right|^{2}\right)
 +C_{19}s^2\lambda^6\varphi^2\left({\vartheta}^{2}+{w}^{2}\right),\ \ \mbox{in}\ Q.
 \end{align}
Moreover, by the second equality in (\ref{l12}) and the Cauchy-Schwarz inequality, we have, for any $\varepsilon_{3}>0$,
\begin{align*}
&\nonumber\left|L_{2}\left({\vartheta},{w}\right)\right|^{2}
=\left|\partial_{t}{\vartheta}+\left\{2s\lambda\varphi
\left(\nabla\psi\cdot\left(a_{5}\nabla{\vartheta}+a_{6}\nabla{w}\right)\right)
+A({\vartheta},{w})\right\}\right|^2
\\&=\left|\partial_{t}{\vartheta}\right|^{2}+b^{2}+2\left(\sqrt{\varepsilon_{3}}\partial_{t}{\vartheta}\right)
\frac{b}{\sqrt{\varepsilon_{3}}}
\geq\left(1-\varepsilon_{3}\right)\left|\partial_{t}{\vartheta}\right|^{2}
+\left(1-\varepsilon^{-1}_{3}\right)b^{2},
\end{align*}
 where
 $b=2s\lambda\varphi
 \left(\nabla\psi\cdot\left(a_{5}\nabla{\vartheta}+a_{6}\nabla{w}\right)\right)
 +A({\vartheta},{w})$. Noting
 \begin{align*}
 & b^{2}\leq 8s^{2}\lambda^{2}\varphi^{2}
 \left|a_{5}\left(\nabla\psi\cdot\nabla{\vartheta}\right)+a_{6}\left(\nabla\psi\cdot\nabla{w}\right)\right|^{2}
 +C_{20}s^{2}\lambda^{4}\varphi^{2}\left({\vartheta}^{2}+{w}^{2}\right),
 \end{align*}
 and choosing $\varepsilon_{3}=\frac{4s\varphi M}{4s\varphi M+1}$, we have
 \begin{align}\label{es2}
 \nonumber&\left|L_{2}\left({\vartheta},{w}\right)\right|^{2}
 \geq\frac{1}{4s\varphi M+1}\left|\partial_{t}{\vartheta}\right|^{2}-2s\lambda^2\varphi
 \left|a_{5}\left(\nabla\psi\cdot\nabla{\vartheta}\right)+a_{6}\left(\nabla\psi\cdot\nabla{w}\right)\right|^{2}
 \\\nonumber&\hspace{2.5cm}-C_{21}s\lambda^{4}\varphi\left({\vartheta}^{2}+{w}^{2}\right)
 \\\nonumber&\geq\frac{1}{5s\varphi M}\left|\partial_{t}{\vartheta}\right|^{2}-2s\lambda^2\varphi
 \left|a_{5}\left(\nabla\psi\cdot\nabla{\vartheta}\right)+a_{6}\left(\nabla\psi\cdot\nabla{w}\right)\right|^{2}
 \\&\hspace{0.5cm}-C_{21}s\lambda^{4}\varphi\left({\vartheta}^{2}+{w}^{2}\right),
 \end{align}
 for all $s>s_1\triangleq\max\{\frac{1}{M},1\}$.
 By (\ref{L12}), (\ref{wg}), (\ref{1})--(\ref{es2}),
 we obtain
 \begin{align*}
 &\nonumber\int_{Q}\left|\mathcal{L}\left({\vartheta},{w}\right)\right|^{2}{\rm d}x{\rm d}t
 \\&=\left\|L_{1}\left({\vartheta},{w}\right)\right\|^{2}_{L^2(Q)}+2\int_QL_{1}\left({\vartheta},{w}\right)
L_{2}\left({\vartheta},{w}\right){\rm d}x{\rm d}t+\left\|L_{2}\left({\vartheta},{w}\right)\right\|^{2}_{L^2(Q)}
 \\&\geq-2s\lambda\int_\Sigma \varphi\left(\nabla\psi\cdot\nu\right)\left(a_{5}
 \frac{\partial {\vartheta}}{\partial\nu}+a_{6}
 \frac{\partial{w}}{\partial \nu}\right)^2{\rm d}\sigma{\rm d}t+\left\|L_{1}\left({\vartheta},{w}\right)\right\|^{2}_{L^2(Q)}
 \\\nonumber&\hspace{0.2cm}+2s\lambda^{2}\int_Q\varphi
 \left|a_{5}\left(\nabla\psi\cdot\nabla{\vartheta}\right)+a_{6}\left(\nabla\psi\cdot\nabla{w}\right)\right|^{2}{\rm d}x{\rm d}t
 \\&\hspace{0.2cm}+\left(4\sigma^{2}_{1}-8M^{4}_{0}\varepsilon_{2}\right)
 s^{3}\lambda^{4}\int_{Q}\varphi^{3}\left|\nabla\psi\right|^{4}{\vartheta}^{2}{\rm d}x{\rm d}t
 -\frac{8M^{4}_{0}s^{3}\lambda^{4}}{\varepsilon_{2}}\int_{Q}\varphi^{3}\left|\nabla\psi\right|^{4}
 {w}^{2}{\rm d}x{\rm d}t
 \\ \nonumber &\hspace{0.2cm}
 +\left(\frac{1}{30M}-n\varepsilon_{1}\right)\int_{Q}\frac{1}{s\varphi}
 \left|\partial_{t}{\vartheta}\right|^{2}{\rm d}x{\rm d}t
 -4Ms\lambda\int_{Q}\varphi
 \left|\partial_{t}{{w}}\right|^{2}{\rm d}x{\rm d}t
 \\ \nonumber &\hspace{0.2cm}-\left(
 228nM^{3}_{0}M_{1}M^{\frac{1}{2}}+20M_0^4-8\sigma_1^2
 \right)
 s\lambda \int_{Q}\varphi\left|\nabla{\vartheta}\right|^{2}{\rm d}x{\rm d}t
 \\&\hspace{0.5cm}-\left(132nM_{1}M^{\frac{1}{2}}+16M_0
 \right)M^{3}_{0}s\lambda\int_{Q}\varphi\left|\nabla{w}\right|^{2}{\rm d}x{\rm d}t
 \\\nonumber&\hspace{0.2cm}-6M_0^2s\lambda\beta T\int_{Q}\varphi\left|\nabla{\vartheta}\right|^{2}{\rm d}x{\rm d}t
 -2M_0^2s\lambda\beta T\int_{Q}\varphi\left|\nabla{w}\right|^{2}{\rm d}x{\rm d}t
 \\&\hspace{0.2cm}-\frac{4M^{3}_{0}}{\sigma_{1}}s\lambda\beta T\int_{Q}\varphi\left|\partial_{t}{w}\right|^{2}{\rm d}x{\rm d}t
 -C_{22}\int_{Q}s\lambda^{\frac{1}{2}}\varphi\left(\left|\nabla{{\vartheta}}\right|^{2}
+\left|\nabla{{w}}\right|^{2}
+\left|\partial_{t}{{w}}\right|^{2}\right){\rm d}x{\rm d}t
\\&\hspace{0.2cm}-C_{22}\int_{Q}\left(s^{3}\lambda^{\frac{7}{2}}\varphi^{3}+s^{2}\lambda^{6}\varphi^{2}\right)
\left({\vartheta}^{2}+{w}^{2}\right){\rm d}x{\rm d}t-C_{22}\int_{Q}\frac{1}{s\lambda\varphi}
\left|\partial_{t}{{\vartheta}}\right|^{2}{\rm d}x{\rm d}t
\\&\hspace{0.2cm}
-C_{22}\int_{Q}g^{2}{\rm e}^{2s\varphi}{\rm d}x{\rm d}t,
\end{align*}
  for all $s>s_1$. Choosing $\varepsilon_{1}>0$ and $\varepsilon_{2}>0$ sufficiently small such that
 \begin{align*}
 4\sigma^{2}_{1}-8M^{4}_{0}\varepsilon_{2}=2\sigma^{2}_{1} \quad\mbox{and}\quad
 \frac{1}{30M}-n\varepsilon_{1}=\frac{1}{60M},
 \end{align*}
 that is,
 $\varepsilon_{1}=\frac{1}{60nM}$ and
 $\varepsilon_{2}=\frac{\sigma^{2}_{1}}{4M^{4}_{0}}$,
 we  arrive at
 \begin{align}\label{equation1}
 &\nonumber2\sigma^{2}_{1}s^{3}\lambda^{4}\int_{Q}\varphi^{3}\left|\nabla\psi\right|^{4}{\vartheta}^{2}{\rm d}x{\rm d}t
 +\frac{1}{60M}\int_{Q}\frac{1}{s\varphi}\left|\partial_{t}{\vartheta}\right|^{2}{\rm d}x{\rm d}t
 \\\nonumber&\hspace{0.5cm}+2 s\lambda^{2}\int_{Q}
 \varphi\left|a_{5}\left(\nabla\psi\cdot\nabla{\vartheta}\right)+a_{6}\left(\nabla\psi\cdot\nabla{w}\right)\right|^{2}{\rm d}x{\rm d}t
 \\\nonumber&\hspace{0.5cm}-2s\lambda\int_\Sigma \varphi\left(\nabla\psi\cdot\nu\right)\left(a_{5}
 \frac{\partial {\vartheta}}{\partial\nu}+a_{6}
 \frac{\partial{w}}{\partial \nu}\right)^2{\rm d}\sigma{\rm d}t +\left\|L_{1}\left({\vartheta},{w}\right)\right\|^{2}_{L^2(Q)}
 \\\nonumber&\leq C_{23}\int_{Q}\left|\mathcal{L}\left({\vartheta},{w}\right)\right|^{2}{\rm d}x{\rm d}t
 +\alpha_{1}s\lambda\int_{Q}\varphi\left|\nabla{{\vartheta}}\right|^{2}{\rm d}x{\rm d}t+\alpha_{2}s\lambda\int_{Q}\varphi\left|\nabla{{w}}\right|^{2}{\rm d}x{\rm d}t
 \\\nonumber&\hspace{0.5cm}+4M s\lambda\int_{Q}\varphi
 \left|\partial_{t}{{w}}\right|^{2}{\rm d}x{\rm d}t
 +\frac{32M^{8}_{0}}{\sigma^{2}_{1}}
 s^{3}\lambda^{4}\int_{Q}\varphi^{3}\left|\nabla\psi\right|^{4}{w}^{2}{\rm d}x{\rm d}t
 \\\nonumber&\hspace{0.2cm}+6M_0^2s\lambda\beta T\int_{Q}\varphi\left|\nabla{\vartheta}\right|^{2}{\rm d}x{\rm d}t
 +2M_0^2s\lambda\beta T\int_{Q}\varphi\left|\nabla{w}\right|^{2}{\rm d}x{\rm d}t
 \\\nonumber&\hspace{0.2cm}+\frac{4M^{3}_{0}}{\sigma_{1}}s\lambda\beta T\int_{Q}\varphi\left|\partial_{t}{w}\right|^{2}{\rm d}x{\rm d}t
 +C_{23}s\lambda^{\frac{1}{2}}\int_{Q}\varphi\left(\left|\nabla{{\vartheta}}\right|^{2}
 +\left|\nabla{{w}}\right|^{2}
  +\left|\partial_{t}{{w}}\right|^{2}\right){\rm d}x{\rm d}t
 \\\nonumber &\hspace{0.2cm}+C_{23}\int_{Q}\left(s^{3}\lambda^{\frac{7}{2}}\varphi^{3}+s^{2}\lambda^{6}\varphi^{2}\right)
 \left({\vartheta}^{2}+{w}^{2}\right){\rm d}x{\rm d}t+C_{23}\int_{Q}\frac{1}{s\lambda\varphi}
 \left|\partial_{t}{{\vartheta}}\right|^{2}{\rm d}x{\rm d}t
 \\&\hspace{0.2cm}
 +C_{23}\int_{Q}g^{2}{\rm e}^{2s\varphi}{\rm d}x{\rm d}t,
\end{align}
 for all $s>s_1$, where $\alpha_1$ and $\alpha_{2}$ are given by the first and second equality in (\ref{alpha}).

 Next we will estimate $s\lambda\int_{Q}\varphi\left|\nabla v\right|^{2}\me^{2s\varphi}{\rm d}x{\rm d}t$. Noting (\ref{l11}) and (\ref{l12}) and applying the integration by parts, similarly to (\ref{ssss}), we can obtain
 \begin{align}
 &\nonumber\int_Q\mathcal{L}\left({\vartheta},{w}\right)\left(s\lambda\varphi{\vartheta}\right){\rm d}x{\rm d}t
 \\\nonumber&=\int_Q\left\{-s\lambda\varphi a_{5}\left(\Delta{\vartheta}\right){\vartheta}-s\lambda\varphi a_{6}\left(\Delta{w}\right){\vartheta}+\frac{1}{2}s\lambda\varphi\left(\partial_{t}\left({\vartheta}^2\right)\right)\right\}{\rm d}x{\rm d}t
 \\\nonumber&\hspace{0.5cm}+\int_Q\Big\{-s^{3}\lambda^{3}\varphi^{3}\left|\nabla\psi\right|^{2}a_{5}{\vartheta}^{2}
-s^{3}\lambda^{3}\varphi^{3}\left|\nabla\psi\right|^{2}a_{6}{\vartheta}{w}
 \\\nonumber&\hspace{0.5cm}
 +2 s^{2}\lambda^{2}\varphi^{2}a_{5}\left(\nabla\psi\cdot\nabla{\vartheta}\right){\vartheta}
 \\\nonumber&\hspace{0.5cm}+2s^{2}\lambda^{2}\varphi^{2}a_{6}\left(\nabla\psi\cdot\nabla{w}\right){\vartheta}
 +s^{2}\lambda^{3}\varphi^{2}\left|\nabla\psi\right|^{2}\left(a_{5}{\vartheta}+a_{6}{w}\right){\vartheta}
 \\\nonumber&\hspace{0.5cm}+2ns^{2}\lambda^{2}\varphi^{2}\left(a_{5}{\vartheta}+a_{6}{w}\right){\vartheta}
 -s^{2}\lambda^{2}\varphi^{2}\left(\partial_{t}\psi\right){\vartheta}^{2}\Big\}{\rm d}x{\rm d}t
 \\\nonumber&\geq
 s\lambda\int_Q\varphi\left(\nabla a_{5}\cdot\nabla{\vartheta}\right){\vartheta}{\rm d}x{\rm d}t
 +s\lambda^{2}\int_Q\varphi a_{5}\left(\nabla\psi\cdot\nabla{\vartheta}\right){\vartheta}{\rm d}x{\rm d}t
 \\\nonumber&\hspace{0.5cm}+s\lambda\int_Q\varphi a_{5}\left|\nabla{\vartheta}\right|^{2}{\rm d}x{\rm d}t
  +s\lambda\int_Q\varphi\left(\nabla a_{6}\cdot\nabla{w}\right){\vartheta}{\rm d}x{\rm d}t
 \\\nonumber&\hspace{0.5cm}+s\lambda^{2}\int_Q\varphi a_{6}\left(\nabla\psi\cdot\nabla{w}\right){\vartheta}{\rm d}x{\rm d}t
 +s\lambda\int_Q\varphi a_{6}\left(\nabla{w}\cdot\nabla{\vartheta}\right){\rm d}x{\rm d}t
 \\\nonumber&\hspace{0.5cm}-\frac{1}{2}s\lambda^2\int_Q\varphi\left(\partial_t\psi\right){\vartheta}^2{\rm d}x{\rm d}t
 -C_{24}s^{3}\lambda^{3}\int_Q\varphi^{3}{w}^{2}{\rm d}x{\rm d}t
 \\\nonumber&\hspace{0.5cm}-C_{24}s^{3}\lambda^{\frac{7}{2}}\int_Q\varphi^{3}{\vartheta}^{2}{\rm d}x{\rm d}t
 -C_{24}s\lambda^{\frac{1}{2}}\int_Q\varphi\left(|\nabla{\vartheta}|^2+|\nabla{w}|^2\right){\rm d}x{\rm d}t
 \\\nonumber&\geq \sigma_1s\lambda\int_Q\varphi \left|\nabla{\vartheta}\right|^{2}{\rm d}x{\rm d}t
 -C_{25}\int_Q\left(|\nabla{\vartheta}|^2+|\nabla{w}|^2\right){\rm d}x{\rm d}t
 \\\nonumber&\hspace{0.5cm}-C_{25}s^2\lambda^{2}\int_Q\varphi^2{\vartheta}^{2}{\rm d}x{\rm d}t
 -C_{25}s\lambda^{\frac{7}{2}}\int_Q\varphi{\vartheta}^{2}{\rm d}x{\rm d}t
 \\\nonumber&\hspace{0.5cm}-C_{25}s\lambda^{\frac{1}{2}}\int_Q\varphi\left(|\nabla{\vartheta}|^2+|\nabla{w}|^2\right){\rm d}x{\rm d}t
 \\\nonumber&\hspace{0.5cm}-s\lambda\int_Q\varphi \left(\sqrt{\sigma_1}\left|\nabla{\vartheta}\right|\right)\left(\frac{M_0^2}{\sqrt{\sigma_1}}
 \left|\nabla{w}\right|\right){\rm d}x{\rm d}t
 \\\nonumber&\hspace{0.5cm}-C_{25}s^{3}\lambda^{3}\int_Q\varphi^{3}{w}^{2}{\rm d}x{\rm d}t
 -C_{25}s^{3}\lambda^{\frac{7}{2}}\int_Q\varphi^{3}{\vartheta}^{2}{\rm d}x{\rm d}t
 \\\nonumber&\geq\frac{1}{2}\sigma_{1}s\lambda\int_Q\varphi\left|\nabla{\vartheta}\right|^{2}{\rm d}x{\rm d}t
 -\frac{M^{4}_{0}}{2\sigma_{1}}s\lambda\int_Q\varphi\left|\nabla{w}\right|^{2}{\rm d}x{\rm d}t
 \\\nonumber&\hspace{0.5cm}-C_{26}s\lambda^{\frac{1}{2}}\int_Q\varphi\left(|\nabla{\vartheta}|^2+|\nabla{w}|^2\right){\rm d}x{\rm d}t
 -C_{26}s^{3}\lambda^{\frac{7}{2}}\int_Q\varphi^{3}{\vartheta}^{2}{\rm d}x{\rm d}t
 \\\nonumber&\hspace{0.5cm}-C_{26}s^{3}\lambda^{3}\int_Q\varphi^{3}{w}^{2}{\rm d}x{\rm d}t.
\end{align}
Since
\begin{align}
&\nonumber\int_Q\left|\mathcal{L}\left({\vartheta},{w}\right)\left(s\lambda\varphi{\vartheta}\right)\right|{\rm d}x{\rm d}t
\leq\frac{1}{2}\int_Q\left|\mathcal{L}\left({\vartheta},{w}\right)\right|^{2}{\rm d}x{\rm d}t
+\frac{1}{2}s^{2}\lambda^{2}\int_Q\varphi^{2}{\vartheta}^{2}{\rm d}x{\rm d}t,
\end{align}
 we obtain
 \begin{align}\label{equation2}
 &\nonumber \sigma_{1}\int_{Q}s\lambda\varphi\left|\nabla{\vartheta}\right|^{2}{\rm d}x{\rm d}t
 \\\nonumber&\leq \frac{M^{4}_{0}}{\sigma_{1}}\int_{Q}s\lambda\varphi\left|\nabla{w}\right|^{2}{\rm d}x{\rm d}t +C_{27}\int_Q\left|\mathcal{L}\left({\vartheta},{w}\right)\right|^{2}{\rm d}x{\rm d}t
 \\\nonumber&\hspace{0.5cm}
 +C_{27}s\lambda^{\frac{1}{2}}\int_Q\varphi\left(|\nabla{\vartheta}|^2+|\nabla{w}|^2\right){\rm d}x{\rm d}t
 +C_{27}s^{3}\lambda^{\frac{7}{2}}\int_Q\varphi^{3}{\vartheta}^{2}{\rm d}x{\rm d}t
 \\&\hspace{0.5cm}+C_{27}s^{3}\lambda^{3}\int_Q\varphi^{3}{w}^{2}{\rm d}x{\rm d}t.
 \end{align}

 Multiplying (\ref{equation2}) by $\frac{2\alpha_{1}}{\sigma_{1}}$, adding it to (\ref{equation1}), and noting the equality in (\ref{L12}), we obtain
 \begin{align}\label{equation12}
 &\nonumber 2\sigma^{2}_{1}s^{3}\lambda^{4}\int_{Q}\varphi^{3}\left|\nabla\psi\right|^{4}{\vartheta}^{2}{\rm d}x{\rm d}t
 +\frac{1}{60M}\int_{Q}\frac{1}{s\varphi}\left|\partial_{t}{\vartheta}\right|^{2}{\rm d}x{\rm d}t
 \\\nonumber&\hspace{0.5cm}+\alpha_{1}s\lambda\int_{Q}\varphi\left|\nabla{\vartheta}\right|^{2}{\rm d}x{\rm d}t
 +2 s\lambda^{2}\int_{Q}\varphi
 \left|a_{5}\left(\nabla\psi\cdot\nabla{\vartheta}\right)+a_{6}\left(\nabla\psi\cdot\nabla{w}\right)\right|^{2}{\rm d}x{\rm d}t
  \\\nonumber&\hspace{0.5cm}-2s\lambda\int_\Sigma \varphi\left(\nabla\psi\cdot\nu\right)\left(a_{5}
 \frac{\partial {\vartheta}}{\partial\nu}+a_{6}
 \frac{\partial{w}}{\partial \nu}\right)^2{\rm d}\sigma{\rm d}t +\left\|L_{1}\left({\vartheta},{w}\right)\right\|^{2}_{L^2(Q)}
 \\\nonumber&\leq C_{28}\int_{Q}\left(f^{2}+g^{2}\right){\rm e}^{2s\varphi}{\rm d}x{\rm d}t
 +\left(\frac{2M^{4}_{0}\alpha_{1}}{\sigma^{2}_{1}}
 +\alpha_{2}\right)s\lambda\int_{Q}\varphi\left|\nabla{{w}}\right|^{2}{\rm d}x{\rm d}t
 \\\nonumber&\hspace{0.2cm}+4M s\lambda\int_{Q}\varphi
 \left|\partial_{t}{{w}}\right|^{2}{\rm d}x{\rm d}t
 +\frac{32M^{8}_{0}}{\sigma^{2}_{1}}
 s^{3}\lambda^{4}\int_{Q}\varphi^{3}\left|\nabla\psi\right|^{4}{w}^{2}{\rm d}x{\rm d}t
 \\\nonumber&\hspace{0.2cm}+6M_0^2s\lambda\beta T\int_{Q}\varphi\left|\nabla{\vartheta}\right|^{2}{\rm d}x{\rm d}t
 +2M_0^2s\lambda\beta T\int_{Q}\varphi\left|\nabla{w}\right|^{2}{\rm d}x{\rm d}t
 \\\nonumber&\hspace{0.2cm}+\frac{4M^{3}_{0}}{\sigma_{1}}s\lambda\beta T\int_{Q}\varphi
 \left|\partial_{t}{w}\right|^{2}{\rm d}x{\rm d}t
 +C_{28}s\lambda^{\frac{1}{2}}\int_{Q}\varphi\left(\left|\nabla{{\vartheta}}\right|^{2}
 +\left|\nabla{{w}}\right|^{2}
  +\left|\partial_{t}{{w}}\right|^{2}\right){\rm d}x{\rm d}t
 \\ &\hspace{0.2cm}+C_{28}\int_{Q}\left(s^{3}\lambda^{\frac{7}{2}}\varphi^{3}+s^{2}\lambda^{6}\varphi^{2}\right)
 \left({\vartheta}^{2}+{w}^{2}\right){\rm d}x{\rm d}t
 +C_{28}\int_{Q}\frac{1}{s\lambda\varphi}
 \left|\partial_{t}{{\vartheta}}\right|^{2}{\rm d}x{\rm d}t,
 \end{align}
 for all $s>s_1$.

 Next we estimate $s\lambda\int_{Q}\varphi\left|\nabla_{x,t}w\right|^{2}{\rm d}x{\rm d}t$ and $s^3\lambda^4\int_{Q}\varphi^3 w^{2}{\rm d}x{\rm d}t$. By (\ref{fg}), we can get the following equation
 \begin{equation}\label{a_78}
 a_7\partial^{2}_{t}p-a_{1}\Delta p-a_8\partial_{t}\Theta=f-a_8g,
 \end{equation}
 where
 \begin{equation}\label{a78}
 a_7= 1+\frac{a_{2}a_{4}}{a_{3}}\qquad
 \mbox{and}\quad a_8=\frac{a_{2}}{a_{3}}.
 \end{equation}
Denote
\begin{equation}\label{l3}
{L}_{3}(\Theta, p)={a_{7}}\partial^{2}_{t}p-a_{1}\Delta p-a_8\partial_{t}\Theta.
\end{equation}
 We recall $\left({\vartheta},{w}\right)={\rm e}^{s\varphi}\left(\Theta,p\right)$. We set
 \begin{align}\label{L_ast}
 \mathcal{L}_{\ast}\left({\vartheta},{w}\right)&={\rm e}^{s\varphi}{L}_{3}\left({\rm e}^{-s\varphi}({\vartheta},{w})\right).
 \end{align} Then we have
 \begin{align*}
 &\nonumber\left(f-{a_{8}}g\right){\rm e}^{s\varphi}=\mathcal{L}_{\ast}\left({\vartheta},{w}\right)={\rm e}^{s\varphi}\left\{{a_{7}}\partial^{2}_{t}\left({\rm e}^{-s\varphi}{w}\right)-a_{1}\Delta\left({\rm e}^{-s\varphi}{w}\right)-a_8\partial_t\left({\rm e}^{-s\varphi}{\vartheta}\right)\right\}
 \\&=\left\{s^{2}\lambda^{2}\varphi^{2}\left(\partial_{t}\psi\right)^{2}
 -s\lambda^{2}\varphi\left(\partial_{t}\psi\right)^{2}-s\lambda\varphi\partial^{2}_{t}\psi\right\}{a_{7}}{w}
 +{a_{7}}\partial^{2}_{t}{w}
 \\&\hspace{0.5cm}-2s\lambda\varphi\left(\partial_{t}\psi\right){a_{7}}\partial_{t}{w}
 -\left(s^{2}\lambda^{2}\varphi^{2}\left|\nabla\psi\right|^{2}
 -s\lambda^{2}\varphi\left|\nabla\psi\right|^{2}-2ns\lambda\varphi\right)a_{1}{w}
 \\&\hspace{0.5cm}-a_{1}\Delta{w}+2s\lambda\varphi a_{1}\left(\nabla\psi\cdot\nabla{w}\right)
 -{a_{8}}\partial_{t}{\vartheta}+s\lambda\varphi\left(\partial_{t}\psi\right)a_8{\vartheta}.
 \end{align*}
 We set
 \begin{align}\label{j0}
 J_{0}=s\lambda\varphi\gamma_{0}{w}+\sum_{j=1}^{n}\gamma_{j}\partial_{j}{w},
 \end{align}
 where $\gamma_{0}=\gamma_{0}(x)$ and $\gamma_{j}=\gamma_{j}(x)\in  C^2\left(\overline{\Omega}\right)$, ($j=1$, $\cdots$, $n$) are functions of $x$, but are independent of $t$, $s$, $\lambda$,
 $\varphi$, and will be suitably chosen later.
 We divide $\mathcal{L}_{\ast}\left({\vartheta},{w}\right)$ as follows:
 \begin{align}\label{last}
  \mathcal{L}_{\ast}\left({\vartheta},{w}\right)=L_{1\ast}\left({\vartheta},{w}\right)+L_{2\ast}\left({w}\right),
 \end{align}
 where
 \begin{align}\label{l1ast}
 \nonumber L_{1\ast}\left({\vartheta},{w}\right)&={a_{7}}\partial^{2}_{t}{w}-a_{1}\Delta{w}
 +s^{2}\lambda^{2}\varphi^{2}\left(\partial_{t}\psi\right)^{2}{a_{7}}{w}
 \\&\hspace{0.5cm}-s^{2}\lambda^{2}\varphi^{2}\left|\nabla\psi\right|^{2}a_{1}{w}
 -{a_{8}}\partial_{t}{\vartheta}+s\lambda\varphi\left(\partial_{t}\psi\right)a_8{\vartheta}+J_{0},
 \end{align}
 \begin{align}\label{l2ast}
 \nonumber L_{2\ast}\left({w}\right)&=-2s\lambda\varphi
 \left(\partial_{t}\psi\right){a_{7}}\partial_{t}{w}
 +2s\lambda\varphi a_{1}\left(\nabla\psi\cdot\nabla{w}\right)
 \\\nonumber&\hspace{0.5cm}-s\lambda^{2}\varphi\left(\partial_{t}\psi\right)^{2}{a_{7}}{w}
 -s\lambda\varphi\left(\partial^{2}_{t}\psi\right){a_{7}}{w}
 \\&\hspace{0.5cm}+s\lambda^{2}\varphi\left|\nabla\psi\right|^{2}a_{1}{w}+2ns\lambda\varphi a_{1}{w}-J_{0}.
 \end{align}
 Then we have
 \begin{align}\label{eq}
 \nonumber&\left\|\left(f-{a_{8}}g\right){\rm e}^{s\varphi}\right\|_{L^2(Q)}^2=\left\|\mathcal{L}_{\ast}\left({\vartheta},{w}\right)\right\|_{L^2(Q)}^2
 \\&=\left\|\mathcal{L}_{1\ast}\left({\vartheta},{w}\right)\right\|_{L^2(Q)}^2+2\int_QL_{1\ast}\left({\vartheta},{w}\right)L_{2\ast}\left({w}\right){\rm d}x{\rm d}t+\left\|\mathcal{L}_{2\ast}\left({w}\right)\right\|_{L^2(Q)}^2.
 \end{align}
 We have
 \begin{align}\label{ik}
 \nonumber L_{1\ast}&\left({\vartheta}, {w}\right)L_{2\ast}\left({w}\right)
 ={a_{7}}\left(\partial^{2}_{t}{w}\right)L_{2\ast}\left({w}\right)
 -a_{1}\left(\Delta{w}\right)L_{2\ast}\left({w}\right)
 \\\nonumber&\hspace{0.5cm}+s^{2}\lambda^{2}\varphi^{2}\left(\partial_{t}\psi\right)^{2}{a_{7}}{w}
 L_{2\ast}\left({w}\right)
 -s^{2}\lambda^{2}\varphi^{2}\left|\nabla\psi\right|^{2}a_{1}{w}L_{2\ast}\left({w}\right)
 \\&\hspace{0.5cm}-{a_{8}}\left(\partial_{t}{\vartheta}\right)L_{2\ast}\left({w}\right)
 +s\lambda\varphi\left(\partial_{t}\psi\right)a_8{\vartheta}L_{2\ast}\left({w}\right)
 +J_{0}L_{2\ast}\left({w}\right)
 \triangleq\sum_{k=1}^{7}I_{k}.
\end{align}
 We integrate $I_{k}$, $k=1,2,...,6$ over $Q$. We calculate them by applying the integration by parts and $\big(a_1, a_2, a_3, a_4\big)\in \mathcal{U}$. By (\ref{Tb2}) and
 $\partial_tw=\partial_t\left({\rm e}^{s\varphi}{p}\right)={\rm e}^{s\varphi}\partial_t{p}+s\left(\partial_t\varphi\right)p$, we have
  $$
 {\vartheta}(x, 0)={\vartheta}(x,T)=0, \  \partial_t^j{w}(x,0)=\partial_t^j{w}(x,T)=0,\  x\in\overline{\Omega}, \ j=0,1.
 $$
 We further recall that, by (\ref{Db2}), we have ${\vartheta}=0$ and ${w}=0$ on $\Sigma$, so that we have
 $\nabla{\vartheta}=\frac{\partial {\vartheta}}{\partial\nu}\nu$ and $\nabla{w}=\frac{\partial {w}}{\partial\nu}\nu$ on $\Sigma$.
 \begin{align*}
 \nonumber &\int_QI_{1}{\rm d}x{\rm d}t=\int_Q{a_{7}}\left(\partial^{2}_{t}{w}\right)
 L_{2\ast}\left({w}\right){\rm d}x{\rm d}t
 =-\int_Q{a_{7}}\left(\partial_{t}{w}\right)\partial_{t}
 \left\{L_{2\ast}\left({w}\right)\right\}{\rm d}x{\rm d}t
 \\\nonumber&=2 s\lambda^{2}\int_Q\varphi
 \left(\partial_{t}\psi\right)^{2}{a_{7}^{2}}\left|\partial_{t}{w}\right|^{2}{\rm d}x{\rm d}t
 +2s\lambda\int_Q\varphi
 \left(\partial^{2}_{t}\psi\right){a_{7}^{2}}\left|\partial_{t}{w}\right|^{2}{\rm d}x{\rm d}t
 \\\nonumber&\hspace{0.2cm}+
 s\lambda\int_Q\varphi\left(\partial_{t}\psi\right){a_{7}^{2}}\partial_{t}\left\{\left(\partial_{t}{w}\right)^2\right\}{\rm d}x{\rm d}t
 \\\nonumber&\hspace{0.2cm}-2s\lambda^{2}\int_Q\varphi\left(\partial_{t}\psi\right)a_{1}{a_{7}}
 \left(\nabla\psi\cdot\nabla{w}\right)\left(\partial_{t}{w}\right){\rm d}x{\rm d}t
 \\\nonumber&\hspace{0.2cm}-s\lambda\int_Q\varphi a_{1}{a_{7}}\left(\nabla\psi\cdot\nabla\left\{\left(\partial_t{w}\right)^2\right\}\right)
 {\rm d}x{\rm d}t
 \\\nonumber&\hspace{0.2cm}+\int_Q\Big\{s\lambda^{2}\varphi
 \left(\partial_{t}\psi\right)^{2}{a_{7}^{2}}\left|\partial_{t}{w}\right|^{2}
 +s\lambda\varphi
 \left(\partial^{2}_{t}\psi\right){a_{7}^{2}}\left|\partial_{t}{w}\right|^{2}
 \\\nonumber&\hspace{0.7cm}-s\lambda^{2}\varphi\left|\nabla\psi\right|^{2}a_{1}{a_{7}}\left|\partial_{t}{w}\right|^{2}
 -2ns\lambda\varphi a_{1}{a_{7}}\left|\partial_{t}{w}\right|^{2}
 +s\lambda\varphi{a_{7}}\gamma_{0}\left|\partial_{t}{w}\right|^{2}\Big\}{\rm d}x{\rm d}t
 \\\nonumber&\hspace{0.2cm}+\frac{1}{2}\int_Q{a_{7}}\sum_{j=1}^{n}\gamma_{j}\partial_j
 \left\{\left(\partial_{t}{w}\right)^2\right\}
 {\rm d}x{\rm d}t
  \\\nonumber&\hspace{0.2cm}+\int_Q\Big\{s\lambda^{3}\varphi
 \left(\partial_{t}\psi\right)^{3}{a_{7}^{2}}\left(\partial_{t}{w}\right){w}
 +2s\lambda^{2}\varphi
 \left(\partial_{t}\psi\right)\left(\partial^{2}_{t}\psi\right){a_{7}^{2}}\left(\partial_{t}{w}\right){w}
 \\\nonumber&\hspace{0.7cm}
 +s\lambda^{2}\varphi
 \left(\partial_{t}\psi\right)\left(\partial^{2}_{t}\psi\right){a_{7}^{2}}\left(\partial_{t}{w}\right){w}-
 s\lambda^{3}\varphi\left(\partial_{t}\psi\right)\left|\nabla\psi\right|^{2}a_{1}{a_{7}}\left(\partial_{t}{w}\right){w}
 \\&\hspace{0.7cm}
 -2ns\lambda^{2}\varphi\left(\partial_{t}\psi\right)a_{1}{a_{7}}
 \left(\partial_{t}{w}\right){w}
 +s\lambda^{2}\varphi\left(\partial_{t}\psi\right){a_{7}}\gamma_{0}
 \left(\partial_{t}{w}\right){w}
 \Big\}{\rm d}x{\rm d}t.
 \end{align*}
 We integrate the third and  fifth terms by parts and then combine them with the sixth term.
 We integrate the seventh term by parts.
 We further estimate the eighth term by
 the Cauchy-Schwarz inequality. Then we obtain
 \begin{align}\label{i1}
 \nonumber &\int_QI_{1}{\rm d}x{\rm d}t
 \geq
 2 s\lambda^{2}\int_Q\varphi
 \left(\partial_{t}\psi\right)^{2}{a_{7}^{2}}\left|\partial_{t}{w}\right|^{2}{\rm d}x{\rm d}t
 +2s\lambda\int_Q\varphi
 \left(\partial^{2}_{t}\psi\right){a_{7}^{2}}\left|\partial_{t}{w}\right|^{2}{\rm d}x{\rm d}t
 \\\nonumber&\hspace{0.2cm}-2s\lambda^{2}\int_Q\varphi\left(\partial_{t}\psi\right)a_{1}{a_{7}}
 \left(\nabla\psi\cdot\nabla{w}\right)\left(\partial_{t}{w}\right){\rm d}x{\rm d}t
 \\\nonumber&\hspace{0.2cm}
 +s\lambda\int_Q\varphi \left(\nabla\psi\cdot\nabla\left(a_{1}{a_{7}}\right)\right)\left|\partial_t{w}\right|^2
 {\rm d}x{\rm d}t
 +s\lambda\int_Q \varphi{a_{7}}\gamma_{0}\left|\partial_{t}{w}\right|^{2}{\rm d}x{\rm d}t
 \\\nonumber &\hspace{0.2cm}
 -\frac{1}{2}\int_Q\left\{\sum_{j=1}^{n}\partial_j\left({a_{7}}\gamma_{j}\right)\right\}
 \left|\partial_{t}{w}\right|^2
 {\rm d}x{\rm d}t
  \\&\hspace{0.2cm}-C_{29}\int_Q\left|\partial_{t}{w}\right|^{2}{\rm d}x{\rm d}t
 -C_{29}s^2\lambda^6\int_Q \varphi^2{w}^{2}{\rm d}x{\rm d}t.
 \end{align}
 \begin{align*}
 \nonumber &\int_QI_{2}{\rm d}x{\rm d}t=-\int_Qa_{1}\left(\Delta{w}\right)L_{2\ast}\left({w}\right){\rm d}x{\rm d}t
 \\\nonumber &=-2s\lambda\int_\Sigma\varphi(\nabla\psi\cdot\nu)a_{1}^2\left|\frac{\partial{w}}{\partial\nu}\right|^2{\rm d}\sigma{\rm d}t
 +\int_Q a_{1}\left(\nabla {w}\cdot \nabla\left\{L_{2\ast}\left({w}\right)\right\}\right){\rm d}x{\rm d}t
 \\\nonumber&\hspace{0.2cm}+\int_Q\left(\nabla a_{1}\cdot\nabla{w}\right)L_{2\ast}\left({w}\right){\rm d}x{\rm d}t
 \\\nonumber &\geq -2s\lambda\int_\Sigma\varphi(\nabla\psi\cdot\nu)a_{1}^2\left|\frac{\partial{w}}{\partial\nu}\right|^2{\rm d}\sigma{\rm d}t
 -2s\lambda^2\int_Q\varphi\left(\partial_t\psi\right)a_1a_7\left(\nabla\psi\cdot\nabla {w}\right)\left(\partial_t {w}\right) {\rm d}x{\rm d}t
 \\\nonumber&\hspace{0.2cm}-2s\lambda\int_Q\varphi\left(\partial_t\psi\right)a_1\left(\nabla a_7\cdot\nabla {w}\right)\left(\partial_t{w}\right) {\rm d}x{\rm d}t
 -s\lambda\int_Q\varphi\left(\partial_t\psi\right)a_1a_7\partial_t \left(\left|\nabla {w}\right|^2\right) {\rm d}x{\rm d}t
 \\\nonumber&\hspace{0.2cm} +2s\lambda^2\int_Q\varphi a_1^2\left(\nabla\psi\cdot\nabla{w}\right)^2{\rm d}x{\rm d}t
 +2s\lambda\int_Q\varphi a_1\left(\nabla a_1\cdot\nabla{w}\right)\left(\nabla\psi\cdot\nabla{w}\right){\rm d}x{\rm d}t
 \\\nonumber&\hspace{0.2cm}+4s\lambda\int_Q\varphi a_1^2\left|\nabla{w}\right|^2{\rm d}x{\rm d}t
 +s\lambda\int_Q\varphi a_1^2\left(\nabla\psi\cdot\nabla\left\{\left|\nabla{w}\right|^2\right\}\right)
 {\rm d}x{\rm d}t
 \\\nonumber&\hspace{0.2cm}-s\lambda^2\int_Q\varphi\left(\partial_t\psi\right)^2a_1a_7\left|\nabla{w}\right|^2{\rm d}x{\rm d}t
 -s\lambda\int_Q\varphi\left(\partial_t^2\psi\right)a_1a_7\left|\nabla{w}\right|^2{\rm d}x{\rm d}t
 \\\nonumber&\hspace{0.2cm}+s\lambda^2\int_Q\varphi\left|\nabla\psi\right|^2a_1^2\left|\nabla{w}\right|^2{\rm d}x{\rm d}t
 +2ns\lambda\int_Q\varphi a_1^2\left|\nabla{w}\right|^2{\rm d}x{\rm d}t
 \\\nonumber&\hspace{0.2cm}-s\lambda\int_Q\varphi \gamma_0a_1\left|\nabla{w}\right|^2{\rm d}x{\rm d}t
 -\frac{1}{2}\int_Qa_1\sum_{j=1}^n\gamma_j\partial_j\left(\left|\nabla{w}\right|^2\right){\rm d}x{\rm d}t
 \\\nonumber&\hspace{0.2cm}-2s\lambda\int_Q\varphi\left(\partial_t\psi\right)a_7\left(\nabla a_1\cdot\nabla {w}\right)\left(\partial_t{w}\right){\rm d}x{\rm d}t
 \\\nonumber&\hspace{0.2cm}+2s\lambda\int_Q\varphi a_1\left(\nabla a_1\cdot\nabla {w}\right)\left(\nabla \psi\cdot\nabla {w}\right){\rm d}x{\rm d}t
 \\&\hspace{0.2cm}-C_{30}\int_Q\left|\nabla{w}\right|^2{\rm d}x{\rm d}t-C_{30}s^2\lambda^6\int_Q\varphi^2\left|{w}\right|^2{\rm d}x{\rm d}t.
 \end{align*}
 We combine the third and fifteenth terms.  We integrate the fourth term by parts and then combine them with the ninth
 and tenth terms. We combine the sixth and sixteenth terms.
 We integrate the eighth term by parts and then combine them with the first, eleventh and  twelfth terms.
 We further integrate the fourteenth term  by parts. Thus we obtain
 \begin{align}\label{i2}
 \nonumber &\int_QI_{2}{\rm d}x{\rm d}t\geq
 -s\lambda\int_\Sigma\varphi(\nabla\psi\cdot\nu)a_{1}^2\left|\frac{\partial{w}}{\partial\nu}\right|^2{\rm d}\sigma{\rm d}t
 -\frac{1}{2}\int_\Sigma a_1\left(\sum_{j=1}^n\gamma_j\nu_j\right)\left|\frac{\partial{w}}{\partial\nu}\right|^2{\rm d}\sigma{\rm d}t
 \\\nonumber&\hspace{0.2cm}-2s\lambda^2\int_Q\varphi\left(\partial_t\psi\right)a_1a_7\left(\nabla\psi\cdot\nabla {w}\right)\left(\partial_t {w}\right) {\rm d}x{\rm d}t
 \\\nonumber&\hspace{0.2cm}-2s\lambda\int_Q\varphi\left(\partial_t\psi\right)
    \left(\nabla \left(a_1 a_7\right)\cdot\nabla {w}\right)\left(\partial_t{w}\right) {\rm d}x{\rm d}t
 \\\nonumber&\hspace{0.2cm} +2s\lambda^2\int_Q\varphi a_1^2\left(\nabla\psi\cdot\nabla{w}\right)^2{\rm d}x{\rm d}t
 +4s\lambda\int_Q\varphi a_1\left(\nabla a_1\cdot\nabla{w}\right)\left(\nabla\psi\cdot\nabla{w}\right){\rm d}x{\rm d}t
 \\\nonumber&\hspace{0.2cm}+4s\lambda\int_Q\varphi a_1^2\left|\nabla{w}\right|^2{\rm d}x{\rm d}t
 -2s\lambda\int_Q\varphi a_1\left(\nabla\psi\cdot\nabla a_1\right)\left|\nabla{w}\right|^2{\rm d}x{\rm d}t
 \\\nonumber&\hspace{0.2cm}-s\lambda\int_Q\varphi \gamma_0a_1\left|\nabla{w}\right|^2{\rm d}x{\rm d}t
 +\frac{1}{2}\int_Q\left\{\sum_{j}^n\partial_j\left(a_1\gamma_j\right)\right\}\left|\nabla{w}\right|^2{\rm d}x{\rm d}t
 \\&\hspace{0.2cm}-C_{31}\int_Q\left|\nabla{w}\right|^2{\rm d}x{\rm d}t-C_{31}s^2\lambda^6\int_Q\varphi^2\left|{w}\right|^2{\rm d}x{\rm d}t.
 \end{align}
 \begin{align*}
 \nonumber \int_Q&I_{3}{\rm d}x{\rm d}t=s^{2}\lambda^{2}\int_Q\varphi^{2}\left(\partial_{t}\psi\right)^{2}{a_{7}}
 {w}L_{2\ast}\left({w}\right){\rm d}x{\rm d}t
 \\\nonumber&\geq- s^{3}\lambda^{3}\int_Q\varphi^{3}\left(\partial_{t}\psi\right)^{3}
 {a_{7}^{2}}\left\{\partial_{t}\left({w}^2\right)\right\}{\rm d}x{\rm d}t
 \\\nonumber&\hspace{0.2cm}+s^{3}\lambda^{3}\int_Q\varphi^{3}
 \left(\partial_{t}\psi\right)^{2}a_{1}{a_{7}}\left(\nabla\psi\cdot\nabla\left({w}^2\right)\right){\rm d}x{\rm d}t
 \\\nonumber&\hspace{0.2cm}-s^{3}\lambda^{4}\int_Q\varphi^{3}\left(\partial_{t}\psi\right)^{4}{a_{7}^{2}}
 {w}^{2}{\rm d}x{\rm d}t+s^{3}\lambda^{4}\int_Q\varphi^{3}\left(\partial_{t}\psi\right)^{2}
 \left|\nabla\psi\right|^{2}a_{1}{a_{7}}{w}^{2}{\rm d}x{\rm d}t
 \\\nonumber&\hspace{0.2cm}-C_{32}s^{3}\lambda^{3}\int_Q\varphi^{3}
 {w}^{2}{\rm d}x{\rm d}t-\frac{1}{2}s^{2}\lambda^{2}\int_Q\varphi^{2}
 \left(\partial_{t}\psi\right)^{2}{a_{7}}\sum_{j=1}^{n}\gamma_{j}\left\{\partial_{j}\left({w}^2\right)\right\}{\rm d}x{\rm d}t.
 \end{align*}
 We integrate the first term by parts and then combine it with the third term.
 We integrate the second term by parts and then combine it with the fourth term.
 We further integrate the sixth term by parts and estimate it. Therefore we can obtain
 \begin{align*}
 \nonumber \int_Q&I_{3}{\rm d}x{\rm d}t\geq
2s^{3}\lambda^{4}\int_Q\varphi^{3}\left(\partial_{t}\psi\right)^{4}{a_{7}^{2}}
 {w}^{2}{\rm d}x{\rm d}t+3s^{3}\lambda^{3}\int_Q\varphi^{3}\left(\partial_{t}\psi\right)^{2}
 \left(\partial_{t}^2\psi\right){a_{7}^{2}}
 {w}^{2}{\rm d}x{\rm d}t
 \\\nonumber&\hspace{0.2cm}-2s^{3}\lambda^{4}\int_Q\varphi^{3}\left(\partial_{t}\psi\right)^{2}
 \left|\nabla\psi\right|^{2}a_{1}{a_{7}}{w}^{2}{\rm d}x{\rm d}t
 \\\nonumber&\hspace{0.2cm}-s^{3}\lambda^{3}\int_Q\varphi^{3}\left(\partial_{t}\psi\right)^{2}\left(\nabla\psi\cdot\nabla\left(a_1a_7\right)\right)
 {w}^{2}{\rm d}x{\rm d}t
 \\\nonumber&\hspace{0.2cm}-C_{33}s^{3}\lambda^{3}\int_Q\varphi^{3}{w}^{2}{\rm d}x{\rm d}t
 +s^{2}\lambda^{3}\int_Q\varphi^{2}
 \left(\partial_{t}\psi\right)^{2}{a_{7}}\left\{\sum_{j=1}^{n}\gamma_{j}\left(\partial_j\psi\right)\right\}{w}^2{\rm d}x{\rm d}t
 \\\nonumber&\hspace{0.2cm}+\frac{1}{2}s^{2}\lambda^{2}\int_Q\varphi^{2}
 \left(\partial_{t}\psi\right)^{2}\left\{\sum_{j=1}^{n}\partial_{j}\left({a_{7}}\gamma_{j}\right)\right\}{w}^2{\rm d}x{\rm d}t.
  \end{align*}
  Noting the first term is nonnegative, we obtain
 \begin{align}\label{i3}
  \int_Q&I_{3}{\rm d}x{\rm d}t\geq
 -2s^{3}\lambda^{4}\int_Q\varphi^{3}\left(\partial_{t}\psi\right)^{2}
 \left|\nabla\psi\right|^{2}a_{1}{a_{7}}{w}^{2}{\rm d}x{\rm d}t
 -C_{34}s^{3}\lambda^{3}\int_Q\varphi^{3}{w}^{2}{\rm d}x{\rm d}t.
 \end{align}
 \begin{align*}
 \nonumber &\int_QI_{4}{\rm d}x{\rm d}t=-s^{2}\lambda^{2}\int_Q\varphi^{2}\left|\nabla\psi\right|^{2}a_{1}
 {w}L_{2\ast}\left({w}\right){\rm d}x{\rm d}t
 \\\nonumber&=s^{3}\lambda^{3}\int_Q\varphi^{3}\left|\nabla\psi\right|^{2}
 \left(\partial_{t}\psi\right)a_{1}{a_{7}}\left\{\partial_{t}\left({w}^2\right)\right\}{\rm d}x{\rm d}t
 \\\nonumber&\hspace{0.2cm}- s^{3}\lambda^{3}\int_Q\varphi^{3}\left|\nabla\psi\right|^{2}a_{1}^{2}
 \left(\nabla\psi\cdot\nabla\left({w}^2\right)\right){\rm d}x{\rm d}t
 \\\nonumber&\hspace{0.2cm}+s^{3}\lambda^{4}\int_Q\varphi^{3}\left|\nabla\psi\right|^{2}
 \left(\partial_{t}\psi\right)^{2}a_{1}{a_{7}}{w}^{2}{\rm d}x{\rm d}t
 +s^{3}\lambda^{3}\int_Q\varphi^{3}
 \left|\nabla\psi\right|^{2}\left(\partial^{2}_{t}\psi\right)a_{1}{a_{7}}{w}^{2}{\rm d}x{\rm d}t
 \\\nonumber&\hspace{0.2cm}-s^{3}\lambda^{4}\int_Q\varphi^{3}\left|\nabla\psi\right|^{4}a_{1}^{2}{w}^{2}{\rm d}x{\rm d}t
 -2n s^{3}\lambda^{3}\int_Q\varphi^{3}
 \left|\nabla\psi\right|^{2}a_{1}^{2}{w}^{2}{\rm d}x{\rm d}t
 \\\nonumber&\hspace{0.2cm}+s^{3}\lambda^{3}\int_Q\varphi^{3}\left|\nabla\psi\right|^{2}a_{1}\gamma_{0}{w}^{2}{\rm d}x{\rm d}t
 +\frac{1}{2}s^{2}\lambda^{2}\int_Q\varphi^{2}\left|\nabla\psi\right|^{2}a_{1}\sum_{j=1}^{n}\gamma_{j}\left\{\partial_{j}\left({w}^2\right)
 \right\} {\rm d}x{\rm d}t.
 \end{align*}
 We integrate the first term by parts and combine it with the third and fourth term.
 We integrate the second term by parts and combine it with the fifth and sixth term.
 We integrate the eighth term by parts.
 Then we obtain
 \begin{align}\label{i4}
 \nonumber &\int_QI_{4}{\rm d}x{\rm d}t\geq
 -2s^{3}\lambda^{4}\int_Q\varphi^{3}\left|\nabla\psi\right|^{2}
 \left(\partial_{t}\psi\right)^{2}a_{1}{a_{7}}{w}^{2}{\rm d}x{\rm d}t
 \\\nonumber&\hspace{0.2cm}+2s^{3}\lambda^{4}\int_Q\varphi^{3}\left|\nabla\psi\right|^{4}a_{1}^{2}{w}^{2}{\rm d}x{\rm d}t
 +s^{3}\lambda^{3}\int_Q\varphi^{3}\left(\nabla\psi\cdot\nabla\left(\left|\nabla\psi\right|^{2}a_1^2\right)\right)
 {w}^{2}{\rm d}x{\rm d}t
 \\\nonumber&\hspace{0.2cm}+s^{3}\lambda^{3}\int_Q\varphi^{3}\left|\nabla\psi\right|^{2}a_{1}\gamma_{0}{w}^{2}{\rm d}x{\rm d}t
 -s^{2}\lambda^{3}\int_Q\varphi^{2}\left|\nabla\psi\right|^{2}a_{1}
 \left(\sum_{j=1}^{n}\gamma_{j}\partial_{j}\psi\right)
 {w}^2 {\rm d}x{\rm d}t
 \\\nonumber&\hspace{0.2cm}-\frac{1}{2}s^{2}\lambda^{2}\int_Q\varphi^{2}\left\{\sum_{j=1}^{n}
 \partial_{j}\left(\left|\nabla\psi\right|^{2}a_{1}\gamma_{j}\right)\right\}
 {w}^2 {\rm d}x{\rm d}t
 \\\nonumber&\geq
 -2s^{3}\lambda^{4}\int_Q\varphi^{3}\left|\nabla\psi\right|^{2}
 \left(\partial_{t}\psi\right)^{2}a_{1}{a_{7}}{w}^{2}{\rm d}x{\rm d}t
 +2s^{3}\lambda^{4}\int_Q\varphi^{3}\left|\nabla\psi\right|^{4}a_{1}^{2}{w}^{2}{\rm d}x{\rm d}t
 \\&\hspace{0.5cm}-C_{35}s^{3}\lambda^{3}\int_Q\varphi^{3}{w}^{2}{\rm d}x{\rm d}t.
 \end{align}
 \begin{align}\label{i50}
 \nonumber I_{5}&=-{a_{8}}\left(\partial_{t}{\vartheta}\right)
 L_{2\ast}\left({w}\right)
 \\ \nonumber &=2 s\lambda\varphi\left(\partial_{t}\psi\right){a_{7}}{a_{8}}
 \left(\partial_{t}{\vartheta}\right)\left(\partial_{t}{w}\right)
 -2 s\lambda\varphi{a_{1}a_{8}}\left(\nabla\psi\cdot\nabla{w}\right)
 \left(\partial_{t}{\vartheta}\right)
 \\ \nonumber &\hspace{0.2cm}+s\lambda^{2}\varphi\left(\partial_{t}\psi\right)^{2}{a_{7}}{a_{8}}
 \left(\partial_{t}{\vartheta}\right){w}
+\big\{s\lambda\varphi\left(\partial^{2}_{t}\psi\right){a_{7}}{a_{8}}
  \left(\partial_{t}{\vartheta}\right){w}
 \\ \nonumber &\hspace{0.2cm}-s\lambda^{2}\varphi\left|\nabla\psi\right|^{2}{a_{1}a_{8}}\left(\partial_{t}{\vartheta}\right){w}
 -2n s\lambda\varphi{a_{1}a_{8}}\left(\partial_{t}{\vartheta}\right){w}
 \\&\hspace{0.2cm}+s\lambda\varphi{a_{8}}\gamma_{0}\left(\partial_{t}{\vartheta}\right){w}\big\}
 +{a_{8}}\left(\partial_{t}{\vartheta}\right) \sum_{j=1}^{n}\gamma_{j}\partial_{j}{w}\triangleq\sum_{k=1}^{5}I_{5k}.
 \end{align}
 Integrating by parts and using the Cauchy-Schwarz inequality, we have
 \begin{align*}
 \nonumber &\int_QI_{51}{\rm d}x{\rm d}t=2 s\lambda\int_Q\varphi\left(\partial_{t}\psi\right){a_{7}}{a_{8}}
 \left(\partial_{t}{\vartheta}\right)\left(\partial_{t}{w}\right) {\rm d}x{\rm d}t
 \\\nonumber&=-2 s\lambda\int_Q\varphi\left(\partial_{t}\psi\right){a_{7}}{a_{8}}
 {\vartheta}\left(\partial^{2}_{t}{w}\right) {\rm d}x{\rm d}t
 \\\nonumber&\hspace{0.5cm}-2 s\lambda^{2}\int_Q\varphi\left(\partial_{t}\psi\right)^{2}{a_{7}}{a_{8}}
 {\vartheta}\left(\partial_{t}{w}\right) {\rm d}x{\rm d}t
 -2 s\lambda\int_Q\varphi\left(\partial^{2}_{t}\psi\right){a_{7}}{a_{8}}
 {\vartheta}\left(\partial_{t}{w}\right) {\rm d}x{\rm d}t
 \\\nonumber&\geq -2 s\lambda\int_Q\varphi\left(\partial_{t}\psi\right){a_{7}}{a_{8}}
 {\vartheta}\left(\partial^{2}_{t}{w}\right){\rm d}x{\rm d}t
 \\&\hspace{0.5cm}-C_{36}s^{2}\lambda^{4}\int_Q\varphi^{2}{\vartheta}^{2}
 {\rm d}x{\rm d}t-C_{36}\int_Q\left|\partial_{t}{w}\right|^{2}{\rm d}x{\rm d}t.
 \end{align*}
 We calculate the first term as follows.
 By (\ref{l1ast}), we have
 \begin{align*}
 a_7\partial^{2}_{t}{w}
 =&L_{1\ast}\left({\vartheta},{w}\right)+a_{1}\Delta{w}
 -s^{2}\lambda^{2}\varphi^{2}\left(\partial_{t}\psi\right)^{2}{a_{7}}{w}
 \\&+s^{2}\lambda^{2}\varphi^{2}\left|\nabla\psi\right|^{2}a_{1}{w}
 +{a_{8}}\partial_{t}{\vartheta}-s\lambda\varphi\left(\partial_{t}\psi\right)a_8{\vartheta}-J_{0},\ \mbox{in} \ Q.
 \end{align*}
 Therefore we have
 \begin{align*}
 \nonumber &\int_QI_{51}{\rm d}x{\rm d}t
 \geq -2 s\lambda\int_Q\varphi\left(\partial_{t}\psi\right){a_{1}}{a_{8}}
 \vartheta\left(\triangle{w}\right){\rm d}x{\rm d}t
 \\ \nonumber&\hspace{0.5cm}+2 s^3\lambda^3\int_Q\varphi^3\left(\partial_{t}\psi\right){a_{8}}
 \left\{\left(\partial_{t}\psi\right)^{2}{a_{7}}
 -\left|\nabla\psi\right|^{2}a_{1}
 \right\}\vartheta{w}{\rm d}x{\rm d}t
 \\ \nonumber&\hspace{0.5cm}- s\lambda\int_Q\varphi\left(\partial_{t}\psi\right){a_{8}^2}
 \left\{\partial_{t}\left({\vartheta}^2\right)\right\}{\rm d}x{\rm d}t
 +2 s^2\lambda^2\int_Q\varphi^2\left(\partial_{t}\psi\right)^2{a_{8}^2}
 {\vartheta}^2{\rm d}x{\rm d}t
 \\ \nonumber&\hspace{0.5cm}-2 s\lambda\int_Q\varphi\left(\partial_{t}\psi\right){a_{8}}
 L_{1\ast}\left({\vartheta},{w}\right){\vartheta}{\rm d}x{\rm d}t
 \\ \nonumber&\hspace{0.5cm}+2 s^2\lambda^2\int_Q\varphi^2\left(\partial_{t}\psi\right){a_{8}}
 \gamma_{0}{w}{\vartheta}{\rm d}x{\rm d}t
+2 s\lambda\int_Q\varphi\left(\partial_{t}\psi\right){a_{8}}{\vartheta}
 \sum_{j=1}^{n}\gamma_{j}\partial_{j}{w}{\rm d}x{\rm d}t
 \\&\hspace{0.5cm}-C_{37}s^{2}\lambda^{4}\int_Q\varphi^{2}{\vartheta}^{2}
 {\rm d}x{\rm d}t-C_{37}\int_Q\left|\partial_{t}{w}\right|^{2}{\rm d}x{\rm d}t.
 \end{align*}
 Integrating the first and third terms by parts and estimating the second and fifth terms by the Cauchy-Schwarz inequality, we
 have
 \begin{align}\label{i51}
 \nonumber &\int_QI_{51}{\rm d}x{\rm d}t
 \geq 2 s\lambda\int_Q\varphi\left(\partial_{t}\psi\right){a_{1}}{a_{8}}
 \left(\nabla{\vartheta}\cdot\nabla{w}\right){\rm d}x{\rm d}t
 \\ \nonumber&\hspace{0.5cm}+2 s\lambda^2\int_Q\varphi\left(\partial_{t}\psi\right){a_{1}}{a_{8}}
 \left(\nabla\psi\cdot\nabla{w}\right){\vartheta}{\rm d}x{\rm d}t
 \\ \nonumber&\hspace{0.5cm}+2 s\lambda\int_Q\varphi\left(\partial_{t}\psi\right)
 \left(\nabla\left({a_{1}}{a_{8}}\right)\cdot\nabla{w}\right){\vartheta}{\rm d}x{\rm d}t
 \\ \nonumber&\hspace{0.5cm}
 + s\lambda^2\int_Q\varphi\left(\partial_{t}\psi\right)^2{a_{8}^2}
 {\vartheta}^2{\rm d}x{\rm d}t
 + s\lambda\int_Q\varphi\left(\partial_{t}^2\psi\right){a_{8}^2}
 {\vartheta}^2{\rm d}x{\rm d}t
 \\ \nonumber&\hspace{0.5cm}+2 s^2\lambda^2\int_Q\varphi^2\left(\partial_{t}\psi\right)^2{a_{8}^2}
 {\vartheta}^2{\rm d}x{\rm d}t
 -5s^2\lambda^2\int_Q\varphi^2\left(\partial_{t}\psi\right)^2{a_{8}}^2
 {\vartheta}^2{\rm d}x{\rm d}t
 \\ \nonumber&\hspace{0.5cm}-\frac{1}{4}\int_Q
 \left|L_{1\ast}\left({\vartheta},{w}\right)\right|^2{\rm d}x{\rm d}t
 -C_{38} s^3\lambda^3\int_Q\varphi^3\left({\vartheta}^2+{w}^2\right){\rm d}x{\rm d}t
 \\\nonumber&\hspace{0.5cm}-C_{38}s^{2}\lambda^{4}\int_Q\varphi^{2}{\vartheta}^{2}
 {\rm d}x{\rm d}t
 -C_{38}\int_Q\left(\left|\nabla{w}\right|^2+\left|\partial_{t}{w}\right|^{2}\right){\rm d}x{\rm d}t
 \\\nonumber &
 \geq 2 s\lambda\int_Q\varphi\left(\partial_{t}\psi\right){a_{1}}{a_{8}}
 \left(\nabla{\vartheta}\cdot\nabla{w}\right){\rm d}x{\rm d}t
 -\frac{1}{4}\int_Q
 \left|L_{1\ast}\left({\vartheta},{w}\right)\right|^2{\rm d}x{\rm d}t
 \\\nonumber&\hspace{0.5cm}-C_{39} s^3\lambda^3\int_Q\varphi^3\left({\vartheta}^2+{w}^2\right){\rm d}x{\rm d}t
 -C_{39}s^{2}\lambda^{4}\int_Q\varphi^{2}{\vartheta}^{2}
 {\rm d}x{\rm d}t
 \\&\hspace{0.5cm}-C_{39}\int_Q\left(\left|\nabla{w}\right|^2+\left|\partial_{t}{w}\right|^{2}\right){\rm d}x{\rm d}t.
 \end{align}
 Integrating by parts two times, we have
 \begin{align}\label{i52}
 \nonumber \int_Q&I_{52}{\rm d}x{\rm d}t=-2 s\lambda\int_Q\varphi{a_{1}a_{8}}\left(\nabla\psi\cdot\nabla{w}\right)
 \left(\partial_{t}{\vartheta}\right){\rm d}x{\rm d}t
  \\ \nonumber&=2 s\lambda\int_Q\varphi{a_{1}a_{8}}\left(\nabla\psi\cdot
 \nabla\left(\partial_{t}{w}\right)\right){\vartheta}{\rm d}x{\rm d}t
 \\\nonumber&\hspace{0.5cm}+2 s\lambda^{2}\int_Q\varphi\left(\partial_{t}\psi\right){a_{1}a_{8}}
 \left(\nabla\psi\cdot\nabla{w}\right){\vartheta}{\rm d}x{\rm d}t
 \\\nonumber&
 =-2 s\lambda\int_Q\varphi{a_{1}a_{8}}\left(\nabla\psi\cdot
 \nabla{\vartheta}\right)\left(\partial_{t}{w}\right){\rm d}x{\rm d}t
 \\\nonumber&\hspace{0.5cm}-2s\lambda\int_Q\varphi\left\{\mbox{div}\left({a_{1}a_{8}}\nabla\psi\right)\right\}
 \left(\partial_{t}{w}\right){\vartheta}{\rm d}x{\rm d}t
 \\\nonumber&\hspace{0.5cm}-2 s\lambda^{2}\int_Q\varphi\left|\nabla\psi\right|^{2}{a_{1}a_{8}}\left(\partial_{t}{w}\right){\vartheta}{\rm d}x{\rm d}t
\\\nonumber&\hspace{0.5cm}
 +2 s\lambda^{2}\int_Q\varphi\left(\partial_{t}\psi\right){a_{1}a_{8}}
 \left(\nabla\psi\cdot\nabla{w}\right){\vartheta}{\rm d}x{\rm d}t
 \\\nonumber&
 \geq -2 s\lambda\int_Q\varphi{a_{1}a_{8}}\left(\nabla\psi\cdot
 \nabla{\vartheta}\right)\left(\partial_{t}{w}\right){\rm d}x{\rm d}t
 \\&\hspace{0.5cm}-C_{40}s^{2}\lambda^{4}\int_Q\varphi^{2}{\vartheta}^{2}
 {\rm d}x{\rm d}t-C_{40}\int_Q\left(\left|\nabla{w}\right|^2+\left|\partial_{t}{w}\right|^{2}\right){\rm d}x{\rm d}t.
 \end{align}
 Integrating by parts, we have
 \begin{align}\label{i53}
 \nonumber \int_Q&I_{53}{\rm d}x{\rm d}t=s\lambda^{2}\int_Q\varphi\left(\partial_{t}\psi\right)^{2}{a_{7}}{a_{8}}
 \left(\partial_{t}{\vartheta}\right){w} {\rm d}x{\rm d}t
 \\\nonumber&=-s\lambda^{2}\int_Q\varphi
 \left(\partial_{t}\psi\right)^{2}{a_{7}}{a_{8}}{\vartheta}\left(\partial_{t}{w}\right){\rm d}x{\rm d}t
 \\\nonumber&\hspace{0.5cm}-2s\lambda^{2}\int_Q\varphi\left(\partial_{t}\psi\right)\left(\partial_{t}^2\psi\right){a_{7}}{a_{8}}
 {\vartheta}{w}{\rm d}x{\rm d}t
 -s\lambda^{3}\int_Q\varphi\left(\partial_{t}\psi\right)^{3}{a_{7}}{a_{8}}
 {\vartheta}{w}{\rm d}x{\rm d}t
 \\&\geq -C_{41}s^{2}\lambda^{4}\int_Q\varphi^{2}{\vartheta}^{2}
 {\rm d}x{\rm d}t-C_{41}\int_Q\left|\partial_{t}{w}\right|^{2}{\rm d}x{\rm d}t-C_{41}s\lambda^{3}\int_Q\varphi{w}^{2}
 {\rm d}x{\rm d}t.
 \end{align}
 Similarly, integrating by parts, we can obtain
 \begin{align}\label{i54}
 \int_Q&I_{54}{\rm d}x{\rm d}t
 \geq -C_{42}s^{2}\lambda^{4}\int_Q\varphi^{2}{\vartheta}^{2}
 {\rm d}x{\rm d}t-C_{42}\int_Q\left|\partial_{t}{w}\right|^{2}{\rm d}x{\rm d}t-C_{42}s\lambda^{3}\int_Q\varphi{w}^{2}
 {\rm d}x{\rm d}t.
 \end{align}
 Integrating by parts two times, we have
 \begin{align}\label{i55}
 \nonumber \int_Q&I_{55}{\rm d}x{\rm d}t=\int_Q{a_{8}}\sum_{j=1}^{n}\gamma_{j}\left(\partial_{j}{w}\right)\left(\partial_{t}{\vartheta}\right) {\rm d}x{\rm d}t
 =-\int_Q{a_{8}}\sum_{j=1}^{n}\gamma_{j}\left(\partial_{t}\partial_{j}{w}\right)
 {\vartheta}{\rm d}x{\rm d}t
 \\\nonumber&=\int_Q{a_{8}}\sum_{j=1}^{n}\gamma_{j}\left(\partial_{j}{\vartheta}\right)\left(\partial_{t}{w}\right){\rm d}x{\rm d}t
 +\int_Q\left\{\sum_{j=1}^{n}\partial_j
 \left({a_{8}}\gamma_{j}\right)\right\}\left(\partial_{t}{w}\right){\vartheta}{\rm d}x{\rm d}t
 \\ &\geq -C_{43}\int_Q\left(\left|\partial_{t}{w}\right|^{2}+\left|\nabla{\vartheta}\right|^{2}+|{\vartheta}|^2\right)
 {\rm d}x{\rm d}t.
 \end{align}
 Combining (\ref{i50})--(\ref{i55}), we arrive at
\begin{align}\label{i5}
&\nonumber \int_QI_{5}{\rm d}x{\rm d}t\geq
2 s\lambda\int_Q\varphi\left(\partial_{t}\psi\right){a_{1}}{a_{8}}
 \left(\nabla{\vartheta}\cdot\nabla{w}\right){\rm d}x{\rm d}t
 -\frac{1}{4}\int_Q
 \left|L_{1\ast}\left({\vartheta},{w}\right)\right|^2{\rm d}x{\rm d}t
 \\\nonumber&\hspace{0.5cm} -2 s\lambda\int_Q\varphi{a_{1}a_{8}}\left(\nabla\psi\cdot
 \nabla{\vartheta}\right)\left(\partial_{t}{w}\right){\rm d}x{\rm d}t
 \\\nonumber&\hspace{0.5cm}-C_{44} s^3\lambda^3\int_Q\varphi^3\left({\vartheta}^2+{w}^2\right){\rm d}x{\rm d}t
 -C_{44}s^{2}\lambda^{4}\int_Q\varphi^{2}{\vartheta}^{2}
 {\rm d}x{\rm d}t
 \\&\hspace{0.5cm}-C_{44}\int_Q\left(\left|\nabla{w}\right|^2+\left|\partial_{t}{w}\right|^{2}+\left|\nabla{\vartheta}\right|^{2}
 \right){\rm d}x{\rm d}t.
 \end{align}
 \begin{align}\label{i60}
 \nonumber \int_Q&I_{6}{\rm d}x{\rm d}t=s\lambda\int_Q\varphi\left(\partial_{t}\psi\right)a_8{\vartheta}L_{2\ast}\left({w}\right){\rm d}x{\rm d}t
 \\\nonumber&
 =-2s^2\lambda^2\int_Q\varphi^2\left(\partial_{t}\psi\right)^2{a_{7}}{a_{8}}{\vartheta}\left( \partial_t{w}\right){\rm d}x{\rm d}t
 \\\nonumber&\hspace{0.5cm}+2s^2\lambda^2\int_Q\varphi^2\left(\partial_{t}\psi\right)a_{1}a_8\left(\nabla\psi\cdot\nabla{w}\right){\vartheta}{\rm d}x{\rm d}t-s^2\lambda^{3}\int_Q\varphi^2\left(\partial_{t}\psi\right)^{3}{a_{7}}{a_{8}}{w}{\vartheta}{\rm d}x{\rm d}t
 \\\nonumber&\hspace{0.5cm}-s^2\lambda^2\int_Q\varphi^2\left(\partial_{t}\psi\right)
 \left(\partial^{2}_{t}\psi\right){a_{7}}{a_{8}}{\vartheta}{w}{\rm d}x{\rm d}t
 +s^2\lambda^{3}\int_Q\varphi^2\left(\partial_{t}\psi\right)\left|\nabla\psi\right|^{2}a_{1}{a_{8}}{\vartheta}{w}{\rm d}x{\rm d}t
 \\\nonumber&\hspace{0.5cm}+2ns^2\lambda^2\int_Q\varphi^2 \left(\partial_{t}\psi\right) {a_{1}}{a_{8}}{\vartheta}{w}{\rm d}x{\rm d}t
-s^2\lambda^2\int_Q\varphi^2\left(\partial_{t}\psi\right)
 {a_{8}}\gamma_0{\vartheta}{w}{\rm d}x{\rm d}t
 \\&\hspace{0.5cm}-s\lambda\int_Q\varphi\left(\partial_{t}\psi\right)a_8{\vartheta}\sum_{j=1}^{n}\gamma_{j}\partial_{j}{w}{\rm d}x{\rm d}t.
 \end{align}
 Similarly to (\ref{trans1}), we can estimate the first and second terms. Therefore we have
 \begin{align}\label{i6}
 \nonumber \int_Q&I_{6}{\rm d}x{\rm d}t
 \geq- C_{45}s\lambda^{\frac{1}{2}}\int_Q\varphi\left(\left|\nabla{w}\right|^{2}
 +\left|\partial_{t}{{w}}\right|^{2}\right){\rm d}x{\rm d}t
 -C_{45}s^{3}\lambda^{\frac{7}{2}}\int_Q\varphi^{3}{\vartheta}^{2}{\rm d}x{\rm d}t
 \\&\hspace{0.5cm}-C_{45} s^2\lambda^3\int_Q\varphi^2\left({\vartheta}^2+{w}^2\right){\rm d}x{\rm d}t
 -C_{45}\int_Q\left|\nabla{w}\right|^2{\rm d}x{\rm d}t.
 \end{align}
 By (\ref{j0}), we have
 \begin{align*}
 \nonumber \int_Q&I_{7}{\rm d}x{\rm d}t=\int_QJ_{0}L_{2\ast}\left({w}\right){\rm d}x{\rm d}t
 \\\nonumber&=s\lambda\int_Q\varphi \gamma_0 {w}L_{2\ast}\left({w}\right){\rm d}x{\rm d}t+\int_Q\left(\sum_{j=1}^n \gamma_j \partial_j{w}\right)L_{2\ast}\left({w}\right){\rm d}x{\rm d}t.
 \end{align*}
 Similarly to (\ref{i60})-(\ref{i6}), we can estimate the first term.
 \begin{align*}
 \nonumber&(\mbox{second term})=
 -2 s\lambda\int_Q\varphi\left(\partial_{t}\psi\right){a_{7}}\left(\partial_{t}{w}\right)\sum_{j=1}^{n}\gamma_{j}
 \left(\partial_{j}{w}\right){\rm d}x{\rm d}t
 \\\nonumber&\hspace{0.5cm}+2 s\lambda\int_Q\varphi a_{1}\left(\nabla\psi\cdot\nabla{w}\right)\left(\sum_{j=1}^{n}\gamma_{j}
 \partial_{j}{w}\right){\rm d}x{\rm d}t
 \\\nonumber&\hspace{0.5cm}-s\lambda^{2}\int_Q\varphi\left(\partial_{t}\psi\right)^{2}{a_{7}}\sum_{j=1}^{n}\gamma_{j}
 \left(\partial_{j}{w}\right){w}{\rm d}x{\rm d}t
 \\\nonumber&\hspace{0.5cm}-s\lambda\int_Q\varphi\left(\partial^{2}_{t}\psi\right){a_{7}}\sum_{j=1}^{n}\gamma_{j}
 \left(\partial_{j}{w}\right){w}{\rm d}x{\rm d}t
 \\\nonumber&\hspace{0.5cm}+s\lambda^{2}\int_Q\varphi\left|\nabla\psi\right|^{2}a_{1}\sum_{j=1}^{n}\gamma_{j}\left(\partial_{j}{w}\right){w}
 +2n s\lambda\int_Q\varphi a_{1}\sum_{j=1}^{n}\gamma_{j}\left(\partial_{j}{w}\right){w}{\rm d}x{\rm d}t
 \\\nonumber&\hspace{0.5cm}
 -s\lambda\int_Q\varphi\gamma_{0}\sum_{j=1}^{n}\gamma_{j}\left(\partial_{j}{w}\right){w}{\rm d}x{\rm d}t
 -\int_Q\left(\sum_{j=1}^{n}\gamma_{j}\partial_{j}{w}\right)^{2}{\rm d}x{\rm d}t
 \\\nonumber&\geq
 -2 s\lambda\int_Q\varphi\left(\partial_{t}\psi\right){a_{7}}\left(\partial_{t}{w}\right)\sum_{j=1}^{n}\gamma_{j}
 \left(\partial_{j}{w}\right){\rm d}x{\rm d}t
 \\\nonumber&\hspace{0.5cm}+2 s\lambda\int_Q\varphi a_{1}\left(\nabla\psi\cdot\nabla{w}\right)\left(\sum_{j=1}^{n}\gamma_{j}
 \partial_{j}{w}\right){\rm d}x{\rm d}t
 \\&\hspace{0.5cm}-C_{46} s^2\lambda^4\int_Q\varphi^2{w}^2{\rm d}x{\rm d}t
 -C_{46}\int_Q\left|\nabla{w}\right|^2{\rm d}x{\rm d}t.
 \end{align*}
 Therefore we can obtain
 \begin{align}\label{i7}
 \nonumber \int_Q&I_{7}{\rm d}x{\rm d}t\geq-2 s\lambda\int_Q\varphi\left(\partial_{t}\psi\right){a_{7}}\left(\partial_{t}{w}\right)\sum_{j=1}^{n}\gamma_{j}
 \left(\partial_{j}{w}\right){\rm d}x{\rm d}t
 \\\nonumber&\hspace{0.5cm}+2 s\lambda\int_Q\varphi a_{1}\left(\nabla\psi\cdot\nabla{w}\right)\left(\sum_{j=1}^{n}\gamma_{j}
 \partial_{j}{w}\right){\rm d}x{\rm d}t
 \\\nonumber&\hspace{0.5cm}- C_{47}s\lambda^{\frac{1}{2}}\int_Q\varphi\left(\left|\nabla{w}\right|^{2}
 +\left|\partial_{t}{{w}}\right|^{2}\right){\rm d}x{\rm d}t
 -C_{47}s^{3}\lambda^{\frac{7}{2}}\int_Q\varphi^{3}{w}^{2}{\rm d}x{\rm d}t
 \\&\hspace{0.5cm}-C_{47} s^2\lambda^4\int_Q\varphi^2{w}^2{\rm d}x{\rm d}t
 -C_{47}\int_Q\left|\nabla{w}\right|^2{\rm d}x{\rm d}t.
 \end{align}
 We combine of (\ref{eq}), (\ref{ik}), (\ref{i1})--(\ref{i4}), (\ref{i5}), (\ref{i6}), and (\ref{i7}).
 We note that the sum of the first and third terms in (\ref{i1}) and the third and fifth terms in (\ref{i2}) equals
 \begin{align*}
 2s\lambda^{2}\int_Q\varphi\left|{a_{7}}\left(\partial_{t}\psi\right)
 \left(\partial_{t}{w}\right)-a_{1}\left(\nabla\psi\cdot\nabla{w}\right)\right|^{2}{\rm d}x{\rm d}t.
 \end{align*}
 Moreover, we take
 \begin{align}
 \gamma_{0}=a_{1}-\frac{\left(\nabla\left(a_{1}{a_{7}}\right)\cdot\nabla\psi\right)}
 {a_7},\quad \gamma_{j}=-2\left(\partial_{j}a_{1}\right),\ j=1,2,\cdots, n.
 \end{align}
 Then the sum of the fourth and fifth terms in (\ref{i1}) equals to $s\lambda\int_Q \varphi {a_{1}}{a_{7}}\left|\partial_{t}{w}\right|^{2}{\rm d}x{\rm d}t$. The sum of the sixth term in (\ref{i2}) and the second term in (\ref{i7}) equals $0$.
 The sum of the fourth term in (\ref{i2}) and the first term in (\ref{i7}) equals to
 $$
  2s\lambda\int_Q\varphi\left(\partial_t\psi\right)
    a_7\left(\nabla a_1\cdot\nabla {w}\right)\left(\partial_t{w}\right) {\rm d}x{\rm d}t
    -2s\lambda\int_Q\varphi\left(\partial_t\psi\right)
    a_1\left(\nabla a_7\cdot\nabla {w}\right)\left(\partial_t{w}\right) {\rm d}x{\rm d}t.
 $$
 We combine the first term in (\ref{i3}) and the first term in (\ref{i4}).
 Thus we arrive at
 \begin{align}\label{g3}
 \nonumber&\left\|\left(f-{a_{8}}g\right){\rm e}^{s\varphi}\right\|_{L^2(Q)}^2=\left\|\mathcal{L}_{\ast}\left({\vartheta},{w}\right)\right\|_{L^2(Q)}^2
 \\\nonumber&=\left\|\mathcal{L}_{1\ast}\left({\vartheta},{w}\right)\right\|_{L^2(Q)}^2+2\int_QL_{1\ast}\left({\vartheta},{w}\right)L_{2\ast}\left({w}\right){\rm d}x{\rm d}t+\left\|\mathcal{L}_{2\ast}\left({w}\right)\right\|_{L^2(Q)}^2
 \\\nonumber&\geq \frac{1}{2}\left\|\mathcal{L}_{\ast}\left({\vartheta},{w}\right)\right\|_{L^2(Q)}^2
 +\left\|\mathcal{L}_{2\ast}\left({w}\right)\right\|_{L^2(Q)}^2
 \\\nonumber&\hspace{0.2cm}-2s\lambda\int_\Sigma\varphi(\nabla\psi\cdot\nu)a_{1}^2\left|\frac{\partial{w}}{\partial\nu}\right|^2{\rm d}\sigma{\rm d}t
 +2\int_\Sigma a_1\left(\nabla{a_1}\cdot\nu\right)\left|\frac{\partial{w}}{\partial\nu}\right|^2{\rm d}\sigma{\rm d}t
 \\\nonumber&\hspace{0.2cm}+4s\lambda^{2}\int_Q\varphi\left|{a_{7}}\left(\partial_{t}\psi\right)
 \left(\partial_{t}{w}\right)-a_{1}\left(\nabla\psi\cdot\nabla{w}\right)\right|^{2}{\rm d}x{\rm d}t
 +2s\lambda\int_Q \varphi {a_{1}}{a_{7}}\left|\partial_{t}{w}\right|^{2}{\rm d}x{\rm d}t
 \\\nonumber&\hspace{0.2cm}+2s\lambda\int_Q\varphi a_{1}\left\{3a_{1}-\left(\nabla a_{1}\cdot\nabla\psi\right)+\frac{a_{1}\left(\nabla{a_{7}}\cdot\nabla\psi\right)}
 {a_{7}}\right\}\left|\nabla{w}\right|^{2}{\rm d}x{\rm d}t
 \\\nonumber&\hspace{0.2cm}+4s^{3}\lambda^{4}\int_Q\varphi^{3}\left|\nabla\psi\right|^{4}a_{1}^{2}{w}^{2}{\rm d}x{\rm d}t
 -4 s\lambda\int_Q\varphi{a_{1}a_{8}}\left(\nabla\psi\cdot
 \nabla{\vartheta}\right)\left(\partial_{t}{w}\right){\rm d}x{\rm d}t
 \\\nonumber&\hspace{0.2cm}+4s\lambda\int_Q\varphi
 \left(\partial^{2}_{t}\psi\right){a_{7}^{2}}\left|\partial_{t}{w}\right|^{2}{\rm d}x{\rm d}t
 +4s\lambda\int_Q\varphi\left(\partial_t\psi\right)
    a_7\left(\nabla a_1\cdot\nabla {w}\right)\left(\partial_t{w}\right) {\rm d}x{\rm d}t
    \\\nonumber&\hspace{0.2cm}-4s\lambda\int_Q\varphi\left(\partial_t\psi\right)
    a_1\left(\nabla a_7\cdot\nabla {w}\right)\left(\partial_t{w}\right) {\rm d}x{\rm d}t
    \\\nonumber&\hspace{0.2cm}
    +4 s\lambda\int_Q\varphi\left(\partial_{t}\psi\right){a_{1}}{a_{8}}
  \left(\nabla{\vartheta}\cdot\nabla{w}\right){\rm d}x{\rm d}t
 \\\nonumber&\hspace{0.2cm}-8s^{3}\lambda^{4}\int_Q\varphi^{3}\left(\partial_{t}\psi\right)^{2}
 \left|\nabla\psi\right|^{2}a_{1}{a_{7}}{w}^{2}{\rm d}x{\rm d}t
 \\\nonumber&\hspace{0.2cm}-C_{48}s\lambda^{\frac{1}{2}}\int_Q\varphi\left(\left|\nabla{w}\right|^{2}
 +\left|\partial_{t}{w}\right|^{2}\right){\rm d}x{\rm d}t
 -C_{48}s^{3}\lambda^{\frac{7}{2}}\int_Q\varphi^{3}
 \left({\vartheta}^{2}+{w}^{2}\right){\rm d}x{\rm d}t
 \\&\hspace{0.5cm}-C_{48}\int_Q\left(\left|\nabla{\vartheta}\right|^{2}
 +s^2\lambda^{4}\varphi^2{\vartheta}^{2}+s^2\lambda^{6}\varphi^2{w}^{2}
 \right){\rm d}x{\rm d}t.
 \end{align}
 By $\big(a_1(x), a_2(x), a_3(x), a_4(x)\big)\in \mathcal{U}$ and the Cauchy-Schwarz inequality, we estimate every term of (\ref{g3}). We note that $\|a_{7}\|_{C(\overline{\Omega})}\leq 1+\frac{M^{2}_{0}}{\sigma_{1}}$, $\|a_{8}\|_{C(\overline{\Omega})}\leq \frac{M_{0}}{\sigma_{1}}$,
 $\|\partial_ja_{7}\|_{C(\overline{\Omega})}\leq \frac{3M^{2}_{0}M_1}{\sigma_{1}^2}$, and $\|\partial_ja_{8}\|_{C(\overline{\Omega})}\leq \frac{2M_1M_{0}}{\sigma_{1}^2}$, $j=1$, $2$, $\cdots$, $n$.
 \begin{align}\label{eq1}
 2a_{1}{a_{7}}s\lambda\varphi\left|\partial_{t}{w}\right|^{2}
 \geq 2\sigma_{1}s\lambda\varphi\left|\partial_{t}{w}\right|^{2},
 \end{align}
 \begin{align}\label{eq3}
 2s\lambda\varphi a_{1}\left\{3a_{1}-\left(\nabla a_{1}\cdot\nabla\psi\right)+\frac{a_{1}\left(\nabla{a_{7}}\cdot\nabla\psi\right)}
 {a_{7}}\right\}\left|\nabla{w}\right|^{2}
 \geq2\sigma_{0}\sigma_{1}s\lambda\varphi\left|\nabla{w}\right|^{2},
 \end{align}
 and
 \begin{align}\label{eq6}
 4 s^{3}\lambda^{4}\varphi^{3}\left|\nabla\psi\right|^{4}a_{1}^{2}{w}^{2}
 \geq 4\sigma^{2}_{1}s^{3}\lambda^{4}\varphi^{3}\left|\nabla\psi\right|^{4}{w}^{2},\ \mbox{in} \ Q.
 \end{align}
 By the Cauchy-Schwarz inequality, we have, for any $\varepsilon_{4}>0$,
\begin{align*}
 &\nonumber\left|-4 s\lambda\varphi{a_{1}a_{8}}\left(\nabla\psi\cdot\nabla{\vartheta}\right)
 \partial_{t}{{w}}\right|
 \leq\frac{2M^{2}_{0}}{\sigma_{1}}s\lambda\varphi\left(\varepsilon_{4}\left|\partial_{t}{w}\right|^{2}
 +\varepsilon^{-1}_{4}\left|\left(\nabla\psi\cdot\nabla{\vartheta}\right)\right|^{2}\right),\ \mbox{in} \ Q.
 \end{align*}
 Choosing $\varepsilon_{4}>0$ small such that $\frac{2M^{2}_{0}}{\sigma_{1}}
 \varepsilon_{4}=\sigma_{1}$, that is,
 $\varepsilon_{4}=\frac{\sigma^{2}_{1}}{2M^{2}_{0}}$, we have
 \begin{align}\label{ee2}
 \left|-4 s\lambda\varphi{a_{1}a_{8}}\left(\nabla\psi\cdot\nabla{\vartheta}\right)
 \partial_{t}{{w}}\right|
 \leq \sigma_{1}s\lambda\varphi\left|\partial_{t}{w}\right|^{2}
 +\frac{4M^{4}_{0}}{\sigma^{3}_{1}}s\lambda\varphi\left|\left(\nabla\psi\cdot\nabla{\vartheta}\right)\right|^{2},\ \mbox{in} \ Q.
 \end{align}
 Next we estimate the second term of (\ref{ee2}).
 \begin{align}\label{ee3}
 &\nonumber s\lambda\varphi\left|\left(\nabla\psi\cdot\nabla{\vartheta}\right)\right|^{2}
 =\frac{1}{a^{2}_{5}}s\lambda\varphi\left|a_{5}\left(\nabla\psi\cdot\nabla{\vartheta}\right)
 +a_{6}\left(\nabla\psi\cdot\nabla{w}\right)-a_{6}\left(\nabla\psi\cdot\nabla{w}\right)\right|^{2}
 \\\nonumber&\leq\frac{2}{a^{2}_{5}}s\lambda\varphi\left|a_{5}\left(\nabla\psi\cdot\nabla{\vartheta}\right)
 +a_{6}\left(\nabla\psi\cdot\nabla{w}\right)\right|^{2}
 +\frac{2a^{2}_{6}}{a^{2}_{5}}s\lambda\varphi\left|\left(\nabla\psi\cdot\nabla{w}\right)\right|^{2}
 \\\nonumber&=\frac{2}{a^{2}_{5}}s\lambda\varphi\left|a_{5}\left(\nabla\psi\cdot\nabla{\vartheta}\right)
 +a_{6}\left(\nabla\psi\cdot\nabla{w}\right)\right|^{2}
 \\\nonumber&\hspace{0.5cm}+\frac{2a^{2}_{6}}{a^{2}_{1}a^{2}_{5}}s\lambda\varphi
 \left|-a_{1}\left(\nabla\psi\cdot\nabla{w}\right)+{a_{7}}
 \left(\partial_{t}\psi\right)\left(\partial_{t}{{w}}\right)
 -{a_{7}}\left(\partial_{t}\psi\right)\left(\partial_{t}{{w}}\right)\right|^{2}
 \\\nonumber&\leq\frac{2}{a^{2}_{5}}s\lambda\varphi\left|a_{5}\left(\nabla\psi\cdot\nabla{\vartheta}\right)
 +a_{6}\left(\nabla\psi\cdot\nabla{w}\right)\right|^{2}+\frac{4a^{2}_{6}{a_{7}^{2}}}{a^{2}_{1}a^{2}_{5}}s\lambda\varphi
 \left(\partial_{t}\psi\right)^{2}\left|\partial_{t}{{w}}\right|^{2}
 \\\nonumber&\hspace{0.5cm}+\frac{4a^{2}_{6}}{a^{2}_{1}a^{2}_{5}}s\lambda\varphi
 \left|{a_{7}}\left(\partial_{t}\psi\right)\left(\partial_{t}{{w}}\right)
 -a_{1}\left(\nabla\psi\cdot\nabla{w}\right)\right|^{2}
 \\\nonumber&\leq \frac{4M^{4}_{0}\beta^{2}T^{2}}{\sigma^{4}_{1}}\left(1+\frac{M^{2}_{0}}{\sigma_{1}}\right)^{2}
 s\lambda\varphi\left|\partial_{t}{{w}}\right|^{2}
 +C_{49}s\lambda\varphi\left|a_{5}\left(\nabla\psi\cdot\nabla{\vartheta}\right)
+a_{6}\left(\nabla\psi\cdot\nabla{w}\right)\right|^{2}
\\&\hspace{0.5cm}+C_{49}s\lambda\varphi\left|{a_{7}}\left(\partial_{t}\psi\right)
\left(\partial_{t}{{w}}\right)-a_{1}\left(\nabla\psi\cdot\nabla{w}\right)\right|^{2}.
\end{align}
By (\ref{ee2}) and (\ref{ee3}), we obtain
\begin{align}\label{ee4}
 &\nonumber\left|-4 s\lambda\varphi{a_{1}a_{8}}\left(\nabla\psi\cdot\nabla{\vartheta}\right)
 \partial_{t}{{w}}\right|
 \\\nonumber&\leq \sigma_{1}s\lambda\varphi\left|\partial_{t}{w}\right|^{2}
 +\frac{16M^{8}_{0}\beta^{2}T^{2}}{\sigma^{7}_{1}}\left(1+\frac{M^{2}_{0}}{\sigma_{1}}\right)^{2}
 s\lambda\varphi\left|\partial_{t}{{w}}\right|^{2}
 \\\nonumber&\hspace{0.5cm}+C_{50}s\lambda\varphi\left|a_{5}\left(\nabla\psi\cdot\nabla{\vartheta}\right)
 +a_{6}\left(\nabla\psi\cdot\nabla{w}\right)\right|^{2}
 \\&\hspace{0.5cm}+C_{50}s\lambda\varphi\left|{a_{7}}\left(\partial_{t}\psi\right)
  \left(\partial_{t}{{w}}\right)-a_{1}\left(\nabla\psi\cdot\nabla{w}\right)\right|^{2},\ \mbox{in}\ Q.
 \end{align}
 Furthermore, we have
\begin{align}\label{eq2}
 &\left|4 s\lambda\varphi\left(\partial^{2}_{t}\psi\right){a_{7}^{2}}
  \left|\partial_{t}{w}\right|^{2}\right|
 \leq8\left(1+\frac{M^{2}_{0}}{\sigma_{1}}\right)^{2}\beta s\lambda\varphi\left|\partial_{t}{w}\right|^{2},
 \end{align}
\begin{align}\label{eq4}
&\nonumber\left|4s\lambda\varphi\left(\partial_t\psi\right)
    a_7\left(\nabla a_1\cdot\nabla {w}\right)\left(\partial_t{w}\right)\right|\leq 4\left(1+\frac{M^{2}_{0}}{\sigma_{1}}\right)M_{1}\beta Ts\lambda\varphi\sum_{j=1}^n\left|\partial_t{w}\right|\left|\partial_j {w}\right|
 \\&\leq 2\left(1+\frac{M^{2}_{0}}{\sigma_{1}}\right)M_{1}\beta Ts\lambda\varphi
 \left(n\left|\partial_{t}{{w}}\right|^{2}+\left|\nabla{w}\right|^{2}\right),
 \end{align}
 \begin{align}\label{eq5}
 &\left|-4s\lambda\varphi\left(\partial_t\psi\right)
     a_1\left(\nabla a_7\cdot\nabla {w}\right)\left(\partial_t{w}\right)\right|
 \leq\frac{6M^{3}_{0}M_1}{\sigma_{1}^2}\beta Ts\lambda\varphi\left(n\left|\partial_{t}{w}\right|^{2}
 +\left|\nabla{w}\right|^{2}\right),
 \end{align}
 \begin{align}\label{e2}
 \nonumber\left|4 s\lambda\varphi\left(\partial_{t}\psi\right){a_{1}a_{8}}
 \left(\nabla{\vartheta}\cdot\nabla{{w}}\right)\right|&\leq 4 \beta Ts\lambda\varphi\left|{a_{1}a_{8}}\right| |\nabla{\vartheta}||\nabla{{w}}|
 \\&\leq\frac{2M^{2}_{0}}{\sigma_{1}}\beta Ts\lambda\varphi\left(\left|\nabla{\vartheta}\right|^{2}
 +\left|\nabla{w}\right|^{2}\right),
 \end{align}
 and \begin{align}\label{eq8}
 &\left|-8 s^{3}\lambda^{4}\varphi^{3}\left(\partial_{t}\psi\right)^{2}
 \left|\nabla\psi\right|^{2}a_{1}{a_{7}}{w}^{2}\right|
 \leq 8M_{0}\left(1+\frac{M^{2}_{0}}{\sigma_{1}}\right)\beta^{2} T^{2}s^{3}\lambda^{4}\varphi^{3}\left|\nabla\psi\right|^{2}{w}^{2},
 \end{align}
 in $Q$.
 Combining (\ref{g3})--(\ref{eq6}) and (\ref{ee4})--(\ref{eq8}), we can obtain
 \begin{align}\label{equation3}
 &\nonumber\frac{1}{2}\left\|\mathcal{L}_{\ast}\left({\vartheta},{w}\right)\right\|_{L^2(Q)}^2
 +\left\|\mathcal{L}_{2\ast}\left({w}\right)\right\|_{L^2(Q)}^2
 \\\nonumber&\hspace{0.2cm}-2s\lambda\int_\Sigma\varphi(\nabla\psi\cdot\nu)a_{1}^2\left|\frac{\partial{w}}{\partial\nu}\right|^2{\rm d}\sigma{\rm d}t
 +2\int_\Sigma  a_1\left(\nabla{a_1}\cdot\nu\right)\left|\frac{\partial{w}}{\partial\nu}\right|^2{\rm d}\sigma{\rm d}t
 \\\nonumber&\hspace{0.2cm}+4s\lambda^{2}\int_Q\varphi\left|{a_{7}}\left(\partial_{t}\psi\right)
 \left(\partial_{t}{w}\right)-a_{1}\left(\nabla\psi\cdot\nabla{w}\right)\right|^{2}{\rm d}x{\rm d}t
 +\sigma_{1}s\lambda\int_{Q}\varphi\left|\partial_{t}{w}\right|^{2}{\rm d}x{\rm d}t
 \\\nonumber&\hspace{0.2cm}+2\sigma_{0}\sigma_{1}s\lambda\int_{Q}\varphi\left|\nabla{w}\right|^{2}{\rm d}x{\rm d}t
 +4\sigma^{2}_{1}s^{3}\lambda^{4}\int_{Q}\varphi^{3}
  \left|\nabla\psi\right|^{4}{w}^{2}{\rm d}x{\rm d}t
 \\\nonumber&\leq C_{51} \int_{Q}\left(f^{2}+g^2\right){\rm e}^{2s\varphi}{\rm d}x{\rm d}t
  \\\nonumber&\hspace{0.2cm}+C_{51} s\lambda\int_{Q}\varphi\left|a_{5}\left(\nabla\psi\cdot\nabla{\vartheta}\right)
 +a_{6}\left(\nabla\psi\cdot\nabla{w}\right)\right|^{2}{\rm d}x{\rm d}t
 \\\nonumber&\hspace{0.2cm}+C_{51} s\lambda\int_{Q}\varphi\left|{a_{7}}
 \left(\partial_{t}\psi\right)\left(\partial_{t}{{w}}\right)
 -a_{1}\left(\nabla\psi\cdot\nabla{w}\right)\right|^{2}{\rm d}x{\rm d}t
 \\\nonumber&\hspace{0.2cm}+\frac{16M^{8}_{0}\beta^{2}T^{2}}{\sigma^{7}_{1}}\left(1+\frac{M^{2}_{0}}{\sigma_{1}}\right)^{2}
 s\lambda\int_{Q}\varphi\left|\partial_{t}{{w}}\right|^{2}{\rm d}x{\rm d}t
 \\\nonumber&\hspace{0.2cm}+8\left(1+\frac{M^{2}_{0}}{\sigma_{1}}\right)^{2}\beta s\lambda\int_{Q}\varphi\left|\partial_{t}{w}\right|^{2}{\rm d}x{\rm d}t
 \\\nonumber&\hspace{0.2cm}+2\left(1+\frac{M^{2}_{0}}{\sigma_{1}}+\frac{3M^{3}_{0}}{\sigma^{2}_{1}}\right)M_{1}\beta T s\lambda\int_{Q}\varphi
 \left(n\left|\partial_{t}{{w}}\right|^{2}+\left|\nabla{w}\right|^{2}\right){\rm d}x{\rm d}t
 \\\nonumber&\hspace{0.2cm}+\frac{2M^{2}_{0}}{\sigma_{1}}\beta T s\lambda\int_{Q}\varphi\left(\left|\nabla{\vartheta}\right|^{2}
 +\left|\nabla{w}\right|^{2}\right){\rm d}x{\rm d}t
 \\\nonumber&\hspace{0.2cm}
 +8M_{0}\left(1+\frac{M^{2}_{0}}{\sigma_{1}}\right)\beta^{2} T^{2}s^{3}\lambda^{4}\int_{Q}\varphi^{3}\left|\nabla\psi\right|^{2}{w}^{2}{\rm d}x{\rm d}t
 \\\nonumber&\hspace{0.2cm}+C_{51}s\lambda^{\frac{1}{2}}\int_Q\varphi\left(\left|\nabla{w}\right|^{2}
 +\left|\partial_{t}{w}\right|^{2}\right){\rm d}x{\rm d}t
 +C_{51}s^{3}\lambda^{\frac{7}{2}}\int_Q\varphi^{3}
 \left({\vartheta}^{2}+{w}^{2}\right){\rm d}x{\rm d}t
 \\&\hspace{0.2cm}+C_{51}\int_Q\left(\left|\nabla{\vartheta}\right|^{2}
 +s^2\lambda^{4}\varphi^2{\vartheta}^{2}+s^2\lambda^{6}\varphi^2{w}^{2}
 \right){\rm d}x{\rm d}t.
 \end{align}

 We multiply (\ref{equation12}) by $\alpha_{3}$ and then add it to (\ref{equation3}). By the definition of $\alpha_3$ in (\ref{alpha}), the second, third and fourth terms in the right hand side of (\ref{equation12}) can be absorbed by the
   left hand side of (\ref{equation3}),  so that we arrive at
 \begin{align}\label{f1}
 &\nonumber 2\alpha_{3}\sigma^{2}_{1}s^{3}\lambda^{4}\int_{Q}\varphi^{3}
 \left|\nabla\psi\right|^{4}{\vartheta}^{2}{\rm d}x{\rm d}t
 +\frac{\alpha_{3}}{60M}\int_{Q}\frac{1}{s\varphi}\left|\partial_{t}{\vartheta}\right|^{2}{\rm d}x{\rm d}t
 \\\nonumber&\hspace{0.2cm}+\alpha_{1}\alpha_{3}s\lambda\int_{Q}\varphi
 \left|\nabla{\vartheta}\right|^{2}{\rm d}x{\rm d}t
 \\\nonumber&\hspace{0.2cm}+2\alpha_{3}s\lambda^{2}\int_{Q}\varphi
 \left|a_{5}\left(\nabla\psi\cdot\nabla{\vartheta}\right)+a_{6}\left(\nabla\psi\cdot\nabla{w}\right)\right|^{2}{\rm d}x{\rm d}t
 \\\nonumber&\hspace{0.3cm}+4 s\lambda^{2}\int_{Q}\varphi\left|{a_{7}}
 \left(\partial_{t}\psi\right)
 \left(\partial_{t}{w}\right)-a_{1}\left(\nabla\psi\cdot\nabla{w}\right)\right|^{2}{\rm d}x{\rm d}t
 \\\nonumber&\hspace{0.2cm}+\frac{1}{2}\sigma_{1}s\lambda\int_{Q}\varphi\left|\partial_{t}{w}\right|^{2}{\rm d}x{\rm d}t
 +\sigma_{0}\sigma_{1}s\lambda\int_{Q}\varphi\left|\nabla{w}\right|^{2}{\rm d}x{\rm d}t \\\nonumber&\hspace{0.2cm}+2\sigma^{2}_{1}s^{3}\lambda^{4}\int_{Q}\varphi^{3}
  \left|\nabla\psi\right|^{4}{w}^{2}{\rm d}x{\rm d}t
 +\frac{1}{2}\left\|\mathcal{L}_{\ast}\left({\vartheta},{w}\right)\right\|_{L^2(Q)}^2
 +\left\|\mathcal{L}_{2\ast}\left({w}\right)\right\|_{L^2(Q)}^2
 \\\nonumber&\hspace{0.2cm}+\alpha_{3}\left\|L_{1}\left({\vartheta},{w}\right)\right\|^{2}_{L^2(Q)}-2\alpha_{3}s\lambda\int_\Sigma \varphi\left(\nabla\psi\cdot\nu\right)\left(a_{5}
 \frac{\partial {\vartheta}}{\partial\nu}+a_{6}
 \frac{\partial{w}}{\partial \nu}\right)^2{\rm d}\sigma{\rm d}t\\\nonumber&\hspace{0.2cm}-2s\lambda\int_\Sigma\varphi(\nabla\psi\cdot\nu)a_{1}^2\left|\frac{\partial{w}}{\partial\nu}\right|^2{\rm d}\sigma{\rm d}t
 +2\int_\Sigma  a_1\left(\nabla{a_1}\cdot\nu\right)\left|\frac{\partial{w}}
 {\partial\nu}\right|^2{\rm d}\sigma{\rm d}t
 \\\nonumber&\leq C_{52}\int_{Q}\left(f^{2}+g^{2}\right){\rm e}^{2s\varphi}{\rm d}x{\rm d}t
 +C_{52} s\lambda\int_{Q}\varphi\left|a_{5}\left(\nabla\psi\cdot\nabla{\vartheta}\right)
 +a_{6}\left(\nabla\psi\cdot\nabla{w}\right)\right|^{2}{\rm d}x{\rm d}t
 \\\nonumber&\hspace{0.2cm} +C_{52} s\lambda\int_{Q}\varphi\left|{a_{7}}
 \left(\partial_{t}\psi\right)\left(\partial_{t}{{w}}\right)
 -a_{1}\left(\nabla\psi\cdot\nabla{w}\right)\right|^{2}{\rm d}x{\rm d}t
 \\\nonumber&\hspace{0.2cm}+\beta Ts\lambda\int_{Q}\varphi\left(\alpha_{4}\left|\nabla{{\vartheta}}\right|^{2}
 +\alpha_{5}\left|\nabla{{w}}\right|^{2}+\alpha_{6}\left|\partial_{t}{{w}}\right|^{2}\right){\rm d}x{\rm d}t
 \\\nonumber&\hspace{0.2cm}+\alpha_{7}\beta^{2}T^{2}
 s\lambda\int_{Q}\varphi\left|\partial_{t}{{w}}\right|^{2}{\rm d}x{\rm d}t
 +4\alpha_{8}\beta s\lambda\int_{Q}\varphi\left|\partial_{t}{w}\right|^{2}{\rm d}x{\rm d}t
 \\\nonumber&\hspace{0.2cm}
 +\frac{8M_{0}\left(\sigma_{1}+M^{2}_{0}\right)}{\sigma_{1}}\beta^{2} T^{2}s^{3}\lambda^{4}\int_{Q}\varphi^{3}\left|\nabla\psi\right|^{2}{w}^{2}{\rm d}x{\rm d}t
 \\\nonumber&\hspace{0.2cm}+C_{52} s\lambda^{\frac{1}{2}}\int_{Q}\varphi\left(\left|\nabla{\vartheta}\right|^{2}
 +\left|\nabla{w}\right|^{2}+\left|\partial_{t}{w}\right|^{2}\right){\rm d}x{\rm d}t
 +C_{52} s^{3}\lambda^{\frac{7}{2}}\int_{Q}\varphi^{3}\left({\vartheta}^{2}+{w}^{2}\right){\rm d}x{\rm d}t
 \\&\hspace{0.2cm}+C_{52}\int_{Q}\left(\frac{1}{s\lambda\varphi}\left|\partial_{t}{{\vartheta}}\right|^{2}
 +s^{2}\lambda^{6}\varphi^{2}{\vartheta}^{2}+s^{2}\lambda^{6}\varphi^{2}{w}^{2}\right){\rm d}x{\rm d}t,
 \end{align}
 for all $s>s_1$.
 We note that $4m\leq|\nabla\psi|^{2}\leq4M$ on $\overline{\Omega}$.
 By (\ref{alpha}) and (\ref{beta}), we have
 \begin{equation}\label{beta3}
 2{\mathcal{D}}\alpha_6\sqrt{\beta}+4{\mathcal{D}}^2\alpha_7\beta+4\alpha_8\beta<
 \frac{1}{4}\sigma_1.
 \end{equation}
 By (\ref{beta}) and (\ref{beta3}),  there exists a sufficiently small constant $\eta=\eta(\beta)>0$ such that
 \begin{align*}
 & 0<\sqrt{\beta}\left({\mathcal{D}}+\eta\right)<\frac{\alpha_{1}\alpha_{3}}{4\alpha_{4}},
 \qquad 0<\sqrt{\beta}\left({\mathcal{D}}+\eta\right)<
 \frac{\sigma_{0}\sigma_{1}}{4\alpha_{5}},
 \\&2\left({\mathcal{D}}+\eta\right)\alpha_6\sqrt{\beta}+4\left({\mathcal{D}}+\eta\right)^2\alpha_7\beta+4\alpha_8\beta<
 \frac{1}{4}\sigma_1,
 \\&0<\beta\left({\mathcal{D}}+\eta\right)^{2}<\frac{m^{2}\sigma^{3}_{1}}{2M_{0}\left(\sigma_{1}+M^{2}_{0}\right)}.
 \end{align*}
 For any $T\in\left(0,\frac{2({\mathcal{D}}+\eta)}{\sqrt{\beta}}\right)$, we have $0<\beta T<2\sqrt{\beta}\left({\mathcal{D}}+\eta\right)$. Therefore, the fourth, fifth, sixth and seventh terms in the right hand side of (\ref{f1}) can be absorbed by the left hand side of it. Hence we arrive at
 \begin{align*}
 &\nonumber 32\alpha_{3}m^2\sigma^{2}_{1}s^{3}\lambda^{4}\int_{Q}\varphi^{3}
 {\vartheta}^{2}{\rm d}x{\rm d}t
 +\frac{\alpha_{3}}{60M}\int_{Q}\frac{1}{s\varphi}\left|\partial_{t}{\vartheta}\right|^{2}{\rm d}x{\rm d}t
 \\\nonumber&\hspace{0.2cm}+\frac{1}{2}\alpha_{1}\alpha_{3}s\lambda\int_{Q}\varphi
 \left|\nabla{\vartheta}\right|^{2}{\rm d}x{\rm d}t
 +2\alpha_{3}s\lambda^{2}\int_{Q}\varphi
 \left|a_{5}\left(\nabla\psi\cdot\nabla{\vartheta}\right)+a_{6}\left(\nabla\psi\cdot\nabla{w}\right)\right|^{2}{\rm d}x{\rm d}t
 \\\nonumber&\hspace{0.3cm}+4 s\lambda^{2}\int_{Q}\varphi\left|{a_{7}}
 \left(\partial_{t}\psi\right)
 \left(\partial_{t}{w}\right)-a_{1}\left(\nabla\psi\cdot\nabla{w}\right)\right|^{2}{\rm d}x{\rm d}t
 \\\nonumber&\hspace{0.2cm}+\frac{1}{4}\sigma_{1}s\lambda\int_{Q}\varphi\left|\partial_{t}{w}\right|^{2}{\rm d}x{\rm d}t
 +\frac{1}{2}\sigma_{0}\sigma_{1}s\lambda\int_{Q}\varphi\left|\nabla{w}\right|^{2}{\rm d}x{\rm d}t \\\nonumber&\hspace{0.2cm}+16m^2\sigma^{2}_{1}s^{3}\lambda^{4}\int_{Q}\varphi^{3}
  {w}^{2}{\rm d}x{\rm d}t
 +\frac{1}{2}\left\|\mathcal{L}_{\ast}\left({\vartheta},{w}\right)\right\|_{L^2(Q)}^2
 +\left\|\mathcal{L}_{2\ast}\left({w}\right)\right\|_{L^2(Q)}^2
 \\\nonumber&\hspace{0.2cm}+\alpha_{3}\left\|L_{1}\left({\vartheta},{w}\right)\right\|^{2}_{L^2(Q)}
 \\\nonumber&\hspace{0.2cm}-2\alpha_{3}s\lambda\int_{\left(\partial\Omega\setminus\Gamma_0\right)\times(0,T)} \varphi\left(\nabla\psi\cdot\nu\right)\left(a_{5}
 \frac{\partial {\vartheta}}{\partial\nu}+a_{6}
 \frac{\partial{w}}{\partial \nu}\right)^2{\rm d}\sigma{\rm d}t
 \\\nonumber&\hspace{0.2cm}-2s\lambda\int_{\left(\partial\Omega\setminus\Gamma_0\right)\times(0,T)}
 \varphi(\nabla\psi\cdot\nu)a_{1}^2\left|\frac{\partial{w}}{\partial\nu}\right|^2{\rm d}\sigma{\rm d}t
 \\\nonumber&\hspace{0.2cm}-2\sqrt{n}M_0M_1\int_{\left(\partial\Omega\setminus\Gamma_0\right)\times(0,T)}  \left|\frac{\partial{w}}
 {\partial\nu}\right|^2{\rm d}\sigma{\rm d}t
 \\\nonumber&\leq C_{53}\int_{Q}\left(f^{2}+g^{2}\right){\rm e}^{2s\varphi}{\rm d}x{\rm d}t
 \\\nonumber&\hspace{0.2cm}+2\alpha_{3}s\lambda\int_{\Gamma_0\times(0,T)} \varphi\left(\nabla\psi\cdot\nu\right)\left(a_{5}
 \frac{\partial {\vartheta}}{\partial\nu}+a_{6}
 \frac{\partial{w}}{\partial \nu}\right)^2{\rm d}\sigma{\rm d}t
 \\\nonumber&\hspace{0.2cm}+2s\lambda\int_{\Gamma_0\times(0,T)}\varphi(\nabla\psi\cdot\nu)
 a_{1}^2\left|\frac{\partial{w}}{\partial\nu}\right|^2{\rm d}\sigma{\rm d}t
 +2\sqrt{n}M_0M_1\int_{\Gamma_0\times(0,T)}\left|\frac{\partial{w}}
 {\partial\nu}\right|^2{\rm d}\sigma{\rm d}t
 \\\nonumber&\hspace{0.2cm}+C_{53} s\lambda\int_{Q}\varphi\left|a_{5}\left(\nabla\psi\cdot\nabla{\vartheta}\right)
 +a_{6}\left(\nabla\psi\cdot\nabla{w}\right)\right|^{2}{\rm d}x{\rm d}t
 \\\nonumber&\hspace{0.2cm} +C_{53} s\lambda\int_{Q}\varphi\left|{a_{7}}
 \left(\partial_{t}\psi\right)\left(\partial_{t}{{w}}\right)
 -a_{1}\left(\nabla\psi\cdot\nabla{w}\right)\right|^{2}{\rm d}x{\rm d}t
 \\\nonumber&\hspace{0.2cm}+C_{53} s\lambda^{\frac{1}{2}}\int_{Q}\varphi\left(\left|\nabla{\vartheta}\right|^{2}
 +\left|\nabla{w}\right|^{2}+\left|\partial_{t}{w}\right|^{2}\right){\rm d}x{\rm d}t
 +C_{53} s^{3}\lambda^{\frac{7}{2}}\int_{Q}\varphi^{3}\left({\vartheta}^{2}+{w}^{2}\right){\rm d}x{\rm d}t
 \\&\hspace{0.2cm}+C_{53}\int_{Q}\frac{1}{s\lambda\varphi}\left|\partial_{t}{{\vartheta}}\right|^{2}
 {\rm d}x{\rm d}t
 +C_{53}s^{2}\lambda^{6}\int_{Q}\varphi^{2}\left({\vartheta}^{2}+{w}^{2}\right){\rm d}x{\rm d}t,
 \end{align*}
 for all $s>s_1$.
 By (\ref{epsilon}), there exists a constant $\epsilon>0$ such that $-(\nabla\psi\cdot\nu)\geq \epsilon$ on $\overline{\partial\Omega\setminus\Gamma_0}$. Therefore, there exists a constant $\lambda_{0}>1$ such that for any $\lambda\geq\lambda_{0}$, the fourteenth term in the left hand side can be absorbed by the thirteenth term in the left hand side, and the terms from the fifth to the ninth in the right hand side can be absorbed by the left hand side.
 Moreover, for any fixed $\lambda\geq\lambda_{0}$, there exists a constant
 $s_{0}=s_{0}(\lambda)>s_1$ such that for any $s\geq s_0(\lambda)$, the tenth term in the right hand side can be absorbed by the left hand side. We further note that the twelfth term in the left hand side is positive.  Thus there exists a constant $K_0=K_0(s_{0}$,
     $\lambda_{0}$, $\beta$, $\Omega$, $T$, $m$, $M$, $M_0$, $M_1$, $M_{2}$, $\sigma_0$, $\sigma_1$, $\epsilon)>0$ such that
 \begin{align*}
 &\nonumber\int_{Q}\left(s^{3}\lambda^{4}\varphi^{3}{\vartheta}^{2}
 +\frac{1}{s\varphi}\left|\partial_{t}{\vartheta}\right|^{2}
 +s\lambda\varphi\left|\nabla{\vartheta}\right|^{2}+s^{3}\lambda^{4}\varphi^{3}{w}^{2}
 +s\lambda\varphi\left|\nabla_{x,t}{w}\right|^{2}\right){\rm d}x{\rm d}t
 \\\nonumber&\hspace{0.2cm}+s\lambda^{2}\int_{Q}\varphi
 \left|a_{5}\left(\nabla\psi\cdot\nabla{\vartheta}\right)+a_{6}\left(\nabla\psi\cdot\nabla{w}\right)\right|^{2}{\rm d}x{\rm d}t
 \\\nonumber&\hspace{0.3cm}+ s\lambda^{2}\int_{Q}\varphi\left|{a_{7}}
 \left(\partial_{t}\psi\right)
 \left(\partial_{t}{w}\right)-a_{1}\left(\nabla\psi\cdot\nabla{w}\right)\right|^{2}{\rm d}x{\rm d}t
 +\left\|\mathcal{L}_{\ast}\left({\vartheta},{w}\right)\right\|_{L^2(Q)}^2
 \\\nonumber&\hspace{0.3cm}+\left\|\mathcal{L}_{2\ast}\left({w}\right)\right\|_{L^2(Q)}^2+\left\|L_{1}\left({\vartheta},{w}\right)\right\|^{2}_{L^2(Q)}
 +\epsilon s\lambda\int_{\left(\partial\Omega\setminus\Gamma_0\right)\times(0,T)}
 \varphi \left|\frac{\partial{w}}{\partial\nu}\right|^2{\rm d}\sigma{\rm d}t
 \\\nonumber&\leq K_0\int_{Q}\left(f^{2}+g^{2}\right){\rm e}^{2s\varphi}{\rm d}x{\rm d}t
 +K_0s\lambda\int_{\Gamma_0\times(0,T)}
 \varphi \left(\left|\frac{\partial{\vartheta}}{\partial\nu}\right|^2+\left|\frac{\partial{w}}{\partial\nu}\right|^2\right){\rm d}\sigma{\rm d}t,
 \end{align*}
 for all $\lambda\geq\lambda_0$ and $s\geq s_0(\lambda)$.
 Noting $\left({\vartheta},{w}\right)={\rm e}^{s\varphi}\left(\Theta,p\right)$, we complete the proof of Theorem \ref{carleman2}.

  \section{Proof of Theorem \ref{carleman}}\label{carproof1}
  \setcounter{equation}{0}

 We prove Theorem \ref{carleman} by applying Theorem \ref{carleman2}
 and the argument in \cite{Iman2002}.
 By a usual density argument, it is sufficient to prove (\ref{k2}) for $(\Theta,p)\in C^{\infty}(Q)\times C^{\infty}(Q)$.
 Let $G_{\varepsilon}=\{x\in\omega\big|\mbox{dist}(x, \partial\omega\cap\Omega)\leq\varepsilon\}$ and $\omega_{\varepsilon}=\omega\setminus G_{\varepsilon}$. Let $Q_{\omega_{\varepsilon}}=\omega_{\varepsilon}\times(0,T)$. By (\ref{omega}),
 there exists $\varepsilon>0$ such that
  \begin{equation}\label{omegavarepsilon}
  \overline{\partial\Omega\setminus\partial\omega_{4\varepsilon}}\subset\big{\{}x\in\partial\Omega\big{|} \left((x-x_0)\cdot\nu(x)\right)<0\big{\}}.
  \end{equation}
 We take $\Gamma_0=\partial\omega_{4\varepsilon}\cap\partial\Omega$. By (\ref{omegavarepsilon}), $\Gamma_0$ satisfies (\ref{epsilon}).
 In order to apply Theorem \ref{carleman2}, we take a cut-off function $\chi_1(x)\in C^{\infty}(\overline{\Omega})$ such that
 $0\leq\chi_1(x)\leq 1$ for all  $x\in \overline{\Omega}$,
 $\chi_1(x)= 1$ for all  $x\in \Omega\setminus\omega_{2\varepsilon}$, and $\chi_1(x)= 0$ for all  $x\in\omega_{3\varepsilon}$.
 Let $u(x,t)=\chi_1(x)p(x,t)$ and $v(x,t)=\chi_1(x)\Theta(x,t)$. By (\ref{fg}), (\ref{Db2}) and (\ref{Tb2}), we have
  \begin{equation*}
  \left\{\begin{aligned}
  \partial_{t}^{2}u-a_{1}\Delta u-a_{2}\Delta v=&\chi_1 f -a_{1} \left(\Delta\chi_1\right)p-2a_{1}\left(\nabla\chi_1\cdot\nabla p\right)
  \\&-a_{2} \left(\Delta\chi_1\right)\Theta-2a_{2}\left(\nabla\chi_1\cdot\nabla \Theta\right), \\
 \partial_{t} v-a_{3}\Delta v-a_{4}\partial_{t}^{2}u
 =&\chi_1 g-a_{3} \left(\Delta\chi_1\right)\Theta-2a_{3}\left(\nabla\chi_1\cdot\nabla \Theta\right), \ \ \mbox{in} \ Q,
 \end{aligned}\right.
  \end{equation*}
 ${v}(x,t)={u}(x,t)=0$ for $(x,t)\in\Sigma$,
 and ${v}(x,0)={v}(x,T)=\partial_t^j{{u}}(x,0)=\partial_t^j{{u}}(x,T)=0$ for $x\in\overline{\Omega}$ and $j=0,1$.
 Furthermore, by the definition of $\chi_1$, we have $\nabla v=\nabla u=0$ in some neighborhood of $\Gamma_0\times(0, T)$. By Theorem \ref{carleman2},
 there exists a constant $\eta(\beta)>0$ such that for any
 $T\in\left(0,\frac{2({\mathcal{D}}+\eta)}{\sqrt{\beta}}\right)$, there exists a constant $\lambda_{0}>0$ such that for all $\lambda>\lambda_{0}$, there exist constants $s_{0}(\lambda)>0$ and $K=K(s_{0}$,
     $\lambda_{0}$, $\beta$, $\Omega$, $T$, $m$, $M$, $M_0$, $M_1$, $M_{2}$, $\sigma_0$, $\sigma_1)>0$  such that
 \begin{align*}
 &\nonumber\int_{Q}\left(s^{3}\lambda^{4}\varphi^{3}{v}^{2}+\frac{1}{s\varphi}|\partial_{t}{v}|^{2}
 +s\lambda\varphi|\nabla{v}|^{2}+s^{3}\lambda^{4}\varphi^{3}{u}^{2}+s\lambda\varphi|\nabla_{x,t}{u}|^{2}\right){\rm e}^{2s\varphi}{\rm d}x{\rm d}t
 \\&\leq K\int_{Q}\left(f^2+g^2\right){\rm e}^{2s\varphi}{\rm d}x{\rm d}t
 +K\int_{\left(\omega_{2\varepsilon}\setminus\omega_{3\varepsilon}\right)\times(0,T)}\big\{\left|\left(\Delta\chi_1\right)p\right|^2
 +\left|\left(\nabla\chi_1\cdot\nabla p\right)\right|^2
 \\&\hspace{1cm}+\left|\left(\Delta\chi_1\right)\Theta\right|^2+\left|\left(\nabla\chi_1\cdot\nabla \Theta\right)\right|^2\big\}{\rm e}^{2s\varphi}{\rm d}x{\rm d}t,\quad \mbox{for all}\quad s\geq s_{0}.
 \end{align*}
 Therefore we have
 \begin{align}\label{carleman3}
 \nonumber&\int_{Q}\left\{s^{3}\lambda^{4}\varphi^{3}\left({\Theta}^{2}+{p}^{2}\right)+\frac{1}{s\varphi}|\partial_{t}{\Theta}|^{2}
 +s\lambda\varphi\left(|\nabla{\Theta}|^{2}+|\nabla_{x,t}{p}|^{2}\right)\right\}{\rm e}^{2s\varphi}{\rm d}x{\rm d}t
 \\&\nonumber=\int_{Q\setminus Q_{\omega_{2\varepsilon}}}
 \bigg\{s^{3}\lambda^{4}\varphi^{3}\left({v}^{2}+{u}^{2}\right)+\frac{1}{s\varphi}|\partial_{t}{v}|^{2}
 +s\lambda\varphi\left(|\nabla{v}|^{2}+|\nabla_{x,t}{u}|^{2}\right)\bigg\}{\rm e}^{2s\varphi}{\rm d}x{\rm d}t
 \\\nonumber&\hspace{0.1cm}+\int_{Q_{\omega_{2\varepsilon}}}\left\{s^{3}\lambda^{4}\varphi^{3}\left({\Theta}^{2}+{p}^{2}\right)
 +\frac{1}{s\varphi}|\partial_{t}{\Theta}|^{2}
 +s\lambda\varphi\left(|\nabla{\Theta}|^{2}+|\nabla_{x,t}{p}|^{2}\right)\right\}{\rm e}^{2s\varphi}{\rm d}x{\rm d}t
 \\\nonumber&\leq K\int_{Q}\left(f^2+g^2\right){\rm e}^{2s\varphi}{\rm d}x{\rm d}t
 +K\int_{\left(\omega_{2\varepsilon}\setminus\omega_{3\varepsilon}\right)\times(0,T)}\big\{\left|\left(\Delta\chi_1\right)p\right|^2
 +\left|\left(\nabla\chi_1\cdot\nabla p\right)\right|^2
 \\\nonumber&\hspace{1cm}+\left|\left(\Delta\chi_1\right)\Theta\right|^2+\left|\left(\nabla\chi_1\cdot\nabla \Theta\right)\right|^2\big\}{\rm e}^{2s\varphi}{\rm d}x{\rm d}t
 \\\nonumber&\hspace{0.05cm}+\int_{Q_{\omega_{2\varepsilon}}}\left\{s^{3}\lambda^{4}\varphi^{3}\left({\Theta}^{2}+{p}^{2}\right)
 +\frac{1}{s\varphi}|\partial_{t}{\Theta}|^{2}
 +s\lambda\varphi\left(|\nabla{\Theta}|^{2}+|\nabla_{x,t}{p}|^{2}\right)\right\}{\rm e}^{2s\varphi}{\rm d}x{\rm d}t
 \\\nonumber&\leq K\int_{Q}\left(f^2+g^2\right){\rm e}^{2s\varphi}{\rm d}x{\rm d}t
 +C_{54}\int_{Q_{\omega_{2\varepsilon}}}\bigg\{s^{3}\lambda^{4}\varphi^{3}\left({\Theta}^{2}+{p}^{2}\right)+\frac{1}{s\lambda\varphi}|\partial_{t}{\Theta}|^{2}
 \\&\hspace{1cm}+s\lambda\varphi\left(|\nabla{\Theta}|^{2}+|\nabla_{x,t}{p}|^{2}\right)\bigg\}{\rm e}^{2s\varphi}{\rm d}x{\rm d}t,\ \mbox{for all} \ s\geq s_{0}.
 \end{align}

 Next we estimate $s\lambda\int_{Q_{\omega_{2\varepsilon}}}
\varphi\left(|\nabla{\Theta}|^{2}+|\nabla p|^2\right){\rm e}^{2s\varphi}{\rm d}x{\rm d}t$  in (\ref{carleman3}). We take a cut-off function  $\chi_2(x)\in C^{\infty}(\overline{\Omega})$ such that $0\leq \chi_2(x)\leq 1$ for all  $x\in \overline{\Omega}$, $\chi_2(x)= 0$ for all $x\in\Omega\setminus\omega_{\varepsilon}$, and $\chi_2(x)= 1$ for all $x\in\overline{\omega_{2\varepsilon}}$.
 We multiply the second equation in (\ref{fg}) by $s\lambda\varphi\chi_2\Theta{\rm e}^{2s\varphi}$ and integrate it over $Q$. Integrating by parts and using (\ref{Db2}) and (\ref{Tb2}), we have
 \begin{align}\label{5}
 \nonumber&s\lambda\int_{Q_{\omega_{\varepsilon}}}\varphi \chi_2g\Theta{\rm e}^{2s\varphi}{\rm d}x{\rm d}t=
 \frac{1}{2}s\lambda\int_{Q_{\omega_{\varepsilon}}} \varphi\chi_2\left\{\partial_t\left(\Theta^2\right)\right\}{\rm e}^{2s\varphi}{\rm d}x{\rm d}t
 \\\nonumber&\hspace{0.5cm}-s\lambda\int_{Q_{\omega_{\varepsilon}}} \varphi\chi_2 a_3\left(\bigtriangleup\Theta\right)\Theta{\rm e}^{2s\varphi}{\rm d}x{\rm d}t
 -s\lambda\int_{Q_{\omega_{\varepsilon}}} \varphi\chi_2 a_4\left(\partial_t^2p\right)\Theta{\rm e}^{2s\varphi}{\rm d}x{\rm d}t
 \\\nonumber&=-\frac{1}{2}s\lambda^2\int_{Q_{\omega_{\varepsilon}}} \varphi\chi_2\left(\partial_t\psi\right)\Theta^2{\rm e}^{2s\varphi}{\rm d}x{\rm d}t
-s^2\lambda^2\int_{Q_{\omega_{\varepsilon}}} \varphi^2\left(\partial_t\psi\right)
\chi_2\Theta^2{\rm e}^{2s\varphi}{\rm d}x{\rm d}t
 \\\nonumber&\hspace{0.5cm}+s\lambda^2\int_{Q_{\omega_{\varepsilon}}} \varphi\chi_2 a_3\left(\nabla \psi\cdot\nabla\Theta\right)\Theta
 {\rm e}^{2s\varphi}{\rm d}x{\rm d}t
 \\\nonumber&\hspace{0.5cm}+s\lambda\int_{Q_{\omega_{\varepsilon}}} \varphi\left(\nabla \left(\chi_2 a_3\right)\cdot\nabla\Theta\right)
 \Theta{\rm e}^{2s\varphi}{\rm d}x{\rm d}t
 \\\nonumber&\hspace{0.5cm}+s\lambda\int_{Q_{\omega_{\varepsilon}}} \varphi \chi_2 a_3\left|\nabla\Theta\right|^2{\rm e}^{2s\varphi}{\rm d}x{\rm d}t
 +2s^2\lambda^2\int_{Q_{\omega_{\varepsilon}}} \varphi^2 \chi_2 a_3 \left(\nabla \psi\cdot\nabla\Theta\right)\Theta {\rm e}^{2s\varphi}{\rm d}x{\rm d}t
 \\\nonumber&\hspace{0.5cm}+s\lambda^2\int_{Q_{\omega_{\varepsilon}}} \varphi \left(\partial_t\psi\right)\chi_2a_4\left(\partial_tp\right)\Theta{\rm e}^{2s\varphi}{\rm d}x{\rm d}t
 +s\lambda\int_{Q_{\omega_{\varepsilon}}} \varphi\chi_2a_4\left(\partial_tp\right)
 \left(\partial_t\Theta\right){\rm e}^{2s\varphi}{\rm d}x{\rm d}t
 \\&\hspace{0.5cm}+2s^2\lambda^2\int_{Q_{\omega_{\varepsilon}}} \varphi^2 \left(\partial_t\psi\right)\chi_2a_4\left(\partial_tp\right)\Theta{\rm e}^{2s\varphi}{\rm d}x{\rm d}t\triangleq \sum_{k=1}^{10}\tilde{J}_{k},\quad \mbox{for all}\quad s\geq s_{0}.
 \end{align}
 We recall $s>s_0\geq 1$ and $\lambda>1$. Thus we have
 $$
  \tilde{J}_1+\tilde{J}_2\geq -C_{55}s^2\lambda^2\int_{Q_{\omega_{\varepsilon}}} \varphi^2 \Theta^2{\rm e}^{2s\varphi}{\rm d}x{\rm d}t,\ \mbox{for all}\ s\geq s_{0}.
  $$
  By the Cauchy-Schwarz inequality, we have
  \begin{align*}
  \tilde{J}_3&\geq -C_{56}\int_{Q_{\omega_{\varepsilon}}}\chi_2
  \left(\frac{1}{2}s^{\frac{1}{2}}\lambda^{\frac{1}{2}}\sigma_1^{\frac{1}{2}}\varphi^{\frac{1}{2}}|\nabla\Theta|\right)
  \left(2s^{\frac{1}{2}}\lambda^{\frac{3}{2}}\sigma_1^{-\frac{1}{2}}\varphi^{\frac{1}{2}}|\Theta|\right){\rm e}^{2s\varphi}{\rm d}x{\rm d}t
  \\&\geq -\frac{1}{4}s\lambda\sigma_1\int_{Q_{\omega_{\varepsilon}}}\varphi\chi_2 |\nabla\Theta|^2{\rm e}^{2s\varphi}{\rm d}x{\rm d}t
  -C_{57} s\lambda^3\int_{Q_{\omega_{\varepsilon}}} \varphi\chi_2\Theta^2{\rm e}^{2s\varphi}{\rm d}x{\rm d}t,
  \end{align*}
  for all $s\geq s_{0}$. Integrating by parts, we have
  \begin{align*}
  \tilde{J}_4&=\frac{1}{2}s\lambda\int_{Q_{\omega_{\varepsilon}}} \varphi\left(\nabla \left(\chi_2 a_3\right)\cdot\nabla\left(\Theta^2\right)\right)
 {\rm e}^{2s\varphi}{\rm d}x{\rm d}t
 \\&=-\frac{1}{2} s\lambda^2\int_{Q_{\omega_{\varepsilon}}} \varphi\left(\nabla\left(\chi_2 a_3\right)\cdot\nabla\psi\right) \Theta^2{\rm e}^{2s\varphi}{\rm d}x{\rm d}t
  \\
  &\hspace{0.5cm}-\frac{1}{2} s\lambda\int_{Q_{\omega_{\varepsilon}}} \varphi\left\{\triangle\left(\chi_2 a_3\right)\right\} \Theta^2{\rm e}^{2s\varphi}{\rm d}x{\rm d}t
  \\
  &\hspace{0.5cm}- s^2\lambda^2\int_{Q_{\omega_{\varepsilon}}} \varphi^2\left(\nabla\left(\chi_2 a_3\right)\cdot\nabla\psi\right) \Theta^2{\rm e}^{2s\varphi}{\rm d}x{\rm d}t
  \\
  &\geq -C_{58}s^2\lambda^2\int_{Q_{\omega_{\varepsilon}}} \varphi^2 \Theta^2{\rm e}^{2s\varphi}{\rm d}x{\rm d}t.
  \end{align*}
  $$
  \tilde{J}_5\geq s\lambda\sigma_1\int_{Q_{\omega_{\varepsilon}}} \varphi\chi_2 \left|\nabla\Theta\right|^2{\rm e}^{2s\varphi}{\rm d}x{\rm d}t,\ \mbox{for all}\ s\geq s_{0}.
  $$
  By the Cauchy-Schwarz inequality, we have
  \begin{align*}
  \tilde{J}_6&\geq -C_{59}\int_{Q_{\omega_{\varepsilon}}} \chi_2 \left(\frac{1}{2}s^{\frac{1}{2}}\lambda^{\frac{1}{2}}\sigma_1^{\frac{1}{2}}\varphi^{\frac{1}{2}}|\nabla\Theta|\right)
  \left(2s^{\frac{3}{2}}\lambda^{\frac{3}{2}}\sigma_1^{-\frac{1}{2}}\varphi^{\frac{3}{2}}|\Theta|\right) {\rm e}^{2s\varphi}{\rm d}x{\rm d}t\\
  &\geq-\frac{1}{4}s\lambda\sigma_1\int_{Q_{\omega_{\varepsilon}}} \varphi\chi_2 \left|\nabla\Theta\right|^2{\rm e}^{2s\varphi}{\rm d}x{\rm d}t-C_{60}s^3\lambda^{3}\int_{Q_{\omega_{\varepsilon}}} \varphi^3\chi_2 \Theta^2{\rm e}^{2s\varphi}{\rm d}x{\rm d}t,
  \end{align*}
  $$
  \tilde{J}_7+\tilde{J}_{9}\geq -C_{61}s^2\lambda^2\int_{Q_{\omega_{\varepsilon}}} \varphi^2\chi_2
   \left(\left|\partial_tp\right|^2+\Theta^2\right){\rm e}^{2s\varphi}{\rm d}x{\rm d}t,
  $$
  \begin{align*}
  &\tilde{J}_8\geq -C_{62}\int_{Q_{\omega_{\varepsilon}}}\chi_2
   \left(s^{\frac{3}{2}}\lambda^{\frac{3}{2}}\varphi^{\frac{3}{2}}\left|\partial_tp\right|\right)
   \left(\frac{1}{s^{\frac{1}{2}}\lambda^{\frac{1}{2}}\varphi^{\frac{1}{2}}}\left|\partial_t\Theta\right|\right){\rm e}^{2s\varphi}{\rm d}x{\rm d}t
   \\&\geq -C_{63} s^{3}\lambda^{3}\int_{Q_{\omega_{\varepsilon}}}\varphi^{3}\chi_2
   \left|\partial_tp\right|^2{\rm e}^{2s\varphi}{\rm d}x{\rm d}t
  -\int_{Q_{\omega_{\varepsilon}}}
   \frac{\chi_2 }{s\lambda\varphi}\left|\partial_t\Theta\right|^2{\rm e}^{2s\varphi}{\rm d}x{\rm d}t,
  \end{align*}
 \begin{align*}
  s\lambda\int_{Q_{\omega_{\varepsilon}}}&\varphi \chi_2g\Theta{\rm e}^{2s\varphi}{\rm d}x{\rm d}t
  \leq C_{64}\int_{Q_{\omega_{\varepsilon}}} \chi_2\left(\frac{1}{s^{\frac{1}{2}}\lambda^{\frac{1}{2}}\varphi^{\frac{1}{2}}}|g|\right)
  \left(s^{\frac{3}{2}}\lambda^{\frac{3}{2}}\varphi^{\frac{3}{2}}|\Theta|\right){\rm e}^{2s\varphi}{\rm d}x{\rm d}t
  \\& \leq \int_{Q_{\omega_{\varepsilon}}} \frac{\chi_2}{s\lambda\varphi}g^2{\rm e}^{2s\varphi}{\rm d}x{\rm d}t
  +C_{65} s^{3}\lambda^{3}\int_{Q_{\omega_{\varepsilon}}} \varphi^{3}\chi_2
 \Theta^2{\rm e}^{2s\varphi}{\rm d}x{\rm d}t,
  \end{align*}
 for all $s\geq s_{0}$.
 Therefore, we obtain
 \begin{align}\label{2}
  &\frac{1}{2}s\lambda\sigma_1\int_{Q_{\omega_{\varepsilon}}} \varphi\chi_2 \left|\nabla\Theta\right|^2{\rm e}^{2s\varphi}{\rm d}x{\rm d}t\nonumber
  \\\nonumber
  &\leq\int_{Q_{\omega_{\varepsilon}}}
   \frac{\chi_2 }{s\lambda\varphi}\left|\partial_t\Theta\right|^2{\rm e}^{2s\varphi}{\rm d}x{\rm d}t + \int_{Q_{\omega_{\varepsilon}}} \frac{\chi_2}{s\lambda\varphi}g^2{\rm e}^{2s\varphi}{\rm d}x{\rm d}t
   \\&\hspace{0.5cm}
  +C_{66}s^{3}\lambda^{3}\int_{Q_{\omega_{\varepsilon}}}\varphi^{3}
 \left(\left|\partial_tp\right|^2+\Theta^2\right){\rm e}^{2s\varphi}{\rm d}x{\rm d}t,\ \mbox{for all}\ s\geq s_{0}.
  \end{align}
 By (\ref{fg}), we can get (\ref{a_78}) where $a_7$ and $a_8$ are given by (\ref{a78}).  Similarly to (\ref{5})--(\ref{2}), we multiply (\ref{a_78}) by $s\lambda\varphi\chi_2p{\rm e}^{2s\varphi}$ and integrate it over $Q$, so that, by the integration by parts and the Cauchy-Schwarz inequality, we have
 \begin{align*}
 s&\lambda\int_{Q_{\omega_{\varepsilon}}}\varphi \chi_2\left(f-a_8g\right)p{\rm e}^{2s\varphi}{\rm d}x{\rm d}t=
 s\lambda\int_{Q_{\omega_{\varepsilon}}} \varphi\chi_2a_7\left(\partial_t^2p\right)p{\rm e}^{2s\varphi}{\rm d}x{\rm d}t
 \\\nonumber
  &\hspace{0.5cm}-s\lambda\int_{Q_{\omega_{\varepsilon}}} \varphi\chi_2a_1\left(\triangle p\right)p{\rm e}^{2s\varphi}{\rm d}x{\rm d}t
 -s\lambda\int_{Q_{\omega_{\varepsilon}}} \varphi\chi_2a_8\left(\partial_t \Theta\right)p{\rm e}^{2s\varphi}{\rm d}x{\rm d}t
 \\
  &\geq -C_{67}s\lambda\int_{Q_{\omega_{\varepsilon}}}\varphi\chi_2\left|\partial_t p\right|^2{\rm e}^{2s\varphi}{\rm d}x{\rm d}t
  -C_{67}s^3\lambda^3\int_{Q_{\omega_{\varepsilon}}}\varphi^3 p^2{\rm e}^{2s\varphi}{\rm d}x{\rm d}t
  \\\nonumber
  &\hspace{0.5cm}+\frac{1}{2}s\lambda\sigma_1\int_{Q_{\omega_{\varepsilon}}}\varphi\chi_2\left|\nabla p\right|^2{\rm e}^{2s\varphi}{\rm d}x{\rm d}t-C_{67}s\lambda\int_{Q_{\omega_{\varepsilon}}}\varphi\chi_2\Theta^2{\rm e}^{2s\varphi}{\rm d}x{\rm d}t,
 \end{align*}
 for all $s\geq s_{0}$. Then we have
 \begin{align}\label{6}
 \nonumber\frac{1}{2}&s\lambda\sigma_1\int_{Q_{\omega_{\varepsilon}}}\varphi\chi_2
 \left|\nabla p\right|^2{\rm e}^{2s\varphi}{\rm d}x{\rm d}t
 \\\nonumber&\leq \int_{Q_{\omega_{\varepsilon}}}\frac{\chi_2}{s\lambda\varphi}\left(f^2+g^2\right){\rm e}^{2s\varphi}{\rm d}x{\rm d}t
  +C_{68}s^3\lambda^3\int_{Q_{\omega_{\varepsilon}}}\varphi^3 p^2{\rm e}^{2s\varphi}{\rm d}x{\rm d}t\\
  &\hspace{0.5cm}+C_{68}s\lambda\int_{Q_{\omega_{\varepsilon}}}\varphi\chi_2\left|\partial_t p\right|^2{\rm e}^{2s\varphi}{\rm d}x{\rm d}t
 +C_{68}s\lambda\int_{Q_{\omega_{\varepsilon}}}\varphi\chi_2\Theta^2{\rm e}^{2s\varphi}{\rm d}x{\rm d}t,
 \end{align}
 for all $s\geq s_{0}$.
 Adding (\ref{2}) and (\ref{6}) and noting $\omega_{2\varepsilon}\subseteq\omega_{\varepsilon}\subseteq\omega$ and $\chi_2(x)=1$ for all
 $x\in\overline{\omega_{2\varepsilon}}$, we obtain
 \begin{align}\label{7}
  &s\lambda\int_{Q_{\omega_{2\varepsilon}}} \varphi \left(\left|\nabla\Theta\right|^2+\left|\nabla p\right|^2\right){\rm e}^{2s\varphi}{\rm d}x{\rm d}t\leq
  s\lambda\int_{Q_{\omega_{\varepsilon}}} \varphi\chi_2 \left(\left|\nabla\Theta\right|^2+\left|\nabla p\right|^2\right)
  {\rm e}^{2s\varphi}{\rm d}x{\rm d}t\nonumber
  \\\nonumber
  &\leq\frac{2}{\sigma_1}\int_{Q_{\omega_{\varepsilon}}}
   \frac{\chi_2 }{s\lambda\varphi}\left|\partial_t\Theta\right|^2{\rm e}^{2s\varphi}{\rm d}x{\rm d}t
   +\frac{4}{\sigma_1} \int_{Q_{\omega_{\varepsilon}}} \frac{\chi_2}{s\lambda\varphi}
   \left(f^2+g^2\right){\rm e}^{2s\varphi}{\rm d}x{\rm d}t
   \\&\hspace{0.5cm}
  +C_{69}s^{3}\lambda^{3}\int_{Q_{\omega_{\varepsilon}}}\varphi^{3}
 \left(\left|\partial_tp\right|^2+p^2+\Theta^2\right){\rm e}^{2s\varphi}{\rm d}x{\rm d}t,\ \mbox{for all}\ s\geq s_{0}.
  \end{align}
  Substituting (\ref{7}) into (\ref{carleman3}), we obtain
  \begin{align}\label{10}
 \nonumber&\int_{Q}\left\{s^{3}\lambda^{4}\varphi^{3}\left({\Theta}^{2}+{p}^{2}\right)+\frac{1}{s\varphi}|\partial_{t}{\Theta}|^{2}
 +s\lambda\varphi\left(|\nabla{\Theta}|^{2}+|\nabla_{x,t}{p}|^{2}\right)\right\}{\rm e}^{2s\varphi}{\rm d}x{\rm d}t
 \\\nonumber&\leq C_{70}\int_{Q}\left(f^2+g^2\right){\rm e}^{2s\varphi}{\rm d}x{\rm d}t
 \\&\hspace{0.5cm}+C_{70}\int_{Q_{\omega_{\varepsilon}}}\bigg\{s^{3}\lambda^{4}\varphi^{3}\left({\Theta}^{2}+{p}^{2}\right)
 +\frac{1}{s\varphi}|\partial_{t}{\Theta}|^{2}
+s^3\lambda^3\varphi^3|\partial_t{p}|^{2}\bigg\}{\rm e}^{2s\varphi}{\rm d}x{\rm d}t,
 \end{align}
 for all $s\geq s_{0}$. Noting $\omega_\varepsilon\subseteq\omega$, we obtain (\ref{k}).
 The proof of Theorem \ref{carleman} is complete.

  \section*{Acknowledgements}
   This paper has been partially done during the stay of Li S  at the Institut de Math$\acute{\rm e}$matiques de Marseille, Aix-Marseille University, France, in July 2019, and the stay was supported by Aix-Marseille University.


\begin{thebibliography}{100}


 \bibitem{Akhouayri:2016}
 H. Akhouayri, M. Bergounioux, A. D. Silva, P. Elbau, A. Litman and L. Mindrinos, Quantitative thermoacoustic tomography with microwaves sources, J. Inverse Ill-Posed Probl. 25 (2016), 703-717.

 \bibitem{bell:2018}
 M. Bellassoued and M. Yamamoto, Carleman Estimates and Applications to Inverse Problems for Hyperbolic Systems, Springer, Tokyo, 2018.

  \bibitem{carl:1939}
 T. Carleman, Sur un probl$\grave{\rm e}$me d'unicit$\acute{\rm e}$ pour les syst¨¨mes d'$\acute{\rm e}$quations aux d$\acute{\rm e}$riv$\acute{\rm e}$es partielles $\grave{\rm a}$ deux variables ind$\acute{\rm e}$pendentes [On a problem of the uniqueness for systems of partial differential equations with two independent variables], Ark. Mat.Astr. Fys.  26B (1939), 1-9.

 \bibitem{chae:1996}
  D. Chae, O. Y. Imanuvilov and S. M. Kim,  Exact controllability for semilinear parabolic equations with Neumann boundary conditions, J. Dyn. Control Syst. 2 (1996), 449-483.

 \bibitem{Cox:2009}
 B. Cox and P. C. Beard,  Modeling photoacoustic propagation in tissue using k-space techniques,
 In: Wang L V, editor,
 Photoacoustic Imaging and Spectroscopy, chapter 3,
 pp. 25-34, CRC Press, 2009.

 \bibitem{part II}
  M. Cristofol, S. Li and Y. Shang,  Carleman estimates and inverse problems for the coupled quantitative thermoacoustic equations. Part II: inverse problems, preprint, 2020.


 \bibitem{egor:1986}
 Y. V. Egorov,  Linear Differential Equations of Principal Type, Consultants Bureau, New York, 1986.

 \bibitem{elle:2000}
 M. M. Eller and V. Isakov, Carleman estimates with two large parameters and applications, Contemp. Math. 268 (2000), 117-136.

 \bibitem{furs:2000}
  A. V. Fursikov and O. Y. Imanuvilov, Controllability of Evolution Equations, Lecture Notes Series vol 34, Seoul National University, Korea, 1996.

 \bibitem{horm:1963}
  L. H$\ddot{\rm o}$rmander, Linear Partial Differential Operators, Springer, Berlin, 1963.


 \bibitem{horm:1985}
 L. H$\ddot{\rm o}$rmander, The Analysis of Linear Partial Differential Operators I-IV, Springer, Berlin, 1985.


 \bibitem{Iman2002}
 O. Y. Imanuvilov, On Carleman estimates for hyperbolic equations, Asymptot. Anal. 32 (2002), 185-220.

 \bibitem{iman:2004}
 O. Y. Imanuvilov and M. Yamamoto, Carleman estimate for a stationary isotropic Lam$\acute{\rm e}$ system and the applications, Appl. Anal. 83 (2004), 243-270.

  \bibitem{isak:1986}
  V. Isakov,  A nonhyperbolic Cauchy problem for $\Box_b\Box_c$ and its applications to elasticity theory,
  Comm. Pure Appl. Math. 39 (1986), 747-767.

 \bibitem{isak:1990}
  V. Isakov, Inverse Source Problems, American Mathematical Society, Providence, RI, 1990.

  \bibitem{isak:1993}
  V. Isakov, Carleman type estimates in an anisotropic case and applications, J. Differential Equations 105 (1993), 217-238.

  \bibitem{isak:2004}
   V. Isakov, Carleman estimates and applications to inverse problems, Milan J. Math. 72 (2004), 249-271.

  \bibitem{isak:2006}
  V. Isakov, Inverse problems for partial differential equations, Springer-Verlag, Berlin, 2006.

  \bibitem{klib:1992}
   M. V. Klibanov, Inverse problems and Carleman estimates, Inverse Problems 8 (1992), 575-596.

  \bibitem{klib:2004}
  M. V. Klibanov and A. Timonov, Carleman estimates for coefficient inverse problems and numerical applications, VSP, Utrect, 2004.

  \bibitem{lavr:1986}
   M. M. Lavrent'ev, V. G. Romanov and S. P. Shishatskii, Ill-posed problems of mathematical physics and analysis, American Mathematical Society, Providence, RI, 1986.

 \bibitem{li:2015}
 S. Li, Carleman estimates for second-order hyperbolic systems in anisotropic cases and applications. Part I: Carleman estimates, Appl. Anal. 94(2015), 2261-2286.

 \bibitem{Patch:2007}
 S. K. Patch and O. Scherzer, Guest editors' introduction: Photo-and thermo-acoustic imaging, Inverse Problems 23 (2007), S1-S10.


  \bibitem{roma:2006}
   V. G. Romanov, Carleman estimates for second-order hyperbolic equations, Sib. Math. J. 47 (2006) 135-151.


 \bibitem{Stefanov:2009}
 P. Stefanov and G. Uhlmann, Thermoacoustic tomography with variable sound speed, Inverse Problems 25 (2009), 075011.


  \bibitem{tata:1996}
   D. Tataru,  Carleman estimates and unique continuation for solutions to boundary value problems, J. Math. Pures Appl. 75 (1996), 367-408.

  \bibitem{trev:1970}
  F. Tr$\grave{\rm e}$ves, Linear Partial Differential Equations, Gordon and Breach, New York, 1970.


 \bibitem{Yamamoto:2009}
 M. Yamamoto, Carleman estimates for parabolic equations and applications, Inverse Problems 25 (2009), 123013.

 \bibitem{yuan:2009}
  G. Yuan and M. Yamamoto, Lipshitz stability in the determination of the principal part of a parabolic equation, ESAIM Control Optim. Calc. Var. 15 (2009), 525-554.


\end{thebibliography}
 \end{document}